\documentclass[11pt, twoside]{amsart}
\calclayout
\usepackage{amsmath,amsthm, amscd, amssymb, amsfonts, mathrsfs,pb-diagram,pb-xy}
\usepackage[all]{xy}
\usepackage[dvips, dvipsnames, usenames]{color}
\usepackage{enumitem}
\usepackage[section]{placeins}
\usepackage{hyperref}

\usepackage[T1]{fontenc}

\usepackage{adjustbox}
\usepackage{microtype}

\allowdisplaybreaks

\newcommand{\ot}{\otimes}

\newcommand{\ub}{\mathbf u}

\newcommand{\un}[1]{\overline{#1}}

\newcommand{\ft}{\mathtt{f}}
\newcommand{\gt}{\mathtt{g}}
\newcommand{\kt}{\mathtt{k}}
\newcommand{\ut}{\mathtt{u}}
\newcommand{\vt}{\mathtt{v}}
\newcommand{\wt}{\mathtt{w}}
\newcommand{\xt}{\mathtt{x}}
\newcommand{\yt}{\mathtt{y}}

\newcommand{\Ht}{\mathtt{H}}

\newcommand{\paraA}{\mathbf a}
\newcommand{\paraB}{\mathbf b}
\newcommand{\paraC}{\mathbf c}
\newcommand{\paraD}{\mathbf d}

\usepackage{hhline}
\usepackage[active]{srcltx}

\newcommand{\ord}{\operatorname{ord}}

\newcommand{\co}{\operatorname{co}}
\newcommand{\supp}{\operatorname{supp}}
\newcommand{\ad}{\operatorname{ad}}
\newcommand{\bq}{\mathfrak{q}}
\newcommand{\ba}{\mathbf{a} }

\newcommand{\gr}{\operatorname{gr}}

\newcommand{\R}{{\mathcal R}}

\newcommand{\Zc}{{\mathcal Z}}
\newcommand{\Pc}{{\mathcal P}}

\newcommand{\Jc}{{\mathcal J}}

\newcommand{\Z}{{\mathbb Z}}
\newcommand{\N}{{\mathbb N}}
\newcommand{\I}{{\mathbb I}}
\newcommand{\Jb}{{\mathbb J}}

\newcommand{\G}{{\mathbb G}}

\newcommand{\toba}{\mathscr B}
\newcommand{\wtoba}{\widetilde{\mathscr B}}
\newcommand{\htoba}{\widehat{\mathscr B}}

\newcommand{\ydG}{{}^{\Bbbk G}_{\Bbbk G}\mathcal{YD}}
\newcommand{\ydg}{{}^{\Bbbk \Gamma}_{\Bbbk \Gamma}\mathcal{YD}}

\newcommand{\Ss}{\mathcal S}

\newcommand{\End}{\operatorname{End}}

\newcommand{\Hom}{\operatorname{Hom}}
\newcommand{\id}{\operatorname{id}}

\newcommand{\cL}{\mathcal{L}}
\newcommand{\bsl}{\boldsymbol{\lambda}}
\newcommand{\cG}{{\mathcal G}}
\newcommand{\att}{\mathtt{a}}

\newcommand\cA{\mathcal{A}}
\newcommand\mH{\mathcal{H}}
\newcommand{\cR}{\mathcal R}

\newcommand\Alg{\operatorname{Alg}}
\newcommand{\cE}{{\mathcal E}}

\newcommand\mP{\mathcal{P}}

\numberwithin{equation}{section}
\theoremstyle{plain}

\newtheorem{theorem}{Theorem}[section]
\newtheorem{lemma}[theorem]{Lemma}
\newtheorem{coro}[theorem]{Corollary}

\newtheorem{prop}[theorem]{Proposition}

\theoremstyle{definition}
\newtheorem{definition}[theorem]{Definition}
\newtheorem{exa}[theorem]{Example}

\newtheorem{remark}[theorem]{Remark}
\newtheorem{obs}[theorem]{Remark}

\newtheorem{question}{Question}
\newtheorem{step}{Step}
\newtheorem{stepi}{Step}

\newtheorem{stepv}{Step}

\def\pf{\begin{proof}}
\def\epf{\end{proof}}

\newcounter{tabla}\stepcounter{tabla}

\begin{document}

\renewcommand{\baselinestretch}{1.2}

\thispagestyle{empty}
\title[Pointed Hopf algebras over non-abelian groups]
{Pointed Hopf algebras over non-abelian groups with non-simple standard braidings}

\author[I. Angiono, S. Lentner, G. Sanmarco]
{Iv\'an Angiono, Simon Lentner, Guillermo Sanmarco}

\thanks{The work of I. A. was partially supported by CONICET and Secyt (UNC). The authors were partially supported by the Alexander von Humboldt Foundation through the Research Group Linkage Programme.}

\address{\noindent Facultad de Matem\'atica, Astronom\'{\i}a, F\'{\i}sica y Computaci\'on,
Universidad Nacional de C\'ordoba. CIEM -- CONICET. 
Medina Allende s/n (5000) Ciudad Univ., C\'ordoba,
Argentina}
\email{angiono@famaf.unc.edu.ar}

\address{\noindent Algebra and Number Theory. 
University of Hamburg,
Bundesstra{\ss}e 55,
20146 Hamburg
}
\email{simon.lentner@uni-hamburg.de}

\address{\noindent Department of Mathematics,
Iowa State University,
Ames, IA 50011,
USA}
\email{sanmarco@iastate.edu}

\makeatletter
\@namedef{subjclassname@2020}{%
\textup{2020} Mathematics Subject Classification}
\makeatother


\begin{abstract} 
We construct finite-dimensional Hopf algebras whose coradical is the group algebra of a central extension of an abelian group. They fall into families associated to a semisimple Lie algebra together with a Dynkin diagram automorphism. 
We show conversely that every finite-dimensional pointed Hopf algebra over a non-abelian group with non-simple infinitesimal braiding of rank at least 4 is of this form.

We follow the steps of the Lifting Method by Andruskiewitsch--Schneider. Our starting point is the classification of finite-dimensional Nichols algebras over non-abelian groups by Heckenberger--Vendramin, which consist of low rank exceptions and large rank families. We prove that the large rank families are cocycle twists of Nichols algebras constructed by the second author as foldings of Nichols algebras of Cartan type over abelian groups by outer automorphisms. This enables us to give uniform Lie-theoretic descriptions of the large rank families, prove generation in degree one and construct liftings. We also show that every lifting is a cocycle deformation of the corresponding coradically graded Hopf algebra using an explicit presentation by generators and relations of the Nichols algebra.

On the level of tensor categories, we construct families of graded extensions of the representation category of a quantum group by a group of diagram automorphism.

\end{abstract}
\maketitle

\tableofcontents

\section{Introduction}\label{subsec:Introduction}

\subsection{Background}

Groups and Lie algebras have in common that they both admit a tensor product of representations and a dual representation; in other words, their categories of representations are tensor categories. More generally, the category of representations of a Hopf algebra is also a tensor category. 
Prominent examples of Hopf algebras are the quantum groups $U_q(\mathfrak{g})$ by Drinfeld-Jimbo \cite{D-qg,Ji}, which are deformations of the enveloping algebra of a semisimple Lie algebra $\mathfrak{g}$ by a formal parameter $q$, and the small quantum groups $u_q(\mathfrak{g})$ by Lusztig \cite{Lu}, which are finite-dimensional non-semisimple quotients of $U_q(\mathfrak{g})$ for $q$ a root of unity. 
One of the initial motivations for quantum groups was their relation to the monodromy of certain differential equations in conformal field theory \cite{FBZ} and invariants of knots and $3$-manifolds \cite{Tur}. Small quantum groups, their representation categories and their semisimplification are related to Lie algebras in positive characteristic \cite{L - mod rep, AJS} and affine Lie algebras \cite{Fu}, and they have again applications to topology \cite{KL} and conformal field theory \cite{FGST, Len17}.

\smallbreak

While the classification of finite-dimensional semisimple Hopf algebras is still a very hard problem, we may ask for the classification of Hopf algebras $H$ with a \emph{given} maximal cosemisimple part $H_0$, called the coradical. For example, if ${H_0=\Bbbk G}$ is a group ring, then $H$ is called a pointed Hopf algebra. 
Andruskiewitsch and Schneider proposed a program called the \emph{Lifting method} to classify finite-dimensional Hopf algebras whose coradical is a subalgebra, that we sketch now. Every Hopf algebra comes with a coalgebra filtration $H_0\subseteq H_1\subseteq\cdots$ called the coradical filtration, where $H_0$ is the coradical. If we assume that $H_0$ is a subalgebra, then the associated graded coalgebra $\gr H$ is a Hopf algebra; we pay special attention to the subspace $\bar{H}_1=H_1/H_0$ and the subalgebra $\toba(\bar{H}_1)$ generated by $\bar{H}_1$, which is the Nichols algebra discussed below. As an example, for the (infinite-dimensional) quantum group $H=U_q(\mathfrak{g})$ the coradical $H_0$ is the group algebra spanned by the root lattice, the space $\bar{H}_1$ is spanned by the simple root vectors $E_i,F_i$, the Nichols algebra is the tensor product of the quantum Borel algebras $U_q(\mathfrak{g})^\pm$, and in the graded algebra $\gr H$ the relation $[E_i,F_i]=0$ holds in contrast to the nontrivial relation $[E_i,F_i]=K_i-K_i^{-1}$ in $H$. Then the classification proceeds in three steps: First, one classifies all finite-dimensional Nichols algebras over a fixed coradical $H_0$. Second, one checks if $H_0 \oplus \bar{H}_1$ generate the entire algebra $\gr H$ which is knows as the genereration in degee one problem. Third, one determines all possible Hopf algebras $H$ associated to each Nichols algebra in the first step, the so-called liftings. Using these ideas, Andruskiewitsch and Schneider obtained the whole classification of finite-dimensional Hopf algebras over $\mathbb{C}$ with coradical $H_0$ a group algebra of an abelian group of order not divisible by $2,3,5,7$ \cite{AS4}. The Nichols algebras appearing in this case are the Borel parts of small quantum groups, generation in degree one always holds and the possible liftings are described by deforming relations such as $[E_i,F_i]=0$ and $E_i^\ell=0$.

\smallbreak

The Nichols algebra $\toba(V)$ of a braided vector space $V$, or more generally of an object $V$ in a braided tensor category, is a quotient of the tensor algebra of $V$, which satisfies several universal properties. In particular, it is the unique Hopf algebra (in the braided tensor category) which is generated by $V$ and contains $V$ as the space of primitive elements. It is a difficult problem to determine the structure of the Nichols algebra of a given braided vector space, even to determine whether it is finite-dimensional. 
In view of the Lifting Method for pointed Hopf algebras, we need to consider Yetter-Drinfeld modules $V$ over a group algebra $H_0=\Bbbk G$. 
If $G$ is a finite abelian group, then the braiding is of diagonal type.
For braidings of diagonal type, Heckenberger classified the complex finite-dimensional Nichols algebras \cite{H-classif RS}. Besides the Nichols algebras $u_q(\mathfrak{g})^+$ coming from small quantum groups, he found Nichols algebras associated to Lie superalgebras, and several exceptional Nichols algebras involving roots of unity of order divisible by $2,3,5,7$; many of the later are now known to be related to Lie superalgebras in positive characteristic \cite{AA-diag-survey}. A main structural insight \cite{H-inv,HY,HS-book} is the existence of generalized root systems and Weyl groupoids. Roughly speaking, the generalization comes from the fact that different Weyl chambers might be in fact different, since reflections may change the braiding and even the Cartan matrix - an effect that already appears for contragredient Lie superalgebras. Nevertheless, finite Weyl groupoids can be classified \cite{CH} and show again a pattern of Lie algebra and Lie superalgebra series, plus low rank exceptions. 
The new Nichols algebras lead to new pointed Hopf algebras over abelian groups of order divisible by $2,3,5,7$, which include the small quantum groups $u_q(\mathfrak{g})$ associated to Lie superalgebras $\mathfrak{g}$. Using a presentation by generators and relations of these Nichols algebras of diagonal type, Angiono proved generation in degree one \cite{A-presentation}. The final answer about the liftings was established a few years later \cite{AnG} and involves the proof that every lifting is a cocycle deformation of the associated graded Hopf algebra in the sense of \cite{DT, Masuoka}.

\smallbreak

We now discuss Nichols algebras $\toba(V)$ of Yetter-Drinfeld modules $V$ over a non-abelian group $G$. In this case, irreducible Yetter-Drinfeld modules $V$ are parametrized by a conjugacy class $g^G$ in $G$ and a representation $\rho$ of the centralizer $G^g$ of $g$. A general Yetter-Drinfeld module $V$ is the direct sum of such irreducible modules, their number is called the rank of $V$ and they correspond to the simple roots of a generalized root system. The study of Nichols algebras over non-abelian groups starts with \cite{MS}, where $G$ is a Coxeter group and $g$ is a reflection. They treat two main examples: The symmetric group $\mathbb{S}_n$ has a single conjugacy class of reflections or transpositions with $\binom{n}{2}$ elements, accordingly $V$ is a Yetter-Drinfeld module of rank $1$ (i.e. irreducible) with dimension $\binom{n}{2}$. As it turns out, the associated Nichols algebras for $n=3,4,5$ have dimension $12,576,8294400$ and were considered by Fomin and Kirillov \cite{FK} in a very different context; for $n\geq 6$ they are conjecturally infinite-dimensional. On the other hand, the dihedral group $\mathbb{D}_4$ has two conjugacy classes of reflections, each with two elements. Accordingly, $V=V_1\oplus V_2$ has dimension $2+2$ and rank $2$. As it turns out, the generalized root system is of type $A_2$, indicating roughly that there is a space of braided commutators $V_{12}$, associated to the third conjugacy class with two elements, and all higher commutators vanish. Since the Nichols algebras of the irreducible modules $V_1,V_2,V_{12}$ have each dimension $4$, the Nichols algebra $\toba(V)$ has dimension $4^3$. 

The study now naturally branches into two directions: Nichols algebras of rank $1$, meaning of irreducible Yetter-Drinfeld modules, and Nichols algebras of rank $>1$ composed of the former via root system theory. In rank $1$, more finite-dimensional examples of Nichols algebras were discovered in \cite{Gr00, AG, HLV}, and later the research concentrated on successfully ruling out finite-dimensional Nichols algebras over most simple groups, see \cite{ACG} and the references there.

In rank $>1$, Lentner has developed the folding method for central extensions of Hopf algebras \cite{Len12,Len14}. When applied to Nichols algebras, it takes a Nichols algebra $\toba(V)$ over a group $\Gamma$ with an additional (twisted) symmetry group $Z$, and converts it to a Nichols algebra $\toba(\tilde{V})$ over the central extension of $\Gamma$. In particular, starting with $\Gamma$ an abelian group, $\toba(V)=u_q(\mathfrak{g})^+$ the Borel part of a small quantum group for $q^2=-1$ and $Z$ generated by a diagram automorphism of $\mathfrak{g}$, the folding method produces a Nichols algebra over central extensions $G$ of the abelian group $\Gamma$. There are cases
where $\mathfrak{g}$ is simple, namely ${^2}A_{2n+1}$ and ${^2}E_6$, and cases where $\mathfrak{g}$ consists of two copies of the same simply-laced Lie algebra interchanged by $Z=\mathbb{Z}_2$, which we denote ${^2}A_{n}^2$ and ${^2}D_{n}^2$ and ${^2}E_{n}^2,\;n=6,7,8$. The root system attached to $\toba(V)$ is the folded root system considered in Lie theory, with $Z$-orbits of roots becoming the new roots and $G$-conjugacy classes. For example, the smallest case ${^2}A_{2}^2$ reproduces the Nichols algebra over $\mathbb{D}_4$ discussed above, and the $6$ roots of $A_{2}^2$ become $3$ roots of $A_2$, each attached to a conjugacy class with two elements. 
In view of the exotic examples above, let us mention that these Nichols algebras are close to the quantum group case: By construction, their braiding can still be diagonalized, the group $G$ is close to  abelian, and the Weyl groupoid is again a Weyl group: they are \emph{standard}.

Around the same time, Heckenberger, Schneider and Vendramin started a systematic classification of finite-dimensional Nichols algebras over non-abelian groups over arbitrary fields using the root system theory, starting in rank $2$ in \cite{HS-rank2-1} and culminating in a full classification in rank $\geq 2$ \cite{HV-rank>2} in 2014. Their surprising observation was that the existence of a finite root system severely restricts the possible groups $G$, so that only very few of the (not yet fully classified) Nichols algebras of rank $1$ can appear in rank $\geq 2$. They found (in characteristic zero) three new exceptional Nichols algebras in rank $2$ and $3$, three infinite series $\alpha_{n}$, $\delta_{n}$, $\gamma_{n}$, and exceptional ones $\epsilon_{n}$, $\phi_4$ with standard root systems $A_n$, $D_n$, $C_n$, respectively $E_n$, $F_4$, whose conjugacy classes have two elements for short roots and one element for long roots. Note that these families precisely match the folding above.

\subsection{Goals and content}         

The purpose of this paper is to construct all pointed Hopf algebras whose group of group-like elements is a non-abelian group $G$ and the Nichols algebra is of standard type $\alpha_{n}$, $\delta_{n}$, $\gamma_{n}$, $\epsilon_{n}$ or $\phi_4$. Our first main result is that these series are all twists of the folded Nichols algebras by $2$-cocycles of the group $G$. This opens the possibility for a uniform treatment of these Nichols algebras.
Our second main result is that all Hopf algebras with coradical and a infinitesimal braiding as above are already generated in degree one. 
Since the folding method also applies to Hopf algebras, 
we construct a family of liftings inspired by those of diagonal type, and then fold accordingly.
Our last main result is that these are already all liftings. To achieve this, we give an explicit presentation of the Nichols algebras by generators and relations, which also allows us to prove that these Nichols algebras are rigid.

\smallbreak

We now discuss the content of this paper in more detail. 

In \S \ref{sec:Nichols} we review the preliminaries and the classification result by Heckenberger--Vendramin. 

Next we recall the folding construction in \S \ref{sec:folding}. Because we wish to work also with Hopf algebras, we give an exposition that emphasizes the origins of the folding method in this setting, as in \cite{Len12}: For a given Hopf algebra $H$ and a group of biGalois objects, the direct sum is again a Hopf algebra. Now we specialize to the case where the Hopf algebra is a smash product $H=\toba(V)\#\Bbbk \Gamma$ and the biGalois objects are based on a $2$-cocycle $\sigma$ on $\Gamma$ and a twisted symmetry of $V$.
Then, the result of the folding is again a smash product of the centrally extended group with the folded Nichols algebra. 
Our first main result is Theorem \ref{thm:Doi-Twist}, which states that the infinite families of Nichols algebras are twists of foldings over central extensions $Z\to G\to \Gamma$. By construction, $Z$ acts trivially on the folded Nichols algebras. We thus have to 
go through the classification of Heckenberger and Vendramin, assume a particular action of $Z$ and prove that the resulting groups $G$ possess sufficient cohomology to twist this action to zero. For large rank, the action of $Z$ turns out to be always trivial, but for small rank there are different possibilities depending on the structure of the group. For example, ${^2}A_2^2$ can be realized over central extensions such as $\mathbb{D}_4$ and $\mathbb{Q}_8$, in the former case we find nontrivial $Z$-action and sufficient cohomology, in the latter case we only find trivial $Z$-action. Since these computations are lengthy and involve tools from group cohomology, they are confined to the Appendix \ref{sec:appendix}.

In \S \ref{sec:gen-degree-one-trivial-action} we state our second main result, Theorem \ref{thm:gen-degree-one}, which gives a positive answer to the generation in degree one problem.

The main result of \S \ref{sec:generators-relations} is Theorem \ref{thm:Nichols-presentation-PBW-basis}, where the Nichols algebra is presented by generators and relations, and a PBW basis is obtained. Using the diagonal setting as inspiration, we give $G$-homogeneous generators and relations, which is necessary to deform the Hopf algebra structure. We also show that our Nichols algebras are rigid in the sense of \cite{AKM}.

In \S \ref{sec:liftings} we construct a big family of liftings using Hopf-Galois objects as in \cite{AAG,AnG,AnS}. Our last main result is Theorem \ref{thm:liftings-general}, where we show that this family exhaust all liftings, and that all liftings are cocycle deformation of the associated graded Hopf algebra.

In \S \ref{sec_outlook} we list some open questions and future directions of research. 

\section{Preliminaries}\label{sec:Nichols}
\subsubsection*{Conventions}
We denote $\N = \{1, 2, 3, \dots\}$ and $\N_0 =  \{0\} \cup \N$.
Given $k < \theta$ in $\N_0$ we put $\I_{k, \theta} = \{n\in \N_0: k\le n \le \theta \}$ and $\I_{\theta} = \I_{1, \theta}$. When $\theta$ is clear from the context we just write $\I=\I_\theta$.
The canonical basis of $\Z^{\theta}$ is denoted by $(\alpha_i)_{i\in \I_{\theta}}$. 

We work over an algebraically closed field $\Bbbk$ of characteristic zero and use $\Bbbk^\times$ to denote the group of non-zero elements. If $N\in\N$, we use $\G_N$ to denote the subgroup of $N$-th roots of unity; the subset of those with order $N$ is $\G_N'$.

Given a group $G$ and an element $g$, we use $g^G$ and $G^g$ to denote the conjugacy class and the centralizer of $g$, respectively. By $\widehat{G}$ we mean the group of characters, and $\Bbbk G$ stands for the group algebra. If $K$ is another group, then a pairing (also called a bicharacter) is a map 
$P:G \times K \to \Bbbk^\times$ such that for all $g,g'\in G$, $k,k'\in K$:
\begin{align*}
P(gg',k)&=P(g,k)P(g',k), & P(g,kk')&=P(g,k)P(g,k').
\end{align*}

A skew-polynomial algebra in variables $z_1, \dots, z_k$ is a quotient of the free algebra in these variables by an ideal generated by $z_iz_j-t_{ij}z_jz_i$, $1\leq i ,j\leq k$, for some $t_{ij}\in\Bbbk^\times$.

We denote Hopf algebras by tuples $(H,\mu,\Delta, \Ss)$ where $\mu$ is the multiplication, $\Delta$ the comultiplication and $\Ss$ the antipode, which we always assume bijective. The subspace of primitive elements is $\Pc(H)$.
The group of group-like elements is $G(H)$. If $ \delta \colon V \to H\ot V$ is a left
$H$-comodule we write $\delta(v)=v_{-1}\ot v_0$; for $g\in G(H)$ we put $V_{g} := \{v\in V: \delta(v) = g \otimes v\}$.
We refer to \cite{Mo-libro} for any unexplained terminology on Hopf algebras and to \cite[\S 2]{AAGMV} for preliminaries on Hopf-Galois objects and cocycle deformations.

\subsection{The Nichols algebra of a braided vector space}\label{subsec:nichols-general}
A \emph{braided vector space} is a pair $(V, c)$ where $V$ is a vector space and $c \in GL(V \otimes V)$ satisfies
$$(c\otimes \id)(\id\otimes c)(c\otimes \id) = (\id\otimes c)(c\otimes \id)(\id\otimes c).$$ 

By declaring the elements of $V$ to be primitive, the tensor algebra $T(V)$ becomes an $\N_0$-graded \emph{braided} Hopf algebra. There is a largest coideal $\Jc(V)$ among those that trivially intersect $\Bbbk \oplus V$; it happens to be graded so we denote $\Jc(V)=\bigoplus_{n\geq 2} \Jc^n(V)$. The Nichols algebra of $(V,c)$ is defined as the quotient $\toba(V) = T(V)/\Jc(V)$. This is again an $\N_0$-graded braided Hopf algebra, which is strictly graded as a coalgebra and generated by $V$ as an algebra, see \cite[\S 7]{HS-book}. Any intermediate quotient $\toba= T(V)/\Jc$ by an $\N_0$-homogeneous Hopf ideal $\Jc$ is called a \emph{pre-Nichols algebra} of $V$.

\medspace

The \emph{braided commutator} of $T(V)$ is defined by
\begin{align*}
[-,-]_c &= \text{mult} ( \id - c)\colon T(V) \ot T(V) \to T(V).
\end{align*}
If $u\in V$ and $v\in T(V)$ we denote $(\ad_c u) v =[u,v]_c$. In \S \ref{subsubsec:Yetter-Drinfeld-groups} will define $\ad_c u$ for arbitrary $u\in T(V)$.
For a fixed basis $(x_i)_{i\in I}$ of $V$ and $k\ge 2$ we set
\begin{align}\label{eq:iterated-ad-c}
x_{i_1\cdots i_k} &:= (\ad_c x_{i_1})\cdots (\ad_c x_{i_{k-1}}) x_{i_k}, & & i_j\in\I.
\end{align}

\begin{exa}
Given $\bq=(q_{ij})_{i,j\in \I}$ a matrix of elements of $\Bbbk^\times$, there is a braided vector space $(V,c^\bq)$ where $V$ has basis $(\xt_i)_{i\in \I}$ and $c^\bq$ is given by
\begin{align}\label{eq:def-diagonal-type}
c^\bq(\xt_i \ot \xt_j) &= q_{ij} \, \xt_j \ot \xt_i, & i,j &\in \I.
\end{align}
\end{exa}

A braided vector space is called \emph{of diagonal type} \cite{AS-cambr} if \eqref{eq:def-diagonal-type} holds in some basis of $V$ for some $\bq=(q_{ij})_{i,j \in \I}$. In this case, we denote the Nichols algebra of $(V,c)$ by $\toba_{\bq}$, which is now $\Z^\I$-graded; we refer to $\bq$ as the \emph{braiding matrix}. The Dynkin diagram of $\bq$ is a graph with $\I$ as the set of vertices, each vertex $i$ labelled with $q_{ii}$. There is an edge between $i\ne j$ if and only if $\widetilde{q}_{ij}:=q_{ij}q_{ji} \ne 1$; such an edge is labelled with this scalar. 

We say that $\bq$ is \emph{of Cartan type} \cite{AS-adv} if there is a Cartan matrix 
$\ba= (a_{ij})$ such that
\begin{align*}
q_{ij}q_{ji} &= q_{ii}^{a_{ij}}, & \mbox{for all }i,j &\in \I.
\end{align*}
If some $q_{ii}$ is not a root of unity, then the integers $a_{ij}$ are uniquely determined. Otherwise we impose $-\ord q_{ii} < a_{ij} \leq 0$ for all $j \neq i$. In this case we say that $\bq$ is of Cartan type $\ba$. 

\medskip
Although braided vector spaces of Cartan type seem quite \emph{simple}, the structure of the corresponding Nichols algebra is related either with quantized enveloping algebras (when the entries of $\bq$ are not roots of unity), or with Frobenius-Lusztig kernels. 

The following example of Cartan type will be particularly relevant in later sections.
\begin{exa}\label{ex:typeADE-1}
Fix a finite Cartan matrix $\ba$ with simply-laced Dynkin diagram. Assume that $\bq=(q_{ij})$ satisfies the following conditions:
\begin{align}\label{eq:typeADE-1}
q_{ii}&=-1; & q_{ij}q_{ji} &=
\begin{cases} 
-1, & a_{ij} =-1, \\
1, & a_{ij} =0;
\end{cases}
& i &\neq j\in \I.
\end{align}
Then $\bq$ is of Cartan type $\ba$.

Let $\beta_1<\beta_2<\dots <\beta_M$ be a convex order on the set of positive roots $\Delta_+$ of $\ba$. In \cite{A-presentation} we can find a \emph{root vector} $\xt_{\beta}\in\toba_{\bq}$ for each $\beta\in\Delta_+$, of $\Z^\I$-degree $\beta$, obtained recursively as braided commutator of root vectors with smaller degree.

In this case, the Nichols algebra $\toba_{\bq}$ is presented by generators $(\xt_i)_{i\in\I}$ and relations
\begin{align}\label{eq:typeADE-1-rels-1}
\xt_{\alpha}^2&=0, & \alpha &\in\Delta_+; \\
\label{eq:typeADE-1-rels-2}
[\xt_{ijk}, \xt_j]_c&=0, & a_{ij} &=a_{jk} =-1; \\
\label{eq:typeADE-1-rels-3}
\xt_{ij}&=0, & a_{ij} &=0.
\end{align}
Moreover, by \cite{A-presentation,AA-diag-survey}, a basis for $\toba_{\bq}$ is given by the set
\begin{align}\label{eq:typeADE-1-PBWbasis}
\{ \xt_{\beta_1}^{n_1}\xt_{\beta_2}^{n_2} \dots \xt_{\beta_M}^{n_M} \, | \, n_i\in\{ 0,1 \}  \}.
\end{align}
\end{exa}

\begin{remark}\label{rem:diagonal-pre-Nichols} Let $\bq$ as above. In  \S \ref{sec:generators-relations} we will need the following constructions, due to  \cite{A-pre-Nichols}.
\begin{enumerate}[leftmargin=*,label=\rm{(\Roman*)}]
\item\label{item:diagonal-pre-Nichols-dist} The \emph{distinguished pre-Nichols} algebra $\wtoba_{\bq}$ is the quotient of $T(V)$ by \eqref{eq:typeADE-1-rels-2}, \eqref{eq:typeADE-1-rels-3} and
\begin{align}
\label{eq:typeADE-1-rels-distinguished}
\xt_{iij}&=0, & a_{ij} &=-1.
\end{align}
The set 
$\{ \xt_{\beta_1}^{n_1}\xt_{\beta_2}^{n_2} \dots \xt_{\beta_M}^{n_M} \, | \, n_i\in\N_0  \}$ is a basis of $\wtoba_{\bq}$.

\item\label{item:diagonal-pre-Nichols-coinv} Let $\pi:\wtoba_{\bq}\twoheadrightarrow \toba_{\bq}$ be the canonical projection. 
By \cite{A-pre-Nichols}, the subalgebra of coinvariants $\Zc_\bq:=\wtoba_{\bq}^{\co \pi}$ 
is a $q$-polynomial algebra with generators $\xt_{\beta}^2$, $\beta\in\Delta_+$.

\item\label{item:diagonal-pre-Nichols-quotient} 
Let $\htoba_{\bq}$ be the quotient of $T(V)$ by \eqref{eq:typeADE-1-rels-2}, \eqref{eq:typeADE-1-rels-3} and $\xt_i^2$, $i\in\I$.
Then $\htoba_{\bq}$ is a pre-Nichols algebra, which coincides with the quotient of $\wtoba_{\bq}$ by $\xt_i^2$, $i\in\I$.  The set 
$$ \{ \xt_{\beta_1}^{n_1}\xt_{\beta_2}^{n_2} \dots \xt_{\beta_M}^{n_M} \, | n_i\in\{0,1\} \text{ if }\beta_i \text{ is simple, }n_i\in\N_0 \text{ otherwise} \} $$ 
is a basis of $\htoba_{\bq}$, so its Hilbert series is
\begin{align*}
\mH_{\htoba_{\bq}} = \left(\prod_{\beta \in \Delta_+^{\bq}-\{\alpha_i\} } \frac{1}{1-t^{\beta}}\right) \left(\prod_{i\in\I} 1+t_i \right).
\end{align*}
\end{enumerate}
\end{remark}

\subsection{Nichols algebras over non-abelian groups}\label{subsec:Nichols-non-abelian-groups}
The goal of this subsection is to introduce the notion of Weyl grupoids. These play a fundamental role in the classification achieved in \cite{HV-rank>2}. We refer to the book \cite{HS-book} for details and unexplained terminology.

\subsubsection{Yetter-Drinfeld modules over groups}\label{subsubsec:Yetter-Drinfeld-groups}
For a group $G$, the category of Yetter-Drinfeld modules $\ydG$ consist of $G$-graded vector spaces $V=\bigoplus_{g\in G} V_g$ 
endowed with a $G$-action such that $h \cdot V_g \subset V_{hgh^{-1}}$ for all $h,g \in G$. This is a braided tensor category where the braiding 
$c_{V,W} \colon V\otimes W \to W \ot V$ is given by $c(v\ot w) = g\cdot w \ot v$ for $v\in V_g$ and $w\in W$. 
We recall that a $G$-grading on a vector space $V=\bigoplus_{g\in G} V_g$ is equivalent to a $\Bbbk G$-comodule structure $\delta \colon V \to V\ot \Bbbk G$, 
declaring $\delta(v)=g \ot v$ if and only if $v \in V_g$. The \emph{support} of $V$ is 
\begin{align*}
\operatorname{supp} V=\{g \in G \vert \ V_g\neq0\}.
\end{align*}

Let $(R,\mu,\Delta,\Ss)$ be a  Hopf algebra in $\ydG$. The braided commutator defined above satisfies the following identity:
If $u\in R_g$ and $v\in R_h$ for some $g,h\in G$, then for any $w \in R$, 
\begin{align}\label{eq:braided-commutator-iteration}
\left[[u,v]_c, w\right]&=  \left[u,[v,w]_c\right] - (g\cdot v)[u,w]_c + [u,h\cdot w ]_c v.
\end{align}
$R$ admits a \emph{braided} adjoint representation $\ad_c \colon R \to \End (R)$ given by 
\begin{align*}
(\ad_c u) v&= \mu (\mu \ot \Ss) (\id \ot c) (\Delta \ot \id) (u \ot v) , &u,& v \in R. 
\end{align*}
When $u$ is primitive, $\ad_c u$ and $[u,-]_c$ coincide. Notice also that
\begin{align}\label{eq:group-action-in-adjoint}
g\cdot((\ad_c u) v)&= (\ad_c g\cdot u) (g\cdot v), &g\in G, \ u\in \Pc(R),&\ v \in R. 
\end{align}

\subsubsection{The Nichols algebra of a Yetter-Drinfeld module} \label{subsubsec:Nichols-in-Yetter-Drinfeld}
Each $V\in \ydG$ is a braided vector space, so it has a Nichols algebra $\toba(V)$ as discussed in 
\S \ref{subsec:nichols-general}. In this setting, $T(V)$ and $\toba(V)=T(V) / \Jc(V)$ turn out to be $\N_0$-graded Hopf algebras in $\ydG$. 

\subsubsection{Skew derivations} \label{subsubsec:skew-derivations} There is a criterion to decide if any given element of $T(V)$ belongs to $\Jc(V)$. Fix a basis ${x_1, \dots, x_r}$ of $V$ with $x_i$ of degree $h_i$. For each $i$, we define a skew derivation $\partial_{x_i}\in \End (T(V))$ recursively in $V^{\ot n}$, $n\geq 0$. For $n=0$, put $\partial_{x_i} (1)=0$; in $V^{\ot1}$ put $\partial_{x_i} (x_j)=\delta_{i,j}$ and in general define
\begin{align*}
\partial_{x_i}(xy)&=x \partial_{x_i}(y) + \partial_{x_i}(x)  h_i \cdot y, &x,&y \in T(V).
\end{align*}
Then, for any $n\geq2$, we have $x\in \Jc^n(V)$ if and only if $\partial_{x_i}(x)\in \Jc^{n-1}(V)$ for all $i$.
The compositions of these derivations with the braided adjoint action satisfy
\begin{align}\label{eq:derivation-in-adjoint}
\partial_{x_i}((\ad_c u) v)= u \partial_{x_i}(v) -  \partial_{x_i}(g\cdot v) h_i\cdot u, \quad u\in V_g, \, v \in T(V),  \, i\in \I_r. 
\end{align}

\subsubsection{Simple Yetter-Drinfeld modules}\label{subsubsec:Yetter-Drinfeld-simples}
Fix $g\in G$ and $(V, \rho)$ an irreducible representation of the centralizer $G^g$. We consider the induced $G$-module $ \Bbbk G \ot_{G^g} V$ endowed with the $G$-grading determined by declaring the degree of $x \ot v$ to be $ gxg^{-1}$ for all $x \in G$ and $v \in V$. This is a simple Yetter-Drinfeld module over $G$, which is denoted by $M(g^G, \rho)$ since it depends on the conjugacy class $g^G$ rather than the element $g$ itself. Moreover, all the simple objects of $\ydG$ arise in this way, and if $G$ is finite the category $\ydG$ is semisimple. 
There is a concrete description of $M(g^G, \rho)$ in \cite[Example 24]{A-leyva}, which we will use several times to perform computations in $\ydG$.

\begin{exa}\label{ex:simple-YD-dim2}
Let $G$ be a finite group. Assume that $g\in G$ is such that $g^G=\{g,g\kappa\}$ for some $\kappa \neq e \in G$. Then there exists $g_0\in G$ such that $g_0g=\kappa gg_0$, $G^{g}=G^{g\kappa}=G^{g^{-1}}$ is a subgroup of index two since $G=G^g\cup g_0G^g$, and $\kappa\in Z(G)$, $\kappa^2=e$. 

Let $\chi$ be a one-dimensional representation of $G^g$, that is $V=\Bbbk$ and $\chi\in \widehat{G^{g}}$, and $M= M(g^G, \chi)$. Then $\dim M=2$: we may fix a basis $\{x,y\}$ such that $y=g_0\cdot x$ and the coaction is given by $\rho(x)=g\ot x$, $\rho(y)=g\kappa \ot y$. 
As $hg_0, g_0^{-1}h\in G^g$ for all $h\notin G^g$, the action satisfies
\begin{align*}
h \cdot x &= \begin{cases}
\chi(h) x, & h\in G^g, \\
\chi(g_0^{-1}h) y, & g\notin G^g;
\end{cases}
& h \cdot y &= \begin{cases}
\chi(g_0^{-1}hg_0) y, & h\in G^g, \\
\chi(hg_0) x, & h\notin G^g;
\end{cases}
& &k\in G.
\end{align*}
Thus $M$ is of diagonal type with braiding matrix 
$\left[\begin{smallmatrix} \chi(g) & \chi(g\kappa) \\ \chi(g\kappa) & \chi(g) \end{smallmatrix}\right]$; from this we deduce that  
$\dim \toba(M)<\infty$ if and only if $\chi(g)\in\{-1\}\cup\G_3'$.
\end{exa}

\subsubsection{Weyl groupoid}
Next we recall the definition of the Weyl groupoid of a non-simple element $M\in\ydG$ such that $\dim \toba(M)<\infty$.
Let $M=\oplus_{i\in\I} M_i$, where each summand $M_i$ is simple: $M_i=M(g_i^G,\rho_i) \in \ydG$ for some $g_i\in G$, $\rho_i$ an irreducible representation of $G^{g_i}$. 

For each $i\ne j \in\I$, set $(\ad M_i)^0 M_j := M_j$, and for $n\in\N$,
\begin{align*}
(\ad M_i)^n M_j &:= \{ (\ad_c m_1) \cdots (\ad_c m_n) m | m_{\ell}\in M_i, m\in M_j \} \subset \toba (M).
\end{align*}
The \emph{generalized Cartan matrix} of $M$ is  $C^M=(c_{ij}^M)\in\Z^{\I\times \I}$, where
\begin{align}\label{eq:cij-defn}
c_{ii}^M &= 2, & c_{ij}^M &= -\min \{n\in \N_0| (\ad M_i)^{n+1} M_j=0 \}, & &j\ne i.
\end{align}
The $i$-\emph{reflection} of $M$ is $\rho_i M=\sum_{j\in\I} \underline{M}_j$, where 
\begin{align*}
\underline{M}_j=\begin{cases} (\ad M_i)^{-c_{ij}^M} M_j, & j\ne i; \\ M_i^*, & j=i. \end{cases}
\end{align*}
Each $\underline{M}_j\in\ydG$ is simple and $\dim \toba(M)=\dim \toba(\rho_i M)$.

These reflections generate the \emph{Weyl groupoid} of $M$, see \cite{HS-book}.

\subsection{Heckenberger--Vendramin classification}
Let $M=\oplus_{i\in \I_{\theta}} M_i$ be a finite dimensional Yetter-Drinfeld module over a non-abelian group $G$, where each $M_i$ is simple. Here $\theta$ is called the \emph{rank} of $M$.

In \cite[Theorem 2.5]{HV-rank>2} the authors classify Yetter-Drinfeld modules as above of rank at least $2$, such that the associated Nichols algebra is finite-dimensional. To be precise, one needs to assume that the support of $M$ generates $G$, and impose a mild non-degeneracy condition on the braiding between different summands of $M$. Up to a few exceptions in ranks two and three, the classification consists on families $\alpha_{\theta}$, $\gamma_{\theta}$ and $\delta_{\theta}$ of arbitrary rank, and types $\epsilon_{\theta}$, $\theta=6,7,8$, $\phi_4$, which resemble the classification of finite-dimensional Lie algebras. These Yetter-Drinfeld modules are invariant under the Weyl groupoid action; we say that they are \emph{standard}, adopting the terminology used for diagonal type, see \cite{AA-diag-survey}. In this paper we study this class; next we give an explicit description of each module.

\subsubsection{Types $\alpha_{\theta}$, $\delta_{\theta}$ and $\epsilon_{\theta}$}\label{subsubsec:type-ADE}

Fix a simply-laced indecomposable Cartan matrix $\ba=(a_{ij}) \in \Z^{\theta \times \theta}$ of finite type; i.e. of type $A_{\theta}$, $\theta \geq 2$, $D_{\theta}$, $\theta \geq 4$, or $E_{\theta}$ for $\theta\in\I_{6,8}$. Let $\Delta_+$ be the set of positive roots.
Following  \cite[Lemma 6.2]{HV-rank>2} we describe a Yetter-Drinfeld module  $M=\oplus_{i\in\I} M_i$ over a non-abelian group $G$ with simply-laced skeleton and Cartan graph of standard type $\ba$. Assume there exist 
\begin{itemize}[leftmargin=*]
\item $\kappa \in Z(G)$ such that $\kappa\ne 1$, $\kappa^2=1$,
\item $g_i \in G$ with $g_i^G=\{g_i, \kappa g_i\}$ for all $i\in \I_{\theta}$,
\item $\chi_i\in \widehat{G^{g_i}}$ such that $\chi_i(g_i)=-1$ for all $i\in \I_{\theta}$,
\end{itemize}
satisfying the following: 
\begin{align} 
\label{eq:typeADE-aij=-1}
g_ig_j&=\kappa g_jg_i, &  \chi_i(\kappa g_j^2)\chi_j(\kappa g_i^2)&=1, &a_{ij}&=-1;
\\ \label{eq:typeADE-aij=0}
g_ig_j&=g_jg_i, &  \chi_i(g_j)\chi_j(g_i) &=1, &  a_{ij}&=0;
\\ \label{eq:typeADE-kappa}
&& \chi_i(\kappa)\chi_j(\kappa) &=1, &  a_{ij}&=0;
\end{align}
For $i\in\I_\theta$ let $M_i=M(g_i^G,\chi_i)\in \ydG$, which has a basis $\{x_i, x_{\un{i}}\}$ with coaction $x_i\mapsto g_i \ot x_i$, $x_{\un{i}}\mapsto \kappa g_i \ot x_{\un{i}}$. As in \cite[Example 37]{A-leyva}  the braidings $c_{M_i, M_j}$ for $i,j\in\I_{\theta}$ are determined by 
\begin{align}\label{eq:ADE-self-braiding}
\begin{bmatrix}c(x_i \ot x_i) & c(x_i \ot x_{\un{i}}) \\ c(x_{\un{i}} \ot x_i) & c(x_{\un{i}} \ot x_{\un{i}}) \end{bmatrix}&=\begin{bmatrix} -x_i \ot x_i & -\chi_i(\kappa) x_{\un{i}} \ot x_i \\ -\chi_i(\kappa)  x_i \ot x_{\un{i}} & -x_{\un{i}} \ot x_{\un{i}} \end{bmatrix};
\\ \label{eq:ADE-connected-braiding}
\begin{bmatrix}c(x_i \ot x_j) & c(x_i \ot x_{\un{j}}) \\ c(x_{\un{i}} \ot x_j) & c(x_{\un{i}} \ot x_{\un{j}}) \end{bmatrix}&=
\begin{bmatrix} x_{\un{j}} \ot x_i & \chi_j(g_i^2) x_j \ot x_i \\ \chi_j(\kappa)  x_{\un{j}} \ot x_{\un{i}} & \chi_j(\kappa g_i^2)x_j \ot x_{\un{i}} \end{bmatrix}, && a_{ij}=-1;
\\ \label{eq:ADE-disconnected-braiding}
\begin{bmatrix}c(x_i \ot x_j) & c(x_i \ot x_{\un{j}}) \\ c(x_{\un{i}} \ot x_j) & c(x_{\un{i}} \ot x_{\un{j}}) \end{bmatrix} &=
\begin{bmatrix} \chi_j(g_i) x_j \ot x_i & \chi_j(\kappa g_i) x_{\un{j}} \ot x_i \\ \chi_j(\kappa g_i)  x_j \ot x_{\un{i}} & \chi_j(g_i)x_{\un{j}} \ot x_{\un{i}} \end{bmatrix}, && a_{ij}=0.
\end{align}
The Cartan matrix of $M=\oplus_{i\in\I} M_i$ is $\ba$, and we have $\dim \toba(M) = 2^{2|\Delta_+|}$.

\begin{remark}
\begin{itemize}[leftmargin=*]
\item The braiding on $M_i$ is of diagonal type and $c_{M_i,M_i}^2=\id$.
\item If $a_{ij}=0$, then $M_i\oplus M_j$ is of diagonal type and  $c_{M_i,M_j}c_{M_j,M_i}=\id_{M_j\ot M_i}$.
\item If $a_{ij}=-1$, then $M_i\oplus M_j$ is of diagonal type if and only if $\chi_i(\kappa)=\chi_j(\kappa)=1$, but here 
$c_{M_i,M_j}c_{M_j,M_i}\ne \id_{M_j\ot M_i}$. 
\end{itemize}
\end{remark}

\subsubsection{Type $\gamma_{\theta}$, $\theta\ge 3$}\label{subsubsec:type-C}
Following \cite[Lemma 7.6]{HV-rank>2} we describe, for each rank $\theta \geq 3$, a Yetter-Drinfeld module $M$ of type $\gamma_{\theta}$ over a non-abelian group $G$. 
Assume there exist 
\begin{itemize}[leftmargin=*]
\item $\kappa \in Z(G)$ such that $\kappa\ne 1$, $\kappa^2=1$,
\item $g_1, \dots, g_\theta \in G$ with $g_i^G=\{g_i, \kappa g_i\}$ for $i\in \I_{\theta -1}$ and $g_\theta^G=\{g_\theta\}$,
\item $\chi_i\in \widehat{G^{g_i}}$ such that $\chi_i(g_i)=-1$,
\end{itemize}
satisfying \eqref{eq:typeADE-aij=-1}, \eqref{eq:typeADE-aij=0} and \eqref{eq:typeADE-kappa} for $i,j\in\I_{\theta-1}$, and
\begin{align}
\label{eq:typeC-aij=0}
g_{i}g_{\theta}&= g_{\theta}g_{i}, & \chi_i(g_\theta)\chi_\theta(g_i) &=1, &&i<\theta-1;
\\ 
\label{eq:typeC-theta-1-theta}
&&\chi_{\theta-1}(g_\theta)\chi_{\theta}(g_{\theta-1})&=-1.
\end{align}
Let $M_i=M(g_i^G,\chi_i) \in \ydG$. Then  $M=\oplus_{i\in\I} M_i$ is of type $\gamma_\theta$ and $\dim \toba(M) = 2^{ 2 \theta^2 - \theta}$. 
Notice that $M_1\oplus \dots\oplus M_{\theta-1}$ is of type $\alpha_{\theta-1}$ so, by \S \ref{subsubsec:type-ADE},  we have a basis $\{x_i, x_{\un{i}}\}$ of $M_i$ such that for $i, j \in\I_{\theta-1}$ the braiding $c_{M_i, M_j}$ is determined by \eqref{eq:ADE-self-braiding}, \eqref{eq:ADE-connected-braiding} and \eqref{eq:ADE-disconnected-braiding}. On the other hand $M_\theta=\Bbbk\{x_\theta\}$ is one-dimensional concentrated in degree $g_\theta\in Z(G)$. The braidings $c_{M_{\theta}, M_{\theta}}$, $c_{M_i, M_{\theta}}$, $c_{M_{\theta}, M_i}$, $i\in \I_{\theta-1}$, are determined by
\begin{align}\label{eq:C-xtheta-self-braiding}
c(x_\theta \ot x_\theta) &= - x_\theta \ot x_\theta;
\\ \label{eq:C-xi-xtheta-braiding}
c(x_i \ot x_\theta) &=\chi_\theta(g_i) x_\theta \ot x_i, && c(x_{\un{i}} \ot x_\theta) = \chi_\theta(\kappa g_i) x_\theta \ot x_{\un{i}};
\\ \label{eq:C-xtheta-xj-braiding}
c(x_\theta \ot x_i)&=\chi_i(g_\theta)x_i \ot x_\theta, && c(x_\theta \ot x_{\un{i}})  = \chi_i(g_\theta) x_{\un{i}} \ot x_\theta.
\end{align}

\subsubsection{Type $\phi_{4}$} \label{subsubsec:type-F4}

Following \cite[Lemma 9.2]{HV-rank>2} we describe a Yetter-Drinfeld module $M$ over a non-abelian group $G$ with Cartan matrix of type $F_4$. 
Assume there exist 
\begin{itemize}[leftmargin=*]
\item $\kappa \in Z(G)$ such that $\kappa\ne1$, $\kappa^2=1$,
\item $g_1, \dots, g_4 \in G$ with $g_i^G=\{g_i, \kappa g_i\}$ for $i=1,2$ and $g_i^G=\{g_i\}$ for $i=3,4$,
\item $\chi_i\in \widehat{G^{g_i}}$ such that $\chi_i(g_i)=-1$ 
\end{itemize}
satisfying the following: 
\begin{align}
\label{eq:typeF4-aij=0}
\chi_1(g_3)\chi_3(g_1) &=\chi_1(g_4)\chi_4(g_1)=\chi_2(g_4)\chi_4(g_2)=1;
\\
\label{eq:typeF4-2-3-4}
\chi_3(g_4)\chi_4(g_3)&=\chi_2(g_3)\chi_3(g_2) =-1;
\\ 
\label{eq:typeF4-1-2}
g_1g_2&=\kappa g_2g_1, \quad  \chi_1(\kappa g_2^2)\chi_2(\kappa g_1^2)=1.
\end{align}
Let $M_i=M(g_i^G,\chi_i)$. The Cartan matrix of $M=\oplus_{i\in\I} M_i$ is of type $F_4$  and 
$\dim \toba(M) = 2^{36}$. The braidings $c_{M_i, M_j}$ are given as in types $\alpha_{\theta}$, $\gamma_{\theta}$, depending on $i,j$.

\subsection{On the structure of the group \texorpdfstring{$G$}{}}

Fix an abelian group $\Gamma$ and $x,y \in \Z_2\times \Z_2$ such that $\Z_2\times \Z_2=\langle x,y\rangle$. For every $u,v\in \Gamma$ and every $\kappa\in \Gamma$ such that $\kappa^2=e$, there exists a 2-cocycle $\beta\in H^2(\Z_2\times \Z_2, \Gamma)$ such that
\begin{align*}
\beta(x,x) &= u, & \beta(x,y) &=\kappa, & \beta(y,x) &=e, & \beta(y,y) &=v.
\end{align*}
We denote by $\Gamma_{u,v,\kappa}$ the associated central extension of $\Gamma$ by $\Z_2\times \Z_2$:
\begin{align}\label{eq:defn-gamma-uvkappa}
\Gamma \hookrightarrow & \Gamma_{u,v,\kappa} \twoheadrightarrow \Z_2\times \Z_2=\langle x,y\rangle,
\end{align}
where $[x,y]=\kappa$, $x^2=u$, $y^2=v$.

\medbreak 
Next we describe some general features of a group $G$ realizing the braidings described in \S \ref{subsubsec:type-ADE}, \ref{subsubsec:type-C}, \ref{subsubsec:type-F4}. 
Let $M=\oplus_{i\in\I} M_i$ be of type $\alpha_{\theta}$, $\gamma_{\theta}$, $\delta_{\theta}$, $\epsilon_{\theta}$, or $\phi_4$.
Here $M_i=M(g_i^G,\chi_i) \in \ydG$, hence $\supp M=\bigcup_{i\in\I} g_i^G$. Note that
\begin{itemize}[leftmargin=*]
\item There exists $\kappa\in Z(G)$ such that $\kappa^2=e$ and $g_i^G=\{g_i,\kappa g_i\}$ for all $i$ such that $g_i\notin Z(G)$.
\item For $i$ with $g_i^G=\{g_i,\kappa g_i\}$, there exists $j\ne i$ such that $g_ig_j=\kappa g_jg_i$; hence $g_j^G=\{g_j,\kappa g_j\}$.
\end{itemize}

The relevance of the central extensions constructed above is explained next.

\begin{lemma}\label{lem:structure-G}
Let $i\ne j\in\I$ be such that $g_ig_j=\kappa g_jg_i$, and let $N=G^{g_i}\cap G^{g_j}$. Then
\begin{enumerate}[leftmargin=*,label=\rm{(\alph*)}]
\item\label{item:structure-G-a} The subgroup $N$ is normal, and $G/N\simeq \Z_2\times\Z_2$.

\item\label{item:structure-G-b} If $G=\langle \supp M\rangle$, then $N$ is an abelian subgroup, generated by the elements
\begin{align*}
& g_k, && \text{ for all }k\in\I \text{ such that }g_kg_i=g_ig_k, \, g_kg_j= g_jg_k, \\
& g_kg_i, && \text{ for all }k\in\I \text{ such that }g_kg_i=g_ig_k, \, g_kg_j=\kappa g_jg_k, \text{ and}\\
& g_kg_j, && \text{ for all }k\in\I \text{ such that }g_kg_j=g_jg_k,  \, g_kg_i=\kappa g_ig_k.
\end{align*}
\item\label{item:structure-G-c} If $N$ is abelian, then $G\simeq N_{g_i^2,g_j^2,\kappa}$. 
\end{enumerate}
\end{lemma}
\pf
For \ref{item:structure-G-a}, we note first that $[G:G^{g_i}]=[G:G^{g_j}]=2$, so both $G^{g_i}$ and $G^{g_j}$ are normal subgroups, 
and $[G:N_{ij}]=4$ since $G=G_iG_j$. Let $g\in G$:
\begin{itemize}
\item If $gg_ig^{-1}=g_i$, $gg_jg^{-1}=g_j$, then $gN=N$.
\item If $gg_ig^{-1}=g_i$, $gg_jg^{-1}=\kappa g_j$, then $gN=g_iN$.
\item If $gg_ig^{-1}=\kappa g_i$, $gg_jg^{-1}=g_j$, then $gN=g_jN$.
\item If $gg_ig^{-1}=\kappa g_i$, $gg_jg^{-1}=\kappa g_j$, then $gN=g_ig_jN$.
\end{itemize}
Hence $G/N=\{N,g_iN,g_jN,g_ig_jN\}$. Since $g_i^2,g_j^2\in N$,  we get $G/N\simeq \Z_2\times \Z_2$.

\ref{item:structure-G-b} and \ref{item:structure-G-c} are straightforward.
\epf

A more detailed description of these groups is postponed to the Appendix. 

\subsubsection{The parity vector} An important invariant of $M=\oplus_{i\in\I} M_i \in \ydG$ is 
\begin{align}\label{eq:action-kappa}
\mathtt{P}:=(\chi_i(\kappa))_{1\leq i\leq \theta} \in \{\pm 1\}^\theta.
\end{align}
If $\mathtt{P}=(1,\ldots,1)$, then the $G$-action on $M$ factors to one of $G/\langle \kappa \rangle$. We show next that this happens in most of the cases, in which case $M$ is a braided vector space of diagonal type.

\begin{lemma}\label{lem:trivial-action-kappa-big-rank}
Let $M$ be either of type $\alpha_{\theta}$, $\theta\ge 4$, $\gamma_{\theta}$, $\theta\ge 5$, $\delta_{\theta}$, $\theta\ge 5$ or $\epsilon_{\theta}$, $\theta=6,7,8$. 
Then $\mathtt{P}=(1,\ldots,1)$. In other words, $\kappa$ acts trivially on $M$.
\end{lemma}

\pf
We consider first the case $\alpha_4$. By \eqref{eq:typeADE-aij=0}, $g_1,g_2\in G^{g_4}$, and by \eqref{eq:typeADE-aij=-1}, $\kappa=g_1g_2g_1^{-1}g_{2}^{-1}$; hence
$\chi_4(\kappa)=1$. By \eqref{eq:typeADE-kappa}, $\chi_i(\kappa)=1$ for all $i\in\I_4$.

Now the proof for $\alpha_\theta$, $\theta\ge 5$, $\delta_\theta$, $\theta\ge 5$ or $\epsilon_\theta$, $\theta=6,7,8$ follows from taking into account all the submodules of type $\alpha_4$. The same fact says that, for type $\gamma_\theta$, $\chi_i(\kappa)=1$ for all $i\in\I_{\theta-1}$. Next we use that $g_{\theta-2},g_{\theta-1}\in G^{g_\theta}$ and $\kappa=g_{\theta-1}g_{\theta-2}g_{\theta-1}^{-1}g_{\theta-2}^{-1}$ to deduce that $\chi_\theta(\kappa)=1$.
\epf

\bigbreak

\section{Foldings of Nichols algebras and trivializing the action of the center}\label{sec:folding}
Motivated by Lemma \ref{lem:trivial-action-kappa-big-rank}, we pay special attention to Yetter-Drinfeld modules where $\kappa$ acts trivially. We will show that these examples are related to diagonal braidings (of Cartan type) via the folding construction for Nichols algebras, developed by the second author in \cite{Len12,Len14}. Then we show that for the other cases the action of the central element can be trivialized via a twist. First, we introduce basic notions needed for the folding construction.

\subsection{Categorical action on Yetter-Drinfeld modules} \label{subsec:categorical-action}
Given a group $G$ and a $2$-cocycle $\sigma \in Z^2(G, \Bbbk^\times)$, we get a pairing $b_\sigma \colon G \times G \to \Bbbk^\times$ given by $b_\sigma (g, h) = \sigma(hgh^{-1},h) \sigma^{-1} (h,g)$.
Each $2$-cocycle $\sigma$ yields a tensor functor $F_\sigma \colon \ydG \to \ydG$ as follows:
\begin{itemize}[leftmargin=*]
\item for an object $M$, let $F_\sigma (M)$ denote the same $\Bbbk G$-comodule;
\item the $\Bbbk G$-module structure on $F_\sigma (M)$ is given by
\begin{align*}
g \cdot_\sigma m&= b_\sigma(g, m_{-1} )g\cdot m_0, &g\in  G,  \ m\in M;
\end{align*}
\item  on Hom spaces $F_\sigma$ is the identity, thus $F_\sigma$ is $\Bbbk$-linear, faithful and exact;
\item the monoidal structure $J_\sigma \colon F_\sigma (M) \ot F_\sigma(N) \to F_\sigma(M \ot N) $ is defined by
\begin{align*}
J_\sigma (m\ot n) &= \sigma( m_{-1} , n_{-1} ) m_0\ot n_0 , &\ m\in M, n \in N.
\end{align*}
\end{itemize}
We note that $J_\sigma$ satisfy the hexagon axiom thanks to the cocycle condition on $\sigma$. 

Given $V,W\in\ydG$, one can see that $F_{\sigma}\otimes F_{\sigma}$ commutes with $c_{V,W}$ if and only if
$\sigma(ghg^{-1}, g) \sigma(hgh^{-1},h)= \sigma(g,h)\sigma(h, g)$ for all $g\in\supp V$, $h \in \supp W$.
Hence, $F_\sigma$ is braided if and only if that equality holds for all $g, h \in  G$, 
which certainly happens if $ G$ is abelian.
Anyhow, these functors patch together to an action by tensor autoequivalences of the group $Z^2( G, \Bbbk^\times)$ on $\ydG$:
\begin{itemize}[leftmargin=*]
\item the trivial cocycle acts as the identity, and
\item $F_\eta F_\sigma = F_{\eta\sigma}$ for all $\eta, \sigma \in Z^2( G, \Bbbk^\times)$.
\end{itemize}

Given a Hopf algebra $\toba$ in $\ydG$ we set, by abuse of notation,
\begin{align}\label{eq:cocycle-extension-to-bosonization}
&\sigma:\toba\#\Bbbk G \otimes \toba\#\Bbbk G \to \Bbbk,& \sigma(x\# g \ot y\# h)&=\epsilon(x)\epsilon(y)\sigma(g,h).
\end{align}

\begin{lemma}\label{lem:sigma}
\begin{enumerate}[leftmargin=*,label=\rm{(\roman*)}]
\item\label{item:sigma-cocycle} The map $\sigma$ is a Hopf 2-cocycle for $\toba\#\Bbbk G$.
\item\label{item:sigma-nichols} For all $M\in\ydG$, we have $\left(\toba(M)\# \Bbbk G\right)_{\sigma}\simeq \toba\left(F_{\sigma}(M)\right)\# \Bbbk G$.
\end{enumerate}
\end{lemma}
\pf
\ref{item:sigma-cocycle} is clear. For \ref{item:sigma-nichols}, apply \cite[4.14 (a) \& (b)]{AAGMV}.
\epf

Here is the first notion towards the folding construction.

\begin{definition}\label{def:folding-data}
A \emph{folding datum} is a triple $(\sigma, M, \ub)$ where $\sigma \in Z^2(\Gamma,\Bbbk^\times)$, $M\in\ydg$ and $\ub\colon F_\sigma (M) \to M$ is an isomorphism in $\ydg$.
\footnote{In other words, the pair $(M, \ub)$ is an object in the category $\left( \ydg\right)^\sigma$ of $\sigma$-equivariant objects.} 
\end{definition}

We will mainly deal with folding data coming from the following source.

\begin{exa}\label{ex:folding-diag-type}
Fix $g_i\in\Gamma$, $\chi_i\in\widehat{\Gamma}$, $\sigma \in Z^2(\Gamma,\Bbbk^{\times})$. Let $f\colon \I \to \I$ be a permutation such that 
$$g_{f(i)}=g_i \qquad  \text{for all }i\in\I.$$
Consider $M=\oplus_{i\in \I} M(g_i,\chi_i)\in\ydg$, and let $0\ne x_i\in M(g_i,\chi_i)$. Note that $M$ is of diagonal type with braiding matrix $\bq=(q_{ij})_{i,j \in \I}$, $q_{ij}=\chi_j(g_i)$. 
Then the linear isomorphism $\ub:F_\sigma (M)\to M$, $x_i \mapsto x_{f(i)}$ is in $\ydg$ if and only if 
\begin{align*}
\chi_{f(j)}&=b_{\sigma}(-,g_j)\chi_j && \text{for all }j\in\I. 
\end{align*}
In this case we have that
\begin{align*}
q_{i,f(j)} &=b_\sigma (g_i,g_j) q_{ij} & & \text{ for all }i,j\in\I,
\end{align*}
and $f$ induces an automorphism of the Dynkin diagram of $\bq$ because $b_\sigma (g_i,g_i)=1$ and $b_\sigma (g_j,g_i)b_\sigma (g_i,g_j)=1$ for all $i\ne j\in\I$. 
\end{exa}

Fixed $M$, the folding data form a group with unit $(1,M,\id)$ and product 
$$(\sigma,M,\ub)\ast(\sigma',M,\ub')=(\sigma\sigma',M,\ub\circ F_\sigma(\ub')). $$

The next results are extracted from \cite[Part I]{Len12}. 

\begin{remark}\label{lemma:folding-data}
\begin{enumerate}[leftmargin=*,label=\rm{(\alph*)}]
\item\label{item:folding-data-a}  Let $(\sigma,M,\ub)$ be a folding datum, $H=\toba(M) \# \Bbbk\Gamma$. By Lemma \ref{lem:sigma}, $\ub$ induces a Hopf algebra isomorphism 
$\ub:H_{\sigma}\to H$.
\item\label{item:folding-data-b} The map $(1\otimes \ub)\Delta_{H_{\sigma}}:H{_\sigma}\to H{_\sigma}\times H$ makes $H{_\sigma}$ a right $H$-Galois object. Moreover ${_\sigma}H$ is an $(H,H)$-bi-Galois object.
\item\label{item:folding-data-c} Given two folding data $(\sigma,M,\ub)$, $(\sigma',M,\ub')$, the map
\begin{align*}
(\id\otimes \ub')\Delta_{H_{\sigma\sigma'}}: H_{\sigma\sigma'}\to H_{\sigma}\,\Box\, H_{\sigma'}
\end{align*}
is an isomorphism of bi-Galois objects.
\item The map in \ref{item:folding-data-c} determines a group homomorphism from the group of folding data over $M\in \ydg$ to the group of bi-Galois objects of $H=\toba(M) \# \Bbbk\Gamma$.
\end{enumerate}
\end{remark}

\subsection{Folding construction} \label{subec:folding-data}
Let $1\to\Sigma\to G\to \Gamma\to 1$ be a central extension of a finite abelian group $\Gamma$ by a finite abelian group $\Sigma$. Fix a set-theoretic section $s:\Gamma\to G$, and let 
\begin{align*}
\tau & \in Z^2(\Gamma,\Sigma), & &\tau(g,h)=s(g)s(h)s(gh)^{-1}, && g,h\in\Gamma.
\end{align*}
For each $t\in{\widehat{\Sigma}}$, we denote by $\sigma_t\in Z^2(\Gamma,\Bbbk^{\times})$ the $2$-cocycle $\sigma_t=t\circ \tau$. The assignment $t\mapsto \sigma_t$ is a group homomorphism $\widehat{\Sigma}\to Z^2(\Gamma,\Bbbk^\times )$.
Now we fix
\begin{itemize}[leftmargin=*]
\item a Yetter-Drinfeld module $M$ over $\Bbbk\Gamma$,
\item isomorphisms $\ub_t:F_{\sigma_t}(M)\to M$, $t\in\widehat{\Sigma}$, in $\ydg$ such that the map $t\mapsto (\sigma_t,M,\ub_t)$ is a group homomorphism from $\widehat{\Sigma}$ to the group of folding data for $M$. In particular, we have $\ub_0=\id:\, _{\sigma_0}H=H\to H$.
\end{itemize}

\begin{remark}
$M$ becomes a $\widehat{\Sigma}$-module, where $t\in \widehat{\Sigma}$ acts by the automorphism $\ub_t$. As $\Sigma$ is finite abelian, the  $\widehat{\Sigma}$-action diagonalizes and $M$ decomposes as a direct sum of $\Sigma$-eigenspaces:
\begin{align*}
M&=\bigoplus_{p\in\Sigma} M^{p}, &  M^{p}&=\{m\in M \colon \,  \ub_t(m)=t(p)m \text{ for all } t\in\widehat{\Sigma}\}.
\end{align*}
\end{remark}
\begin{theorem}{\cite[Theorem 3.6]{Len14}} \label{thm_NicholsFolding}
Let $\Gamma$, $G$ and $M$ as above. The following structure define a $\Bbbk G$-Yetter-Drinfeld module $\widetilde{M}$:
\begin{itemize}[leftmargin=*]
\item as a vector space, $\widetilde{M}=M$, 
\item the $G$-action is obtained by pulling back the $\Gamma$-action (hence $\Sigma$ acts trivially),
\item the $G$-grading is given by $\widetilde{M}_g =M_{\overline{g}}^{gs(g)^{-1}}:= M_{\overline{g}} \cap M^{gs(g)^{-1}}$, for each $g\in G$.
\end{itemize}
Also, as a braided vector space, $\widetilde{M}=M$.\qed
\end{theorem}

Next we introduce the folding construction, which produces a Nichols algebra over $G$ starting from folding data on a Nichols algebra over $\Gamma$.
The procedure gives a central extension of Hopf algebras and is related with the \emph{Fourier transform} developed in \cite{AG}.

\begin{theorem}\label{thm:folding}
Let $H:=\toba(M)\# \Bbbk \Gamma$.
There exists a Hopf algebra structure on $\tilde{H}:=\bigoplus\limits_{t\in\widehat{\Sigma}} H_{\sigma_{t}}$ given by
\begin{align*}
\Delta|_{{\sigma_t}H}&=\bigoplus_{t't''=t}(\id\otimes \ub_{t''})\Delta_H: H_{\sigma_t}\to \bigoplus_{t't''=t} H_{\sigma_{t'}}\,\otimes\,H_{\sigma_{t''}};
&&&
\epsilon_{|H}&=\epsilon_H, \qquad \epsilon_{|H_{\sigma_t}}=0, \, t\ne 0.
\end{align*}
Moreover, $\widetilde{H}\cong \toba(\tilde{M})\# \Bbbk G$ as Hopf algebras.
\end{theorem}

\begin{remark}
By the results above, the group of folding data induces a homomorphism of $2$-groups $\underline{\hat{\Sigma}}\to \underline{BiGal}(H)$. 
This in turn defines a homomorphism of $2$-groups $\underline{\hat{\Sigma}}\to \underline{BrPic}(\mathrm{Rep}(H))$ by \cite{Mombelli} and thus defines by \cite{EGNO} a $\Sigma$-extension of the tensor category $\mathrm{Rep}(H)$. This tensor category coincides with $\mathrm{Rep}(\widetilde{H})$. 

A way to see this fact comes from the equivariantization process applied to Hopf algebras since the folding data gives a functor from (the category defined from) $\Sigma$ to the Drinfeld double of $\mathrm{Rep}(H)$. Reciprocally, $\mathrm{Rep}(H)$ is the de-equivariantization of $\mathrm{Rep}(\widetilde{H})$ associated to a central extension of Hopf algebras \cite{AGP}, see also \cite[Theorem 3.6]{Len12}.
\end{remark}

\subsection{Folding data for trivial action of \texorpdfstring{$\kappa$}{}}\label{subsec:folding-examples}
Next we realize most of the examples of types $\alpha_{\theta}$, $\gamma_{\theta}$, $\delta_{\theta}$, $\epsilon_{\theta}$ and $\phi_4$ 
as foldings of braided vector spaces of diagonal type. In all cases we can proceed as in Example \ref{ex:folding-diag-type} with $\Sigma=\Z_2=\{e,\kappa\}$.

\begin{exa}\label{ex:folding-ADE}
Fix a finite Cartan matrix $\ba=(a_{ij})_{i,j\in\I}$ with simply-laced Dynkin diagram. Assume that $\Gamma$ is a finite abelian group generated by $g_i$, $i\in\I$, which admits a $2$-cocycle
\begin{align}\label{eq:2-cocycle-ADE}
\tau \in Z^2(\Gamma,\Sigma) \text { such that } &&
\tau(g_i,g_j)&=\begin{cases} \kappa^{a_{ij}}, & i<j\in\I; \\ e & i\ge j\in\I.
\end{cases}
\end{align}
Let $G$ be the extension of $\Gamma$ by $\Sigma$ associated to $\tau$. Thus $G$ is generated by $g_i$ and $\kappa$; in $G$ we have $g_ig_j=\kappa^{a_{ij}}g_jg_i$ for $i\ne j\in\I$.
Assume further that we have $\chi_i\in\widehat{\Gamma}$, $i\in\I$, satisfying
\begin{align*}
\chi_i(g_i)&=-1, & \chi_i(g_j)\chi_j(g_i)&=(-1)^{a_{ij}}, & &i\ne j\in\I.
\end{align*}
Then $V=\oplus_{i\in I} \Bbbk_{g_i}^{\chi_i}$ is of Cartan type $\ba$ with $q=-1$, as in Example \ref{ex:typeADE-1}. 

Here, $\widehat{\Sigma}=\{e,\mathtt{t}\}$, with $\mathtt{t}(\kappa)=-1$. Set $\sigma:=\mathtt{t}\circ \tau$, and
\begin{align*}
g_{i+\theta} &= g_i & \chi_{i+\theta} := b_{\sigma}(-, g_i)\chi_i & \in\widehat{\Gamma}, & &i\in\I, & M&:=\oplus_{i\in I_{2\theta}} \Bbbk_{g_i}^{\chi_i}.
\end{align*}
Then $M=V\oplus F_{\sigma}(V)$, and is of Cartan type with Cartan matrix $\widetilde{a}:=\left(\begin{smallmatrix} \ba & 0 \\ 0 & \ba \end{smallmatrix}\right)$.

Set also $f:\I_{2\theta}\to\I_{2\theta}$, $i\mapsto i+\theta$ modulo $2\theta$. Then $\ub:F_{\sigma}(M)\to M$ as in Example \ref{ex:folding-diag-type} is a folding datum, and
the map from $\widehat{\Sigma}$ to the group of folding data such that $\mathtt{t}\mapsto (\sigma,M,\ub)$ is a group homomorphim. Following \cite{Len14}, if $\ba$ is of type 
$X_{\theta}\in \{A_{\theta}|\theta\geq 2\} \cup \{D_{\theta}|\theta\geq 4\}\cup\{E_{\theta}|\theta=6,7,8\}$, we use ${^2}X_{\theta}^2$ to denote the corresponding folding of $X_{\theta}\times X_{\theta}$ by $f$, $\ub$ as above. 
\end{exa}

\begin{exa}\label{ex:folding-F4}
Fix a finite abelian group $\Gamma$ generated by $g_i$, $i\in\I_4$, which admits 
\begin{align}\label{eq:2-cocycle-F4}
\tau \in Z^2(\Gamma,\Sigma) \text { such that } && \tau(g_i,g_j)&=\kappa^{\delta_{i3}\delta_{j4}}, &&i,j\in\I_4.
\end{align}
Let $G$ be the extension of $\Gamma$ by $\Sigma$ associated to $\tau$. Now $G$ is generated by $g_i$ and $\kappa$; the relations 
$g_3g_4=\kappa g_4g_3$, and $g_ig_j=g_jg_i$ for $\{i,j\}\ne \{3,4\}$ hold in $G$.

Assume further that $\chi_i\in\widehat{\Gamma}$, $i\in\I$, satisfy
\begin{align*}
\chi_i(g_i)&=-1, & &i\in\I_4; & \chi_i(g_j)\chi_j(g_i)&=(-1)^{\delta_{i+1,j}}, & &i<j\in\I_4.
\end{align*}
Again, fix $\mathtt{t}\in \widehat{\Sigma}$ such that $\mathtt{t}(\kappa)=-1$. Set $\sigma:=\mathtt{t}\circ \tau$, and
\begin{align*}
g_{i+2} &= g_i & \chi_{i+2} := b_{\sigma}(-, g_i)\chi_i & \in\widehat{\Gamma}, & &i\in\{3,4\}, & M&:=\oplus_{i\in I_6} \Bbbk_{g_i}^{\chi_i}.
\end{align*}
Then $M$ is of Cartan type $E_6$. Let $f:\I_6\to\I_6$ be the bijection that exchanges $3\leftrightarrow 5$ and $4\leftrightarrow6$. Then $\ub:F_{\sigma}(M)\to M$ as in Example \ref{ex:folding-diag-type} is a folding datum, and the map $\mathtt{t}\mapsto (\sigma,M,\ub)$ is a group homomorphim from $\widehat{\Sigma}$ to the group of folding data. Following \cite{Len14}, ${^2}E_6$ denotes a folding as above. 
\end{exa}

\begin{exa}\label{ex:folding-C}
Fix a finite abelian group $\Gamma$ generated by $g_i$, $i\in\I$, which admits a $2$-cocycle 
\begin{align}\label{eq:2-cocycle-C}
\tau \in Z^2(\Gamma,\Sigma) \text { such that } && \tau(g_i,g_j)&=\begin{cases} \kappa, & j=i+1<\theta, \\ e, & \text{otherwise}.\end{cases}
\end{align}
Let $G$ be the extension of $\Gamma$ by $\Sigma$ associated to $\tau$. Thus $G$ is generated by $g_i$ and $\kappa$; in $G$, we have the relations $g_ig_{i+1}=\kappa g_{i+1}g_i$ if $i<\theta-1$, and $g_ig_j=g_jg_i$ otherwise.

Assume further that $\chi_i\in\widehat{\Gamma}$, $i\in\I$, satisfy
\begin{align*}
\chi_i(g_i)&=-1, & &i\in\I; & \chi_i(g_j)\chi_j(g_i)&=(-1)^{\delta_{i+1,j}}, & &i<j\in\I.
\end{align*}
Again, fix $\mathtt{t}\in \widehat{\Sigma}$ such that $\mathtt{t}(\kappa)=-1$. Set $\sigma:=\mathtt{t}\circ \tau$, and
\begin{align*}
g_{2\theta-i} &= g_i & \chi_{2\theta-i} := b_{\sigma}(-, g_i)\chi_i & \in\widehat{\Gamma}, & &i\in\I_{\theta-1}, & M&:=\oplus_{i\in I_{2\theta-1}} \Bbbk_{g_i}^{\chi_i}.
\end{align*}
Then $M$ is of Cartan type $A_{2\theta-1}$. Let $f:\I_{2\theta-1}\to\I_{2\theta-1}$, $f(i)=2\theta-i$. Then $\ub:F_{\sigma}(M)\to M$ as in Example \ref{ex:folding-diag-type} is a folding datum, and the map $\mathtt{t}\mapsto (\sigma,M,\ub)$ is a group homomorphim from $\widehat{\Sigma}$ to the group of folding data. Following \cite{Len14}, ${^2}A_{2n-1}$ denotes the folding above. 
\end{exa}

\begin{remark}
In the three examples above, we can take $\Gamma=\Z_2^{\I}$, see \cite[\S 5]{Len14}.
\end{remark}

\begin{theorem}\label{prop:folding-under-trivial-action}
Let $M$ as in \S \ref{subsubsec:type-ADE}, \S \ref{subsubsec:type-C} or \S \ref{subsubsec:type-F4}. Assume that $\kappa$ acts trivially on $M$.
\begin{enumerate}[leftmargin=*,label=\rm{(\alph*)}]
\item If $M$ is of type $\alpha_{\theta}$, $\delta_{\theta}$ or $\epsilon_{\theta}$, then $\toba(M)$ is a folded Nichols algebra as in Example \ref{ex:folding-ADE}.
\item If $M$ is of type $\phi_4$, then $\toba(M)$ is a folded Nichols algebra as in Example \ref{ex:folding-F4}.
\item If $M$ is of type $\gamma_{\theta}$, then $\toba(M)$ is a folded Nichols algebra as in Example \ref{ex:folding-C}.
\end{enumerate}
\end{theorem}

\pf
Since $\kappa$ acts trivially, this follows by \cite[Theorems 5.6, 5.7, 5.8]{Len14}.
\epf

\begin{remark}\label{rem:diagonalizeBraiding}
Fix a braided vector space $M$ of type either $\alpha_{\theta}$, $\theta\ge 4$, $\delta_{\theta}$, $\theta\ge 5$, or $\epsilon_{\theta}$, $\theta=6,7,8$ (type $\gamma_{\theta}$ was already considered in general), with Cartan matrix $\ba=(a_{ij})_{i,j\in\I_\theta}$. By Lemma \ref{lem:trivial-action-kappa-big-rank}, $\chi_i(\kappa)=1$ for all $i\in\I$, so $M$ is of diagonal type. We exhibit a basis in which the braiding \emph{is} of diagonal type, and we give the braiding matrix. 
\begin{enumerate}[leftmargin=*,label=\rm{(\Roman*)}]
\item\label{item:braiding-trivial-action-I} Set $q_{ii}=-1$ for all $i\in\I_\theta$, and $q_{ij}=\chi_j(g_i)$ if $a_{ij}=0$.

\item\label{item:braiding-trivial-action-II} Let $1\leq i<j\leq \theta$ be such that $a_{ij}=-1$. Let $q_{ij}\in\Bbbk^{\times}$ be such that $q_{ij}^2=\chi_j(g_i^2)$, and set $q_{ji}:=-q_{ij}^{-1}$. By Step \ref{lem:alpha3-chi2-square} of Proposition \ref{prop:alpha3-Nichols}, if $k>i$ also satisfies that $a_{ik}=-1$ we may choose $q_{ik}=q_{ij}$. 

\item\label{item:braiding-trivial-action-III} We also set $\bq=(q_{ij})_{i,j\in\I_{2\theta}}$, where
\begin{align*}
q_{ij} &= \begin{cases}
-1, & i\le \theta<j \text{ or }j\le \theta<i; \\
-q_{i-\theta,j-\theta}, & i,j>\theta, a_{i-\theta,j-\theta}=-1; \\
q_{i-\theta,j-\theta}, & i,j>\theta, a_{i-\theta,j-\theta}=0; \\
-1, & i=j>\theta.
\end{cases}
\end{align*}
\end{enumerate}

Given $i\in\I_\theta$, there is $j\neq i$ in $\I_\theta$ such that $a_{ij}\ne 0$. If possible, take $j>i$ such that $a_{ij}\ne 0$; otherwise take $j<i$ with $a_{ij}\ne 0$. By \ref{item:braiding-trivial-action-II} above, we can define
\begin{align*}
\xt_i &:= x_i+q_{ij}x_{\un{i}}, & 
\xt_{\un{i}} &:= x_i-q_{ij}x_{\un{i}}, &
i\in\I.
\end{align*}
Using \eqref{eq:ADE-self-braiding}, \eqref{eq:ADE-connected-braiding} and \eqref{eq:ADE-disconnected-braiding} we verify that
\begin{align*}
c(\xt_i \ot \xt_j) &= \begin{cases}
-1, & i=j; \\
q_{ij}, & i\ne j.
\end{cases}
&
c(\xt_{\un{i}} \ot \xt_{\un{j}}) &= \begin{cases}
-1, & i=j \\
-q_{ij}, &  a_{ij}=-1; \\
q_{ij}, &  a_{ij}=0.
\end{cases}
\\
c(\xt_i \ot \xt_{\un{j}}) &= \begin{cases}
-1, & i=j \\
-q_{ij}, &  a_{ij}=-1; \\
q_{ij}, &  a_{ij}=0.
\end{cases}
&
c(\xt_{\un{i}} \ot \xt_j) &= \begin{cases}
-1, & i=j \\
q_{ij}, &  a_{ij}=-1; \\
q_{ij}, &  a_{ij}=0.
\end{cases}
\end{align*}
so the braiding matrix of $M$ is $\bq$, and $M$ is, respectively, of type $A_{\theta} \times A_{\theta}$, $D_{\theta} \times D_{\theta}$ or $E_{\theta} \times E_{\theta}$, both copies with parameter $q=-1$.
\end{remark}

\subsection{Trivializing the action of \texorpdfstring{$\kappa$}{} via a twist}\label{subsec:trivial-action}

Retain the notation introduced in \S \ref{subsec:folding-examples}. Thus $G$ is a non-abelian group and  $M=\oplus_{i\in\I} M_i \in \ydG$, where $M_i=M(g_i^G,\chi_i)$.

For the cases not covered by Lemma \ref{lem:trivial-action-kappa-big-rank}, we will show the existence of a 2-cocycle $\sigma$ such that $\kappa$ acts trivially on $F_\sigma(M)$. Recall the parity vector $\mathtt{P}=(\chi_1(\kappa),\ldots,\chi_\theta(\kappa))$ from \eqref{eq:action-kappa}.

\smallbreak
Let $\sigma\in H^2(G,\Bbbk^\times)$. Following \S\ref{subsec:categorical-action},
the \emph{twisted Yetter-Drinfeld module} associated to $\sigma$ is $F_\sigma(M)=M^{\sigma}=\oplus_{i\in\I} M(g_i,\chi_i^{\sigma})$, where $\chi_i^{\sigma}\in\widehat{G^{g_i}}$ is given by the following formula:
\begin{align*}
\chi_i^{\sigma}(h)&=\sigma(hg_ih^{-1},h)\sigma^{-1}(h,g_i)\chi_i(h), & & h\in G.
\end{align*}
Since $\kappa$ is central in $G$, we have 
$\chi_i^\sigma(\kappa)= \sigma(g_i,\kappa)\sigma^{-1}(\kappa,g_i)\chi_i(\kappa)$.

\begin{prop}\label{prop:DoiTwist}
Let $G$ be a non-abelian group, $M\in\ydG$ of type $\alpha_2$, $\alpha_3$, $\delta_4$, $\gamma_3$, $\gamma_4$ or $\phi_4$ such that $\supp M$ generates $G$.
There exists $\sigma\in H^2(G,\Bbbk^\times)$ with $\chi_i^{\sigma}(\kappa)=1$ for all $i$.
\end{prop}

Let us outline the strategy that will be used in the Appendix to prove this statement.
\begin{enumerate}[leftmargin=*,label=\rm{(\roman*)}]
\item We go through the cases and list the possible $\mathtt{P}=(\chi_i(\kappa))_i \in\{\pm 1\}^n$. 

\item It is sufficient to consider one $\mathtt{P}$ representing each Weyl groupoid orbit. The $i$-th reflection of  $M$ is
$$\rho_i M= \oplus_{j\in\I} M(g_jg_i^{-c_{ij}},\chi_j\chi_i^{-c_{ij}}),$$
and the parity vector of $\rho_i M$ is $\mathtt{P}'=(\chi_j(\kappa)\chi_i(\kappa)^{-c_{ij}})_{j\in\I}$.

\item Next we introduce an auxiliary \emph{minimal} group $G^{\min}$. Namely $G^{\min} \subseteq \End M$ is generated by (the action of) $g_i$. The definition of $G^{\min}$ depends only on the scalars $\chi_i(g_j)$, $\chi_i(\kappa)$. It is enough to prove Proposition \ref{prop:DoiTwist} for this group, since the asserted $2$-cocycle $\sigma$ on $G^{\min}$ can be pulled back to $G$.

\item  In the next steps we show case-wise that there exists $\sigma\in H^2(G^{\min},\Bbbk^\times)$ such that 
\begin{align*}
\sigma(g_i,\kappa)\sigma(\kappa,g_i)^{-1} &= \chi_i(\kappa), && \text{for all }i\in\I.
\end{align*}

\item For type $\alpha_2$ there are two Weyl groupoid orbits for $\mathtt{P}=(\chi_1(\kappa),\chi_2(\kappa))$, namely $\{(1,1)\}$ and $\{(1,-1),(-1,-1),(-1,1)\}$, the first one corresponding to $\mathtt{P}$ trivial.  For $(-1,1)$ we get three different types of groups according to the order of $\chi_2(g_1^2)$: we find the desired cocycle using semidirect product decompositions.

\item  The cases $\alpha_3$ and $\gamma_3$ are treated using spectral sequences arguments for a central extension. Necessary information about the structure of the group 
(minimal orders of central elements for example) enters conveniently via the existence of a $1$-dimensional representation, constructed from the structure of $M$.

\item The remaining cases are treated using two simultaneous extensions.
\end{enumerate}

We postpone the proof until \S \ref{subsubssection:proof-lemma-Doi-twist} since we need technical results on group cohomology. 

\smallbreak
By assumption $G$ is a central extension of an abelian group $\Gamma$ by $\langle\kappa\rangle \simeq \Z_2$, say:
$$1\longrightarrow\Z_2\longrightarrow G\longrightarrow \Gamma\longrightarrow 1. $$
To illustrate the proof of Proposition \ref{prop:DoiTwist}, we give an example where such central extensions are related to symplectic forms on $\Z_2^\theta$.

\begin{exa} 
Let $\Gamma=\Z_2^\theta$ with generators $g_i$, $i\in\I$. The commutator in $G$ defines a symplectic form on $\Gamma$ such that the radical $\Z_2^r$ is the image of the center $Z(G)$; in particular $g_i,g_j$ commute iff the symplectic form on them is zero; hence the size of the conjugacy class of $g_i$ is $1$ or $2$ depending if $g_i$ is in the radical or not. This symplectic form of type $(\theta,r)$ is uniquely determined.

Central extensions $\Z_2 \to G \to \Gamma$ are classified by quadratic forms with fixed symplectic form. For type $\alpha_2$ we have $(\theta,r)=(2,0)$ and there are two types of central extensions of order $8$ with this commutator structure, namely the extraspecial groups $2_-^{2+1}$ (the quaternion group) and $2_+^{2+1}$ (the dihedral group). The group $2_-^{2+1}$ has defining relations $g_1^2=g_2^2=\kappa$ and trivial second cohomology $H^2(2_-^{2+1},\Bbbk^\times)$, while the group $2_+^{2+1}$ has defining relations $g_1^2=e$, $g_2^2=\kappa$ (depending on a choice of generators) and non-trivial second cohomology  $H^2(2_+^{2+1},\Bbbk^\times)=\{1,\sigma\}$. Our proof works for the group $2_-^{2+1}$ because the relations  \eqref{eq:typeADE-aij=-1} defining $M$ imply 
$$\chi_1(\kappa)=\chi_1(g_1^2)=1, \qquad \chi_2(\kappa)=\chi_2(g_2^2)=1.$$
For $2_+^{2+1}$, both $\mathtt{P}=(1,1)$ and $\mathtt{P}=(-1,1)$ are possible and we have a $2$-cocycle $\sigma$ with 
$$\sigma(g_1,\kappa)\sigma^{-1}(\kappa,g_1)=-1, \qquad \sigma(g_2,\kappa)\sigma^{-1}(\kappa,g_2)=1.$$
The other choices of generators for $2^{2+1}_+$ work similarly and also follow from the first choice by using Weyl groupoid reflections.

The cases $\alpha_3$, $\delta_4$, $\gamma_3$, $\gamma_4$, $\phi_4$ present similar behavior because $\frac{1}{2}(\theta-r)=1$ in all of them. For each case there is an underlying extraspecial group $2_{\pm}^{2+1}$. 
\end{exa}

\begin{remark}
The proof becomes more involved for an arbitrary abelian group $\Gamma$ because there are many central extensions by $\Z_2$, parametrized by the powers of the generators, and the existence of a non-trivial group cohomology is very sensitive to these choices. 
\end{remark}

We are ready to state the main result of this section, which states that each Nichols algebras of type $\alpha_{\theta}$, $\gamma_{\theta}$, $\delta_{\theta}$, $\epsilon_{\theta}$, or $\phi_{4}$ is a twist of the corresponding Nichols algebra of diagonal type as in \S \ref{subsec:folding-examples}. 

\begin{theorem}\label{thm:Doi-Twist}
Let $G$ be a finite non-abelian group, $M\in\ydG$ of type either $\alpha_{\theta}$, $\gamma_{\theta}$, $\delta_{\theta}$, $\epsilon_{\theta}$ or $\phi_{4}$ whose support generates $G$. 
Then there exists $\sigma\in H^2(G,\Bbbk)$ such that $F_{\sigma}(M)$ is of diagonal type.
\end{theorem}

\pf
If $M$ is of type $\alpha_\theta$, $\theta\ge 4$, $\gamma_\theta$, $\theta\ge 5$, $\delta_\theta$, $\theta\ge 5$ or $\epsilon_\theta$, $\theta=6,7,8$, then $\kappa$ acts trivially by Lemma \ref{lem:trivial-action-kappa-big-rank}, so $M$ is naturally a Yetter-Drinfeld module over $\Gamma=G/\langle \kappa\rangle$, a finite abelian group, thus $M$ is itself of diagonal type. For the other cases, we apply Proposition \ref{prop:DoiTwist}.
\epf

Let $G$ and $M=\oplus_{i\in\I} M_i \in \ydG$ as above. Consider
\begin{align}\label{eq:ell}
\ell:=\begin{cases}
2 & \text{for type }\phi_4,
\\
\theta-1 & \text{for type }\gamma_{\theta},
\\
\theta & \text{otherwise}.
\end{cases}
\end{align}

Keeping the notation used for \S \ref{subsubsec:type-ADE}, \ref{subsubsec:type-C}, \ref{subsubsec:type-F4}, we have
\begin{itemize}[leftmargin=*]\renewcommand{\labelitemi}{$\circ$}
\item If $i\le \ell$, then $\dim M_i=2$, $g_i^G=\{g_i,\kappa g_i\}$, where $\kappa\in Z(G)$ satisfies $\kappa^2=1$. We fix a basis $x_i$, $x_{\un{i}}$ of $M_i$ as above.
\item If $i>\ell$, then $\dim M_i=1$, $g_i^G=\{g_i\}$. For a basis, we fix any nonzero element $x_i$ in $M_i$.
\end{itemize}

\begin{remark}\label{rem:support-quotient}
Let $\mathtt{G}$ be the group generated by $g_i$, $i\in \I$ and $\kappa$ with relations
\begin{align*}
g_i\kappa &=\kappa g_i, & \kappa^2 &= e, & g_ig_j &= \begin{cases}
\kappa g_jg_i, & i\ne j \le \ell, \, a_{ij}=-1, \\ g_jg_i & \text{otherwise},
\end{cases}
\end{align*}
where $\ba=(a_{ij})$ denotes the Cartan matrix of $M$. Then $\mathtt{G}$ is a central extension of $\Z^{\I}$ by $\Z_2$, and the braided vector space $M$ has a realization over $\mathtt{G}$. Moreover, 
the subgroup of $G$ generated by $\supp M$ is a quotient of $\mathtt{G}$.

The $\Z^{\I}$-grading on $T(M)$ and its homogeneous quotients (in particular, $\toba(M)$) is given by the induced coaction of $\Bbbk\Z^{\I} \simeq \Bbbk \mathtt{G}/\langle \kappa \rangle$.
\end{remark}

\section{Generation in degree one}\label{sec:gen-degree-one-trivial-action}
Using Theorem \ref{thm:Doi-Twist} and generation-in-degree-one for the diagonal setting \cite{A-presentation}, we get
\begin{theorem}\label{thm:gen-degree-one}
Let $H$ be a finite-dimensional pointed Hopf algebra with infinitesimal braiding $M$ of type $\alpha_{\theta}$, $\gamma_{\theta}$, $\delta_{\theta}$, $\epsilon_{\theta}$ or $\phi_{4}$. 
Then 
$$\gr H\simeq \toba(M)\#\Bbbk G(H).$$
In other words, $H$ is generated by skew-primitive and group-like elements.
\end{theorem}

\pf
Let $R$ be the diagram of $H$, so $\gr H\simeq R\#\Bbbk G(H)$. Put also $\toba:=R^{\ast}$. Notice that $W:=M^{\ast}$ is of the same type as $M$; we fix $g_i$, $i\in\I$, $\kappa$ as in \S \ref{subsec:Nichols-non-abelian-groups}. Let $G$ be the subgroup of $G(H)$ generated by $g_i$, $i\in\I$. Then $\toba\in\ydG$ is a finite-dimensional pre-Nichols algebra of $W$; i.e. $\toba=\oplus_{n \ge 0}\toba^n$, where $\toba^n=R(n)^{\ast}$, is a graded Hopf algebra such that 
$\toba^0=\Bbbk 1$ and is generated as an algebra by $W=\toba^1$.

For $\sigma\in H^2(G,\Bbbk^\times)$ as in Theorem \ref{thm:Doi-Twist}, set $\sigma:\toba\#\Bbbk G \otimes \toba\#\Bbbk G \to \Bbbk$ as in \eqref{eq:cocycle-extension-to-bosonization}, and let $\Ht:=(\toba\#\Bbbk G)_{\sigma}$, $\Ht^n:=\toba^n\# \Bbbk G$. Then $\Ht=\oplus_{n \ge 0}\Ht^n$ is a graded coalgebra, because twisting by $\sigma$ leaves the coalgebra structure unchanged. Thus $\Ht$ is pointed with coradical $\Ht^0 \simeq \Bbbk G$ by \cite[5.3.4]{Mo-libro}. As $\sigma$ is trivial in degree $>0$, $\Ht$ is a graded Hopf algebra: By \cite[4.14 (a)]{AAGMV}, $\Ht\simeq \toba'\# \Bbbk G$, where $\toba'\in \ydG$ is a pre-Nichols algebra of $F_{\sigma}(W)$. 

As $\kappa\in Z(G)$ and it acts trivially on $W^{\sigma}$, we have that $\kappa\in Z(H)$. Set $Q:=\Ht/\Ht(\kappa-1)$, $\Gamma=G/\langle \kappa\rangle$. The $G$-actions on $F_{\sigma}(W)$ and on $\toba'$ induce respective $\Gamma$-actions on them. Also, $F_{\sigma}(W)$ and $\toba'$ become $\Bbbk\Gamma$-comodules via $\pi:G\twoheadrightarrow \Gamma$. Moreover, with these structures both $\toba'$ and $F_{\sigma}(W)$ are $ \ydg$, and $\toba'\in \ydg$ is a pre-Nichols algebra of $F_{\sigma}(W)$. We identify $\Ht \simeq\toba'\# \Bbbk G$ and consider the map
$\Phi:\Ht\to \toba'\#\Bbbk\Gamma$, $\Phi(x\# g)=x\#\pi(g)$. Then $\Phi$ is a surjective Hopf algebra map such that $\Phi(\kappa-1)=0$, hence $\Phi$ induces a surjective Hopf algebra map $\phi:Q\to \toba'\#\Bbbk\Gamma$. As $\Ht$ is a central extension of $Q$ by $\Bbbk\Z_2$ (since $\kappa^2=1$, $\kappa\ne 1$),
\begin{align*}
\dim Q &= \frac{1}{2}\dim \Ht = \frac{1}{2}\dim\toba' |G| = \dim\toba' |\Gamma|,
\end{align*}
so $\phi$ is an isomorphism.

Now $\Gamma$ is an abelian group since for each $i\ne j$ either $g_ig_j=g_jg_i$ or $g_ig_j=\kappa g_jg_i$ in $G$, and $\Gamma$ is generated by the images of the $g_i$'s (see \S \ref{subsec:Nichols-non-abelian-groups}). Hence \cite{A-presentation} implies that $\toba'=\toba(F_{\sigma}(W))$, so $\toba=\toba(W)$ by \cite[4.14 (b)]{AAGMV}. Dualizing, $R=\toba(M)$, as desired.
\epf

\section{Generators and relations for Nichols algebras}\label{sec:generators-relations}
In this section we exhibit a presentation by generators and relations for the Nichols algebras of the Yetter-Drinfeld modules $M=\oplus_{i\in\I} M_i\in \ydG$ as in \S \ref{subsubsec:type-ADE}, \ref{subsubsec:type-C}, \ref{subsubsec:type-F4}, which are standard, i.e. all reflection $\rho_i M$ are of the same type as $M$. In particular, the root system $\varDelta^M$ is a classical one. We may assume that $G$ is the group in Remark \ref{rem:support-quotient}.

\subsection{Types \texorpdfstring{$\alpha_2$}{} and \texorpdfstring{$\alpha_3$}{}}\label{subsec:presentation-type-ADE}

We give a presentation of Nichols algebras of types $\alpha_2$ and $\alpha_3$. They will be a key step toward the presentation for the general case, since all relations which are not powers of root vectors are supported on smaller submodules of these types.

\subsubsection{Type \texorpdfstring{$\alpha_2$}{}}\label{subsubsec:alpha2-presentation}
Let $G$ be a group, $e\ne \kappa$, $g_1$, $g_2\in G$ such that 
\begin{align*}
g_1g_2&=\kappa g_2g_1, & \kappa^2&=e, & g_i^G&=\{g_i,\kappa g_i\}, & i&=1,2. 
\end{align*}
Following \cite{HS-rank2-1}, the subgroup of $G$ generated by $g_1$ and $g_2$ is a quotient of
\begin{align}\label{eq:alpha2-universal-group} 
G_2=\langle g_1, g_2, \kappa \vert \kappa^2=1, \kappa g_1=g_1\kappa, \kappa g_2=g_2\kappa, g_2g_1=\kappa g_1g_2\rangle.
\end{align}
Assume there are $\chi_i\in \widehat{G^{g_i}}$ such that $\chi_i(g_i)=-1$ and $\chi_1(\kappa g_2^2)\chi_2(\kappa g_1^2)=1$.
For $i\in\I_2$ set $M_i=M(g_i^G,\chi_i)\in \ydG$, thus $M=M_1\oplus M_2$ is of type $\alpha_2$ by \cite[Theorem 4.6]{HS-rank2-1}.

\begin{prop}\label{prop:alpha2-Nichols}
Let $M=M_1\oplus M_2\in \ydG$ of type $\alpha_2$ as above. Then $\toba(M)$ is presented by generators $x_1, x_{\un{1}}, x_2, x_{\un{2}}$ and relations
\begin{align}
\label{eq:alpha2-NicholsMi-relations}
&x_i^2=x_{\un{i}}^2=0, && (\ad_c x_i)x_{\un{i}}=0, && i \in \I.
\\
\label{eq:alpha2-M12-dim2-relations}
&x_{1\un{2}}= -\chi_2(g_1^2) x_{\un{1}2}, & 
&x_{\un{1}\un{2}}= -\chi_1(\kappa) x_{12};
\\  \label{eq:alpha2-NicholsM12-relations} 
&x_{12}^2=0, & & [x_{12}, x_{\un{1}2} ]_c=0.
\end{align}

The following set is a PBW basis of $\toba(M)$: 
\begin{align} \label{eq:alpha2-NicholsPBW}
\{x_2^a x_{\un{2}}^b x_{12}^c x_{\un{1}2}^d x_1^e x_{\un{1}}^f \colon a, b, c, d, e, f \in \{0,1\}\}.
\end{align}
\end{prop}

\pf We proceed in several steps.

\begin{step}
By \eqref{eq:ADE-self-braiding} the Nichols algebra of $M_i$ for $i=1,2$ is a quantum linear space, and \eqref{eq:alpha2-nongenuinePBW} implies that relations \eqref{eq:alpha2-NicholsMi-relations} hold in $\toba (M)$.
\end{step}

\begin{step}
The inclusion $M_{12}:=(\ad_c M_1) M_2\hookrightarrow \toba (M)$ extends to a $\Z^2$-graded algebra inclusion $\toba(M_{12})\hookrightarrow \toba (M)$ and the multiplication
\begin{align} \label{eq:alpha2-nongenuinePBW}
\toba(M_2)\ot \toba (M_{12}) \ot \toba (M_1) \longrightarrow \toba (M)
\end{align}
is an isomorphism of $\Z^2$-graded objects in $\ydG$, where $M_1$ sits in degree $\alpha_1$, $M_2$ in degree $\alpha_2$, and $M_{12}$ in degree $\alpha_1 + \alpha_2$.
\end{step}

This follows by \cite[Theorem 4.6]{HS-rank2-1}. In order to find defining relations for $\toba(M_{12}) \subset \toba (M)$, we need a more explicit description of the structure of $M_{12}$.

\begin{step}\label{lem:alpha2-structureM12}
\begin{enumerate}[leftmargin=*,label=\rm{(\alph*)}]
\item \label{item:alpha2-M12-basis} The set $\{x_{12}, x_{\un{1}2}\}$ is a basis of $M_{12}$, and the braiding in this basis is of diagonal type with matrix
\begin{align}\label{eq:alpha2-M12-braiding}
\begin{pmatrix} -1 & -\chi_1(\kappa)\chi_2(\kappa) \\ -\chi_1(\kappa)\chi_2(\kappa) & -1  \end{pmatrix}.
\end{align}

\item \label{item:alpha2-M12-dim2} Relations \eqref{eq:alpha2-M12-dim2-relations} and \eqref{eq:alpha2-NicholsM12-relations} hold in $\toba (M)$.
\end{enumerate}
\end{step}

\pf
Note that $\partial_{x_2} (x_{12})=x_1$ and $\partial_{x_2} (x_{\un{1}2})= x_{\un{1}}$, thus $x_{12}$ and $x_{\un{1}2}$ are linearly independent in $\toba(M)$. Relations \eqref{eq:alpha2-M12-dim2-relations} are verified using the skew-derivations of $T(M_1\oplus M_2)$. The braiding of $M_{12}$ is obtained from a straight-forward computation in $\ydG$.
Now \eqref{eq:alpha2-NicholsM12-relations} follows from \eqref{eq:alpha2-nongenuinePBW} and \eqref{eq:alpha2-M12-braiding}.
\epf

Note that $\toba(M_{12})$ is presented by the relations \eqref{eq:alpha2-NicholsM12-relations} and $(x_{\un{1}2})^2=0$, which has not been included above since it can be deduced from the previous ones, as we show in Step \ref{lem:alpha2-structureM112-M221}.

With \eqref{eq:alpha2-nongenuinePBW} in mind, the next step towards exhibiting a presentation of $\toba(M)$ should be to find braided commutations between $\toba(M_i)$ and $\toba(M_{12})$ for $i=1,2$. Such relations are known to exist, since $(\ad M_i)^2 M_j = 0$ for $i\ne j$ by \cite[Lemma 4.2]{HS-rank2-1}. However, we show next that these can be deduced from some of the already stablished relations.

\begin{step}\label{lem:alpha2-structureM112-M221}
Let $A$ denote the quotient of $T(M_1 \oplus M_2)$ by the ideal generated by \eqref{eq:alpha2-NicholsMi-relations} and \eqref{eq:alpha2-M12-dim2-relations}. In $A$ the following relations hold
\begin{align} \label{eq:alpha2-structureM112}
(\ad_c x_1) x_{12} &=0, & ( \ad_c x_{\un{1}} )x_{12} &= 0, &
( \ad_c x_1) x_{\un{1}2} &=0, & (\ad_c x_{\un{1}})x_{\un{1}2} &= 0;
\\
\label{eq:alpha2-structureM221}
 [x_{12}, x_2]_c&=0, & [x_{12}, x_{\un{2}}]_c&=0, &
[x_{\un{1}2}, x_2]_c&=0, & [x_{\un{1}2}, x_{\un{2}}]_c&=0; \\
\label{eq:alpha2-square-generators-M12}
 x_{\un{1}2}^2 &= 0.
\end{align}
\end{step}

\pf
First, $(\ad_c x_1) x_{12}=(\ad_c x_1^2) x_2=0$ since  $x_1^2=0$.
Analogously $(\ad_c x_{\un{1}}) x_{\un{1}2}=0$ follows from $x_{\un{1}}^2=0$.
Next, using that $x_{12}$ is a scalar multiple of $x_{\un{1}\un{2}}$, we get that $( \ad_c x_{\un{1}} )x_{12} = 0$ follows from $x_{\un{1}}^2=0$. Since $x_{\un{1}2}$ is a scalar multiple of $x_{1\un{2}}$ we get $(\ad_c x_1)x_{\un{1}2} = 0$ from $x_1^2=0$. Thus \eqref{eq:alpha2-structureM112} hold.

For \eqref{eq:alpha2-structureM221}, unpacking the definitions we see that $[x_{12}, x_2]_c=[x_{1\un{2}}, x_{\un{2}}]_c=0$ follow from $x_2^2=x_{\un{2}}^2=0$. 
Now $[x_{12}, x_{\un{2}}]_c=0$ follows using the previous argument, since $x_{12}$ is a scalar multiple of $x_{\un{1}\un{2}}$. The remaining relation holds similarly.

Finally, we show that \eqref{eq:alpha2-square-generators-M12} follows from \eqref{eq:alpha2-structureM112}, \eqref{eq:alpha2-structureM221} and \eqref{eq:alpha2-M12-dim2-relations}. In fact,
\begin{align}\label{eq:x12square}
\begin{aligned}
(x_{12})^2
&= (x_1 x_2 - x_{\un{2}} x_1) (\ad_c x_1) (x_2) \\
&=-(\chi_2(\kappa g_1^2))^{-1} x_1  (\ad_c x_1) (x_2) x_{\un{2}} - \chi_2(g_1^2) x_{\un{2}} (\ad_c x_{\un{1}}) (x_2) x_1 \\
&= -\chi_2(\kappa) (\ad_c x_{\un{1}}) (x_2) x_1 x_{\un{2}} + \chi_2(\kappa g_1^2) (\ad_c x_{\un{1}}) (x_2) x_2 x_1 \\ 
&= -\chi_2(\kappa) (\ad_c x_{\un{1}}) (x_2) (\ad_c x_1) (x_{\un{2}}) = \chi_2(\kappa g_1^2) (x_{\un{1}2})^2,
\end{aligned}
\end{align}
as claimed.
\epf

We are ready to give a presentation of $\toba(M_1\oplus M_2)$.
Let $R$ denote the quotient of $T(M_1\oplus M_2)$ by the ideal generated by \eqref{eq:alpha2-NicholsMi-relations}, \eqref{eq:alpha2-M12-dim2-relations} and \eqref{eq:alpha2-NicholsM12-relations}. We already know from \eqref{eq:alpha2-nongenuinePBW} and Step \ref{lem:alpha2-structureM12} that the canonical projection $T(M_1\oplus M_2)\to \toba(M)$ factors to a surjective algebra map $R\to \toba(M)$. We show that this map is injective. Since $\dim \toba(M)=2^6$, it is enough to verify that the set  \eqref{eq:alpha2-NicholsPBW} linearly generates $R$. Let $J$ denote the subspace spanned by \eqref{eq:alpha2-NicholsPBW} in $R$. Since $J$ contains $1$, it is enough to show that $J$ is a left ideal, which reduces to verify that $x_i I \subset I$ and $x_{\un{i}} I \subset I$ for $i=1,2$. Clearly $x_2 J \subset J$; as $x_{\un{2}}x_2$ is a scalar multiple of $x_2x_{\un{2}}$, it is equaly clear that $x_{\un{2}} J \subset J$. So we need to verify that $x_1 J \subset J$ and $x_{\un{1}} J \subset J$, which follow since the \eqref{eq:alpha2-structureM112} and \eqref{eq:alpha2-structureM221} hold in $R$ by Step \ref{lem:alpha2-structureM112-M221}.
\epf

\subsubsection{Type $\alpha_3$}\label{subsubsec:alpha3-presentation}
Let $G$ denote a non-abelian group and let $M=M_1 \oplus M_2 \oplus M_3$ in $\ydG$ of type $\alpha_3$. By \cite[Lemma 5.2]{HV-rank>2} for $i\in\I_2$, the subgroup  $\langle \kappa, g_i, g_{i+1} \rangle \subset G$ is a quotient of $G_2$, see \eqref{eq:alpha2-universal-group}. Next we describe $\toba(M)$.

\begin{prop}\label{prop:alpha3-Nichols}
For $M$ of type $\alpha_3$, the Nichols algebra $\toba(M)$ is presented by generators $x_i, x_{\un{i}}$, $i\in\I_3$ and relations
\begin{align}\label{eq:alpha3-NicholsMi-relations}
&x_i^2=x_{\un{i}}^2=0, & &(\ad_c x_i)x_{\un{i}}=0, && i \in \I;
\\\label{eq:alpha3-Mij-dim2-relations}
& x_{i\un{j}}=-\chi_j(g_i^2) x_{\un{i}j}, & & x_{\un{i}\un{j}}=-\chi_i(\kappa) x_{ij},  &&i<j,  \, a_{ij}=-1;
\\
&x_{ij}^2=0, & & [x_{ij}, x_{i\un{j}} ]_c=0, &&i<j,  \, a_{ij}=-1;
\\\label{eq:alpha3-NicholsM13-relations}
&x_{13}=x_{1\un{3}}=0, & &x_{\un{1}3}=x_{\un{1}\un{3}} = 0;
\\\label{eq:alpha3-NicholsM123-relations} 
&x_{123}^2=0, & &[x_{123}, x_{\un{1}23}]_c=0;
\\\label{eq:alpha3-structureM2123}
&(\ad_c x_2) x_{123}=0, & &(\ad_c x_2) x_{\un{1}23}=0.
\end{align}
A PBW basis of $\toba(M)$ is given by
\begin{align} \label{eq:alpha3-NicholsPBW}
\{x_3^{a_3} x_{\un{3}}^{b_3} x_{23}^{a_{23}} x_{\un{2}3}^{b_{23}} x_2^{a_{2}} x_{\un{2}}^{b_{2}}  x_{123}^{a_{123}} x_{\un{1}23}^{b_{123}} x_{12}^{a_{12}} x_{\un{1}2}^{b_{12}} x_1^{a_{1}} x_{\un{1}}^{b_{1}} \colon \ a_{\beta}, b_{\beta} \in \{0,1\} \}. 
\end{align}
\end{prop}

\pf Again, we proceed in several steps.

\begin{stepi}\label{step:alpha3-nongenuinePBW}
The multiplication map is an isomorphism of $\Z^3$-graded objects in $\ydG$:
\begin{align*} 
\toba(M_3) \ot \toba (M_{23}) \ot\toba(M_2) \ot \toba (M_{123}) \ot \toba (M_{12}) \ot \toba (M_1) \simeq \toba(M)
\end{align*}
\end{stepi}
This follows by \cite[Theorem 2.6]{HS-rank2-1}.
Next we give some relations that hold in $\toba(M)$.

\begin{stepi}\label{step:alpha3-rels-subalgebras}
The relations \eqref{eq:alpha3-NicholsMi-relations}, \eqref{eq:alpha3-Mij-dim2-relations}
and \eqref{eq:alpha3-NicholsM13-relations} hold in $\toba(M)$.
\end{stepi}

We know that $M_{13}=0$, so \eqref{eq:alpha3-NicholsM13-relations} hold, while \eqref{eq:alpha3-NicholsMi-relations} and \eqref{eq:alpha3-Mij-dim2-relations} follow since $M_1\oplus M_2$ and $M_2\oplus M_3$ are of type $\alpha_2$. \qed

Following the treatment in \S \ref{subsubsec:alpha2-presentation}, we describe $\toba (M_{123})\subset \toba(M)$. As a first step in this direction, we seek for a basis of $M_{123}$.

\begin{stepi}\label{lem:alpha3-M123-dim2}
Let $A$ denote the quotient of $T(M_1 \oplus M_2\oplus M_3)$ by the ideal generated by  \eqref{eq:alpha3-NicholsM13-relations} and \eqref{eq:alpha3-Mij-dim2-relations}. In $A$ the following relations hold
\begin{align} \label{eq:alpha3-M123-dim2-relations}
x_{1\un{2}3}&=-\chi_2(g_1^2) x_{\un{1}23}, &
x_{\un{1}\un{2}3}&=-\chi_1(\kappa)x_{123};
\end{align}
\end{stepi}
\pf
We only verify the first one. We compute
\begin{align*}
x_{1\un{2}3} &=[x_1,  [x_{\un{2}},x_3]_c]_c=[[x_1,  x_{\un{2}}]_c,x_3]_c + (g_1\cdot x_{\un{2}}) [x_1,x_3]_c - \chi_3(\kappa) [x_1, x_{\un{3}}] x_{\un{2}}\\
&=[[x_1,  x_{\un{2}}]_c, x_3]_c = - \chi_2(g_1^2) [[x_{\un{1}},  x_2]_c, x_3]_c = - \chi_2(g_1^2) [x_{\un{1}},  [x_2,x_3]_c]_c,
\end{align*}
where the second equality follows by \eqref{eq:braided-commutator-iteration}, the third from  \eqref{eq:alpha3-NicholsM13-relations}, the fourth from \eqref{eq:alpha3-Mij-dim2-relations}, and the fifth one by \eqref{eq:braided-commutator-iteration} and \eqref{eq:alpha3-NicholsM13-relations}.
\epf

Surprisingly, there are further restrictions on the character $\chi_2$:
\begin{stepi}\label{lem:alpha3-chi2-square}
If $M_1\oplus M_2\oplus M_3$ is of type $\alpha_3$, then $\chi_2(g_1^2)=\chi_2(g_3^2)$.
\end{stepi}
\pf
We compute the action of $g_2$ on $x_{123} \in \toba(M)$ following two different approaches. Applying \eqref{eq:group-action-in-adjoint} first, followed by \eqref{eq:alpha3-Mij-dim2-relations} we get
\begin{align*}
g_2\cdot x_{123}& = - (\ad_c x_{\un{1}}) (\ad_c x_2)x_{\un{3}} = \chi_3(g_2^2)  x_{\un{1}\un{2}3}\\
&= \chi_3(g_2^2)(x_{\un{1}} ((\ad_c x_{\un{2}})x_3) -  \chi_2(\kappa g_1^2) \chi_3(\kappa g_1) ((\ad_c x_2)x_3) x_{\un{1}})
\end{align*}
On the other hand, if we first unpack the definition of $\ad_c x_1$ and then let $g_2$ act, we get
\begin{align*}
g_2 \cdot  x_{123} &= g_2\cdot (x_1 ((\ad_c x_2)x_3) - \chi_3(g_1) ((\ad_c x_{\un{2}})x_3) x_1) \\
&= - x_{\un{1}} ((\ad_c x_2)x_{\un{3}}) + \chi_2(\kappa)  \chi_3(g_1) ((\ad_c x_{\un{2}})x_{\un{3}}) x_{\un{1}} \\
&= \chi_3(g_2^2) x_{\un{1}} ((\ad_c x_{\un{2}})x_3) -  \chi_3(g_1) ((\ad_c x_2)x_3) x_{\un{1}} \\
&= \chi_3(g_2^2) (x_{\un{1}} ((\ad_c x_{\un{2}})x_3) -  \chi_2(\kappa g_3^2)\chi_3(\kappa g_1) ((\ad_c x_2)x_3) x_{\un{1}})
\end{align*}
where the second equality follows by \eqref{eq:group-action-in-adjoint}, the third one from \eqref{eq:alpha3-Mij-dim2-relations} and the last one from $ \chi_2(\kappa g_3^2)\chi_3(\kappa g_2^2)=1$.

These two equations give $(\chi_2(g_1^2) - \chi_2(g_3^2) ) ((\ad_c x_2)x_3) x_{\un{1}} =0$, and the claim follows because $\partial_{x_{\un{1}}}( ((\ad_c x_2)x_3) x_{\un{1}}) = (\ad_c x_2)x_3 \ne 0$.
\epf

The next result is analogue to Step \ref{lem:alpha2-structureM12} in Proposition \ref{prop:alpha2-Nichols}, and its proof follows from similar arguments, so we only give a sketch. 

\begin{stepi}\label{lem:alpha3-structureM123}
\begin{enumerate}[leftmargin=*,label=\rm{(\alph*)}]
\item \label{item:alpha3-M123-basis}  A basis of $M_{123}$ is $\{x_{123}, x_{\un{1}23}\}$, where the braid is diagonal with matrix
\begin{align}\label{eq:alpha3-M123-braiding}
\begin{pmatrix} -1 & -\chi_1(\kappa)\chi_2(\kappa)\chi_3(\kappa) \\ -\chi_1(\kappa)\chi_2(\kappa)\chi_3(\kappa) & -1  \end{pmatrix}.
\end{align}

\item \label{item:alpha3-M123-dim2} Relations \eqref{eq:alpha3-NicholsM123-relations} hold in $\toba (M)$. 
\end{enumerate}
\end{stepi}

\pf
\ref{item:alpha3-M123-basis} The set $\{x_{123}, x_{\un{1}23}\}$ linearly spans $M_{123}$ by Step \ref{lem:alpha3-M123-dim2}, and it is linearly independent since so are $\partial_{x_3}(x_{123})=x_{12}$ and $\partial_{x_3}(x_{\un{1}23})=x_{\un{1}2}$. The braiding is computed using \eqref{eq:group-action-in-adjoint} and Step \ref{lem:alpha3-chi2-square}. Now \ref{item:alpha3-M123-dim2} follows from \ref{item:alpha3-M123-basis} and Step \ref{step:alpha3-nongenuinePBW}.
\epf

The Nichols algebra of $M_{123}$ is presented by the relations \eqref{eq:alpha3-NicholsM123-relations} and $x_{\un{1}23}^2=0$; we will show that this last relation can be deduced from others. 

Braided commutations between $M_i$ and $M_{123}$ for $i=1,3$ can now be deduced:

\begin{stepi}\label{lem:alpha3-structureM1123-M1233}
Let $A$ denote the quotient of $T(M_1 \oplus M_2 \oplus M_3)$ by the ideal generated by \eqref{eq:alpha2-NicholsMi-relations}, \eqref{eq:alpha3-NicholsM13-relations} and \eqref{eq:alpha3-Mij-dim2-relations}. In $A$ the following relations hold
\begin{align} \label{eq:alpha3-structureM1123}
(\ad_c x_1) x_{123} &=0, & (\ad_c x_{\un{1}}) x_{123} &=0, &
(\ad_c x_1) x_{\un{1}23}& =0, & (\ad_c x_{\un{1}}) x_{\un{1}23} &=0;
 \\
\label{eq:alpha3-structureM1233}
[x_{123}, x_3]_c&=0, & [x_{123}, x_{\un{3}}]_c&=0, &
[x_{\un{1}23} , x_3]_c&=0, & [x_{\un{1}23}, x_{\un{3}}]_c&=0;
\end{align}\vspace{-.8cm}
\begin{align}
\label{eq:alpha3-structureM2123-repeated}
(\ad_c x_{\un{2}}) &x_{123}=-\chi_1(g_2^{-2})(\ad_c x_2) x_{\un{1}23}, &
(\ad_c x_{\un{2}}) x_{\un{1}23}&=-\chi_1(g_2^{-2})(\ad_c x_2) x_{123}.
\end{align}
\end{stepi}

\pf
The relations \eqref{eq:alpha3-structureM1123} can be verified using the argument in the proof of Step \ref{lem:alpha2-structureM112-M221} of Proposition \ref{prop:alpha2-Nichols}. For the first relation in \eqref{eq:alpha3-structureM1233}, use \eqref{eq:braided-commutator-iteration} to get
\begin{align*}
[x_{123}, &x_3]_c=[[x_1,x_{23}]_c, x_3]_c= [x_1,[x_{23}, x_3]_c]_c - x_{\un{2}} [x_1, x_3]_c +[x_1, x_{\un{3}}]_c x_2 = 0.
\end{align*}
The three remaining relations follow similarly. For \eqref{eq:alpha3-structureM2123-repeated}, use \eqref{eq:braided-commutator-iteration}, \eqref{eq:alpha3-Mij-dim2-relations} and \eqref{eq:alpha2-structureM112} to get
\begin{align*}
(\ad_c x_{\un{2}}) x_{123}  &= \left[x_{\un{2}1}, x_{23}\right]+(g_2\kappa \cdot x_1)[x_{\un{2}},x_{23}]_c - \chi_3(g_1) [x_{\un{2}},x_{\un{2}3} ]_c x_1
\\
& = -\chi_1(g_2^{-2})\left[x_{2\un{1}}, x_{23}\right] = -\chi_1(g_2^{-2}) (\ad_c x_2) x_{\un{1}23}.
\end{align*}
and the other relation follows analogously.
\epf

We employ skew derivations to verify braided commutations between $M_2$ and $M_{123}$. 

\begin{stepi}\label{lem:alpha3-structureM2123} 
Relations \eqref{eq:alpha3-structureM2123} hold in $\toba(M)$.
\end{stepi}

\pf
We focus on the first relation. One can see directly from \eqref{eq:group-action-in-adjoint} that 
\begin{align*}
\partial_{x_i}(M_{123}) = \partial_{x_{\un{i}}}(M_{123}) &= 0, &\partial_{x_3}(x_{123}) &= x_{12}, & \partial_{x_{\un{3}}}(x_{\un{1}2\un{3}}) &= x_{\un{1}2}.
\end{align*}
Now, using \eqref{eq:group-action-in-adjoint} and also Step \ref{lem:alpha3-chi2-square} we get 
\begin{align*}
\partial_{x_{\un{3}}}(x_{123}) &= -\chi_2(\kappa)x_{1\un{2}}, &
\partial_{x_3}(x_{\un{1}2\un{3}}) &= - \chi_3 (g_2^2) x_{\un{1}\un{2}}.
\end{align*}
For $i\in\I_2$, both $\partial_{x_i}$ and $\partial_{x_{\un{i}}}$ annihilate the first relation, and we compute 
\begin{align*}
\partial_{x_3}((\ad_c x_2)x_{123}) &= x_2 x_{12}  -  \chi_3 (g_2^2) x_{\un{1}\un{2}} x_{\un{2}}
= x_2 x_{12}  + \chi_3 (\kappa g_2^2) x_{12} x_{\un{2}} = \chi_3 (\kappa g_2^2) [x_{12},x_{\un{2}}]_c =0,
\end{align*}
where the third equality follows by \eqref{eq:alpha3-Mij-dim2-relations}. Similarly, 
$\partial_{x_{\un{3}}}((\ad_c x_2)x_{123}) = \chi_2(\kappa) [x_{\un{1}2},x_{\un{2}}]_c =0$.
The other relations follow analogously.
\epf

Next we deduce that braided brackets between $\toba(M_{23})$ and $\toba(M_{12})$ can be rewritten in terms of intermediate factors of the decomposition given in Step \ref{step:alpha3-nongenuinePBW}.

\begin{stepi}\label{lem:alpha3-structure[M12,M23]}
We have $[\toba(M_{12}),\toba(M_{23})]_c \subset \toba(M_2) \ot \toba(M_{123})$ in $\toba(M)$.
\end{stepi}
\pf
This follows from $x_2^2=0$, $x_{\un{2}}^2=0$, \eqref{eq:alpha3-Mij-dim2-relations} and \eqref{eq:alpha3-structureM2123}. As an illustration:
\begin{align*}
\left[ x_{12}, x_{23}\right]_c = (\ad_c x_1)(\ad_c x_2^2) x_3& - x_{\un{2}}  x_{123}  -x_{12\un{3}} x_2,
\end{align*}
which belongs to $\Bbbk x_{\un{2}} x_{123}   + \Bbbk x_2 x_{\un{1}23}$.
\epf

\begin{stepi}\label{lem:alpha3-square-generators-M123}
Let $A$ denote the quotient of $T(M_1\oplus M_2 \oplus M_3)$ by the ideal generated by \eqref{eq:alpha3-Mij-dim2-relations}, \eqref{eq:alpha3-M123-dim2-relations},  \eqref{eq:alpha3-structureM1123}, \eqref{eq:alpha3-structureM1233}, and \eqref{eq:alpha3-structureM2123}. Then in $A$ we have
\begin{align} \label{eq:alpha3-square-generators-M123}
x_{123}^2=  \chi_2(g_1^4)\chi_3(g_2^2) x_{\un{1}23}^2 
\end{align}
\end{stepi}

\begin{proof}
Use \eqref{eq:alpha3-structureM1123}, \eqref{eq:alpha3-structureM1233}, and \eqref{eq:alpha3-structureM2123} several times to get explicit braided commutations between  $x_{123}$ and each of its monomials:
\begin{align*}
x_{123}^2 &= (x_1x_2x_3 - x_1 x_{\un{3}} x_2 - \chi_3(g_1) x_{\un{2}} x_3 x_1 + \chi_3(\kappa g_1) x_{\un{3}} x_{\un{2}} x_1)x_{123} \\
&=\chi_2(g_1^2) x_{\un{1}23}  (x_1 (x_2 x_{\un{3}} - \chi_3(g_2^2) x_3x_2) - \chi_3(\kappa g_1) (x_{\un{2}} x_{\un{3}}  - \chi_3(\kappa g_2^2) x_3 x_{\un{2}} )x_1 )\\
&=\chi_2(g_1^4)\chi_3(g_2^2) x_{\un{1}23}^2,
\end{align*}
where the last equality follows from \eqref{eq:alpha3-Mij-dim2-relations} and \eqref{eq:alpha3-M123-dim2-relations}.
\end{proof}

As in the proof of Proposition \ref{prop:alpha2-Nichols}, by Steps \ref{step:alpha3-rels-subalgebras}, \ref{lem:alpha3-structureM123} and \ref{lem:alpha3-structureM2123}, there exists an algebra surjection from $\toba$, the algebra presented by relations in Proposition \ref{prop:alpha3-Nichols}, onto $\toba(M)$. Now we use Steps \ref{step:alpha3-nongenuinePBW}, \ref{lem:alpha3-M123-dim2}, \ref{lem:alpha3-structureM1123-M1233}, \ref{lem:alpha3-structure[M12,M23]}, \ref{lem:alpha3-square-generators-M123} together with the fact that $\dim \toba(M)=2^{12}$ to conclude that $\toba=\toba(M)$.
\epf

\subsection{A kind of distinguished pre-Nichols algebra}\label{subsec:distinguished}
For Nichols algebras of diagonal type, a presentation by generators and relations was achieved by the first author in \cite{A-presentation}. A fundamental role is played by an intermediate quotient known as the distinguished pre-Nichols algebra. 
Inspired by that construction, we define a pre-Nichols algebra $\htoba(M)$ for each one of the braidings $M$ described in \S \ref{subsubsec:type-ADE}, \ref{subsubsec:type-C}, \ref{subsubsec:type-F4}.
More precisely, our construction resembles the algebra $\htoba_{\bq}$ introduced in Remark \ref{rem:diagonal-pre-Nichols}. 

This pre-Nichols algebra will also play a key role in \S \ref{sec:liftings}, where we describe the liftings.
\medspace

Given $M=\oplus_{i \in \I}M_i$ as in \S \ref{subsubsec:type-ADE}, \ref{subsubsec:type-C}, \ref{subsubsec:type-F4}, let $\sigma$ as in Theorem \ref{thm:Doi-Twist}, and denote by $\bq$ the braiding matrix of $W:=F_{\sigma}(M)$ (see \S \ref{subsec:categorical-action}). Recall the index $\ell$ defined in \eqref{eq:ell}. Consider
\begin{align*}
\un{\,\cdot\,}:&\I_{\ell}\to\I_{\theta+\ell}, & 
\un{i} &:=\begin{cases}
i+\theta, & \text{for types }\alpha_{\theta},\delta_{\theta},\epsilon_{\theta}, \\ \theta+\ell-i+1, & \text{otherwise,}
\end{cases}
\end{align*}
so $\I_{\theta+\ell}=\I_{\theta}\cup\{\un{i}:i\in\I_{\ell}\}$. We fix a basis $\xt_i$, $i\in\I_{\theta+\ell}$, such that
\begin{itemize}[leftmargin=*]
\item $\xt_i$, $\xt_{\un{i}}$ is a basis of $F_{\sigma}(M_i)$ for each $i\in\I_{\ell}$;
\item $\xt_i=x_i$ if $\ell<i\le \theta$;
\item the braiding in this basis is given by $\bq$; i.e. $c(\xt_i\ot \xt_j)=q_{ij} \xt_j\ot \xt_i$.
\end{itemize}
Let $\mathtt{A}=(\att_{ij})_{i,j\in \I_{\theta+\ell}}$ be the Cartan matrix of $\bq$.

\begin{remark}\label{rem:grading-Fsigma-diagonal}
Let $\Xi:\Z^{\theta+\ell}\to\Z^{\theta}$ be the group morphism such that
\begin{align*}
\Xi(\alpha_i) &=\alpha_i, \quad i\le \theta; & \Xi(\alpha_i) &=\alpha_{j}, \quad i>\theta, i=\un{j}.
\end{align*}
This map \emph{identifies} the two elements of the basis above corresponding to each $F_{\sigma}(M_i)$ when $\dim M_i=2$.
Thus, if $\toba$ is an $\N_0^{\theta}$-graded pre-Nichols algebra such that $F_{\sigma}(\toba)$ is $\N_0^{\theta+\ell}$-graded 
(for the usual grading as pre-Nichols algebra of diagonal type), then 
\begin{align*}
\dim \toba_{\beta} &= \sum_{\gamma\in\Xi^{-1}(\beta)} \dim F_{\sigma}(\toba)_{\gamma}, & &\beta\in\N_0^{\theta}.
\end{align*}
Thus the Hilbert series $\mH_{\toba}$ is the image of $\mH_{F_{\sigma}(\toba)}$ under $\Xi$.
\end{remark}

Let $\htoba(M)$ be the algebra generated $x_i$, $i\in \I$, $x_{\un{j}}$, $j\in\I_{\ell}$, subject to the relations
\begin{align}\label{eq:dist-pre-NicholsMi-relations}
&(\ad_c x_i)x_{\un{i}}, \quad x_i^2, \quad x_{\un{i}}^2, & & i \in \I_{\ell};
\\
\label{eq:dist-pre-NicholsMi-dim1}
&x_i^2, & &i >\ell;
\\
\label{eq:dist-pre-Nichols-aij=-1-relations}
& x_{i\un{j}}+\chi_j(g_i^2) x_{\un{i}j}, \quad x_{\un{i}\un{j}}+\chi_i(\kappa) x_{ij},  &&i<j\le \ell,  \, a_{ij}=-1;
\\
\label{eq:dist-pre-Nichols-aij=0-relations}
&x_{ij}, \quad x_{i\un{j}}, \quad x_{\un{i}j}, \quad x_{\un{i}\un{j}}, && i<j\le \ell, \, a_{ij}=0;
\\
\label{eq:dist-pre-Nichols-aij=0-relations-bis}
& x_{ij}, \quad x_{\un{i}j}, && i\le \ell<j, \, a_{ij}=0;
\\
\label{eq:dist-pre-Nichols-rels-gamma}
& (\ad_c x_j) x_{ij\un{i}}=0, && a_{ij}=-2.
\\
\label{eq:dist-pre-Nichols-structureMjijk}
&(\ad_c x_j) x_{ijk}, \quad (\ad_c x_j) x_{\un{i}jk}, && i<j<k, \, a_{ji}=a_{jk}=-1.
\end{align}

\begin{prop}\label{prop:pre-Nichols}
$\htoba(M)$ is a graded pre-Nichols algebra with Hilbert series
\begin{align}\label{eq:hilbert-series-pre-Nichols}
\mH_{\htoba(M)} = & \prod_{\substack{\beta \in \Delta_+^{M}-\{\alpha_i\} \\ \dim M_{\beta}=2}} \left(\frac{1}{1-t^{\beta}}\right)^2 \, 
\prod_{j\in\I_{\ell}}(1+t_i)^2 \prod_{\substack{\beta \in \Delta_+^{M}-\{\alpha_i\} \\ \dim M_{\beta}=1 }} \left(\frac{1}{1-t^{\beta}}\right) \, \prod_{j>\ell}(1+t_j).
\end{align}
\end{prop}
\pf
First, $\htoba(M)$ is a quotient of $T(M)$ by an homogeneous ideal, so it is a graded algebra. 
Also, \eqref{eq:dist-pre-NicholsMi-relations}--\eqref{eq:dist-pre-Nichols-aij=0-relations-bis} are primitive elements of $T(M)$, 
while those in \eqref{eq:dist-pre-Nichols-structureMjijk} are primitive modulo the previous relations, see the proofs of Propositions \ref{prop:alpha2-Nichols} and \ref{prop:alpha3-Nichols}. Thus $\htoba(M)$ is a graded pre-Nichols algebra.

We fix $G=\mathtt{G}$ defined in Remark \ref{rem:support-quotient}, and $M$ has a canonical $G$-Yetter-Drinfeld module structure.
Now $\htoba:=F_{\sigma}\left(\htoba(M)\right)\in\ydG$ is a pre-Nichols algebra of $W$ such that 
$$\left(\htoba(M)\#\Bbbk G \right)_{\sigma} \simeq \htoba \#\Bbbk G .$$
We claim that $\htoba=\htoba_{\bq}$. If so, the statement on the Hilbert series follows from Remarks \ref{rem:diagonal-pre-Nichols} and \ref{rem:grading-Fsigma-diagonal}.
First, note that \eqref{eq:typeADE-1-rels-distinguished} holds in $\htoba$ for all $i\in\I_{\theta+\ell}$ by \eqref{eq:dist-pre-NicholsMi-relations} and \eqref{eq:dist-pre-NicholsMi-dim1}.

To verify \eqref{eq:typeADE-1-rels-3}, let $i<j$ be such that $\att_{ij}=0$. 
We consider 5 cases. 

\begin{itemize}[leftmargin=*]\renewcommand\labelitemi{$\diamond$}
\item $i\le \ell$, $j=\un{i}$. The space of primitive elements of $T(M)$ of degree $2\alpha_i\in\N_0^{\theta}$ is three-dimensional, spanned by $x_i^2$, $x_{\un{i}}^2$ and $x_{i\un{i}}$. Thus $\htoba(M)$ has no primitive elements of degree $2\alpha_i$. 
On the other hand, the spaces of primitive elements of $T(W)$ of degree $\alpha_i+\alpha_j\in\N_0^{\theta+\ell}$ is one-dimensional spanned by $x_{ij}$, and those of degrees $2\alpha_i$, $2\alpha_j$ are also one-dimensional, spanned by $\xt_{i}^2$ and $\xt_{j}^2$ respectively. As the space of homogeneous primitive elements of $\htoba$ coincides with that of $\htoba(M)$, we have that $\xt_{ij}=0$ in $\htoba$.

\item $i\le\ell$, $\theta<j$, $\att_{ik}=0$, where $j=\un{k}$. The space of primitive elements of $T(M)$ of degree $\alpha_i+\alpha_k\in\N_0^{\theta}$ has dimension $4$, spanned by $x_{ik}$, $x_{i\un{k}}$, $x_{\un{i}k}$, $x_{\un{i}\un{k}}$, so the space of primitive elements of $\htoba(M)$ of the same degree is 0.
On the other hand, the space of primitive elements of $T(W)$ of degree $\alpha_i+\alpha_k$ is spanned by $\xt_{ik}$, $\xt_{ij}$, $\xt_{k\un{i}}$, $\xt_{j\un{i}}$. As the space of homogeneous primitive elements of $\htoba$ coincides with that of $\htoba(M)$, we have that $\xt_{ij}=0=\xt_{j\un{i}}$ in $\htoba$. This also shows that $\xt_{ik}=\xt_{k\un{i}}=0$.

\item $i\le\ell$, $\theta<j$, $\att_{ik}=-1$, where $j=\un{k}$. The space of primitive elements of $T(M)$ of degree $\alpha_i+\alpha_k\in\N_0^{\theta}$ is two-dimensional, spanned by \eqref{eq:dist-pre-Nichols-aij=-1-relations}, so $\htoba(M)$ has no primitive elements of this degree.
On the other hand, the spaces of primitive elements of $T(W)$ of degrees $\alpha_i+\alpha_j$ and $\alpha_{\un{i}}+\alpha_k$ are one-dimensional, 
spanned by $\xt_{ij}$ and $\xt_{k\un{i}}$ respectively, and those of degrees $\alpha_i+\alpha_k$ and $\alpha_j+\alpha_{\un{i}}$ are $0$ since $a_{ik}=a_{\un{i}j}=-1$. Hence $\xt_{ij}=\xt_{k\un{i}}=0$ in $\htoba$.

\item $i\le\ell<j\le \theta$. Here $a_{ij}=0$, and the space of primitive elements of $T(M)$ of degree $\alpha_i+\alpha_j\in\N_0^{\theta}$ is two-dimensional, spanned by  $x_{ij}$ and $x_{i\un{j}}$, so the space of primitive elements of $\htoba(M)$ of this degree is 0. 
On the other hand, the spaces of primitive elements of $T(W)$ of degrees $\alpha_i+\alpha_j$ and $\alpha_{\un{i}}+\alpha_j$ are one-dimensional, spanned by $\xt_{ij}$ and $\xt_{j\un{i}}$, respectively. Hence $\xt_{ij}=\xt_{j\un{i}}=0$ in $\htoba$.

\item $i,j\le\ell$. Here $a_{ij}=0$, respectively $a_{\un{i}\un{j}}=0$, and this case is  the second one.
\end{itemize}

Now we check \eqref{eq:typeADE-1-rels-2} in $T(W)$, modulo \eqref{eq:typeADE-1-rels-1} and \eqref{eq:typeADE-1-rels-3}. We have 3 cases:

\begin{itemize}[leftmargin=*]\renewcommand\labelitemi{$\heartsuit$}
\item $j\le \ell$. The space of primitive elements of $T(M)$, modulo \eqref{eq:dist-pre-NicholsMi-relations}, \eqref{eq:dist-pre-NicholsMi-dim1},
\eqref{eq:dist-pre-Nichols-aij=-1-relations}, \eqref{eq:dist-pre-Nichols-aij=0-relations} and \eqref{eq:dist-pre-Nichols-aij=0-relations-bis}, of degree $\alpha_i+2\alpha_j+\alpha_k$ has dimension $\le 2$ and is spanned by \eqref{eq:dist-pre-Nichols-structureMjijk}: we can use skew-derivations as in type $\alpha_3$, see the proof of Proposition \ref{prop:alpha3-Nichols}. 
Thus the space of primitive elements of $\htoba(M)$ of this degree is 0. 
On the other hand, the space of primitive elements of $T(W)$, modulo \eqref{eq:typeADE-1-rels-1} and \eqref{eq:typeADE-1-rels-3}, of each degree in $\Xi^{-1}(\alpha_i+2\alpha_j+\alpha_k)$ is either one-dimensional or $0$. Indeed, the non-zero cases are spanned by
\begin{itemize}
\item  $[\xt_{ijk},\xt_j]_c$, $[\xt_{\un{i}\un{j}\un{k}},\xt_{\un{j}}]_c$ if $j<\ell$;
\item $[\xt_{\theta-2 \, \theta-1 \, \theta},\xt_{\theta-1}]_c$, $[\xt_{\theta+2 \, \theta+1 \, \theta},\xt_{\theta+1}]_c$ if $j=\theta-1$ in type $\gamma_\theta$;
\item $[\xt_{123},\xt_2]_c$, $[\xt_{653},\xt_{5}]_c$ if $j=2$ in type $\phi_4$.
\end{itemize}
As the space of homogeneous primitive elements of $\htoba$ coincides with that of $\htoba(M)$, we deduce that \eqref{eq:typeADE-1-rels-2} hold in $\htoba$.

\item $j>\ell$. There are three possibilities: $i=\theta-1$, $j=\theta$, $k=\un{i}$ in types $\gamma_{\theta}$ or $\phi_4$, and $i\in\{2,5\}$, $j=3$, $k=4$ in type $\phi_4$. The proof is analogous, using 
\eqref{eq:dist-pre-Nichols-rels-gamma} for the first case, and \eqref{eq:dist-pre-Nichols-structureMjijk} for the last one.
\end{itemize}
From the analysis above there exists a surjective Hopf algebra map $\htoba\twoheadrightarrow\htoba_{\bq}$. In a similar way, checking spaces of homogeneous primitive elements of appropiate degree, each defining relation of $\htoba(M)$ annihilates in $F_{\sigma}^{-1}(\htoba_{\bq})=F_{\sigma^{-1}}(\htoba_{\bq})$, so there exists a surjective Hopf algebra map
$F_{\sigma}^{-1}(\htoba_{\bq}) \twoheadrightarrow \htoba(M)$. As $F_{\sigma}$ preserves the $\N_0$-graduation of the pre-Nichols algebras, both surjective maps are indeed isomorphisms. Hence, $\htoba=\htoba_{\bq}$ as we claimed.
\epf

\subsection{The subalgebra of coinvariants}\label{subsec:coinvariants}

Let $\pi:\htoba(M)\twoheadrightarrow\toba(M)$ be the canonical projection, $\Zc(M):=\htoba(M)^{\co \pi}$ the subalgebra of coinvariants.
The next step towards the presentation of $\toba(M)$ is to describe $\Zc(M)$. To uncover the structure of  this subalgebra, we will use a cocycle $\sigma$ as in Theorem \ref{thm:Doi-Twist} to translate the problem to the diagonal setting, where the situation is better understood. In particular, we compute the Hilbert series of $\Zc(M)$. Since we know that of $\htoba(M)$, we will thus obtain the Hilbert series of $\toba(M)$.

To do so, and also to compute a PBW basis of $\toba(M)$ later on, we fix a reduced expression of the element $w_0$ of maximal length (or equivalently, a convex order on $\Delta_+$) for each type. Using this reduced expression, \cite{HS-book} defines a submodule $M_{\beta} \in\ydG$, $\beta\in\Delta_+$.
We exhibit a basis $\{x_{\beta}\}$ or $\{x_{\beta},x_{\un{\beta}}\}$ of the submodule $M_{\beta}$, depending on its dimension.

\begin{description}
\item[$\alpha_{\theta}$] The set of positive roots is $\varDelta_+=\{\alpha_{ij}: i\le j\in\I\}$, and
\begin{align*}
\alpha_1 <\alpha_{12} <\alpha_2 <\cdots <\alpha_{\theta-1} < \alpha_{1\theta} < \alpha_{2\theta} < \cdots < \alpha_{\theta}
\end{align*}
is a convex order on $\varDelta_+$. By \cite{HS-book} the modules $M_{\beta}$ can be defined as
\begin{align}\label{eq:M-alphaij-defn}
M_{\alpha_{ij}} &= (\ad_c M_i)\cdots (\ad_c M_{j-1})M_j.
\end{align}
A basis of $M_{\alpha_{ij}}$ is given by
\begin{align}\label{eq:M-alphaij-basis}
x_{\alpha_{ij}} &= x_{i \, i+1 \cdots j}, & x_{\un{\alpha_{ij}}}&=x_{\un{i} \, i+1 \cdots j}.
\end{align}

\item[$\delta_{\theta}$] The positive roots are $\varDelta_+=\{\alpha_{ij}: i\le j\in\I\}-\{\alpha_{\theta-1\theta}\}\cup \{\alpha_{i\theta-2}+\alpha_{\theta}: i\in\I_{\theta-2}\} \cup \{\alpha_{i\theta-2}+\alpha_{j\theta}: i<j\in\I_{\theta-2}\}$, and a convex order on $\varDelta_+$ is given by
\begin{align*}
\alpha_1 & < \alpha_{12} <\alpha_2 <\cdots <\alpha_{\theta-2} < \alpha_{1\theta-1} < \cdots < \alpha_{\theta-1} < \cdots < \alpha_{\theta-1} 
\\
& < \alpha_{1\theta}+\alpha_{2\theta-2} < \cdots < \alpha_{1\theta}+\alpha_{\theta-2}
< \alpha_{2\theta}+\alpha_{3\theta-2} < \cdots < \alpha_{2\theta}+\alpha_{\theta-2}
\\
&  < \cdots < \alpha_{\theta-3\theta}+\alpha_{\theta-2} < \alpha_{1\theta} < \cdots 
< \alpha_{\theta-2 \, \theta} < \alpha_{\theta}.
\end{align*}
The Yetter-Drinfeld modules $M_{\alpha_{ij}}$, $j\ne \theta$, are defined as for type $\alpha_{\theta}$. For $j=\theta$ let
$$ M_{\alpha_{i\theta}} = (\ad_c M_i)\cdots (\ad_c M_{\theta-3})(\ad_c M_{\theta-1})(\ad_c M_{\theta-2})M_{\theta}. $$
For the other roots we have
\begin{align*}
M_{\alpha_{i\theta-2}+\alpha_{\theta}} &= (\ad_c M_i)\cdots (\ad_c M_{\theta-2})M_{\theta},
\\
M_{\alpha_{i\theta}+\alpha_{j\theta-2}} &= [M_{\alpha_{i\theta-2}+\alpha_{\theta}}, 
M_{\alpha_{j\theta-1}}]_c.
\end{align*}
A basis of $M_{\beta}$ for either $\beta=\alpha_{ij}$ or $\beta=\alpha_{i\theta-2}+\alpha_{\theta}$ is given as in \eqref{eq:M-alphaij-basis}. For $\beta=\alpha_{i\theta}+\alpha_{j\theta-2}$, a basis of $M_{\beta}$ is
\begin{align}\label{eq:M-alpha-deltan-basis}
x_{\beta} &= [x_{i \cdots \theta-2 \theta}, x_{j \cdots \theta-1}]_c, & x_{\un{\beta}}&=[x_{\un{i} i+1 \cdots \theta-2 \theta}, x_{j \cdots \theta-1}]_c.
\end{align}

\item[$\epsilon_{\theta}$] Here one can fix a convex order as in \cite[\S 5]{A-standard}. For braidings of diagonal type, a PBW basis is obtained recursively on the height of the roots, starting with $\xt_{\alpha_i}=\xt_i$ for simple roots, and later $\xt_{\beta}=[\xt_{\beta_1},\xt_{\beta_2}]_c$ for some pair $(\beta_1,\beta_2)$ such that $\beta_1+\beta_2=\beta$, see \cite[Corollary 3.17]{A-jems}. For each non simple root $\beta\in\varDelta_+$ we have, accordingly,
\begin{align*}
M_{\beta} &=[M_{\beta_1},M_{\beta_2}]_c, & 
x_{\beta} &=[x_{\beta_1},x_{\beta_2}]_c, & 
x_{\un{\beta}} &=[x_{\un{\beta_1}},x_{\beta_2}]_c.
\end{align*}

\item[$\gamma_{\theta}$] Now $\varDelta_+=\{\alpha_{ij}: i\le j\in\I\}\cup \{\alpha_{i\theta}+\alpha_{k \, \theta-1}: i\le j\in\I_{\theta-1}\}$, and
\begin{align*}
\alpha_1 & <\alpha_{12} <\alpha_2 <\cdots <\alpha_{\theta-1} < \alpha_{1\theta} < \alpha_{2\theta} < \cdots <
\\
& \alpha_{1\theta}+\alpha_{1 \, \theta-1} <\alpha_{1\theta}+\alpha_{2 \, \theta-1} < \cdots < \alpha_{1\theta}+\alpha_{\theta-1} < \cdots <
\\
& \alpha_{2\theta}+\alpha_{2 \, \theta-1} <\cdots < \alpha_{2\theta}+\alpha_{\theta-1} < \cdots <
2\alpha_{\theta-1}+\alpha_{\theta} <\alpha_{\theta-1}+\alpha_{\theta} < \alpha_{\theta}.
\end{align*}
is a convex order associated to the following reduced expression of $w_0$:
\begin{align*}
s_1(s_2s_1)(s_3s_2s_1)\cdots (s_{\theta-1}\cdots s_1)(s_{\theta}s_{\theta-1}\cdots s_1)(s_{\theta}\cdots s_2)\cdots s_{\theta}.
\end{align*}
The Yetter-Drinfeld modules $M_{\alpha_{ij}}$, $\alpha_{ij}\ne \alpha_{\theta}$, are as in \eqref{eq:M-alphaij-defn}, and \eqref{eq:M-alphaij-basis} is a basis as well. Now, $M_{\alpha_{\theta}}$ is one-dimensional, spanned by $x_{\theta}$, and for the other roots we check that
\begin{align}\label{eq:M-gamma-defn}
M_{\alpha_{i\theta}+\alpha_{k\theta-1}} &= [M_{\alpha_{i\theta}}, M_{\alpha_{k\theta-1}}]_c.
\end{align}
A basis of $M_{\alpha_{i\theta}+\alpha_{k\theta-1}}$ is given by
\begin{align}\label{eq:M-alpha-gamma-n-basis}
x_{\beta} &= [x_{i \cdots\theta}, x_{k \cdots \theta-1}]_c, & x_{\un{\beta}}&=[x_{\un{i} i+1 \cdots \theta}, x_{k \cdots \theta-1}]_c.
\end{align}

\item[$\phi_4$] The element $w_0$ of maximal length has a reduced expression 
\begin{align}\label{eq:phi4-reduced-expression-w0}
s_1s_2s_3s_4s_2s_3s_1s_2s_3s_4s_1s_2s_3s_2s_1s_2s_3s_4s_2s_3s_2s_3s_4s_3,
\end{align}
which induces the following convex order on the set of positive roots:
{\small
\begin{align*}
&1, && 12, && 1^22^23, && 1^22^234, && 123, && 1^22^23^24, && 12^23, && 1^22^33^24,
\\
& 1^22^43^24, && 1^22^43^34, && 1234, && 12^234, && 1^22^43^34^2, && 12^23^24, && 12^33^24, && 2,
\\
& 2^23, && 2^234, && 23, && 2^23^24, && 234, && 3, && 34, && 4.
\end{align*}
}We denote by $\beta_i$ the $i$-th root according with this order.

Next we give, for each non-simple root $\beta=1^a2^b3^c4^d$ such that $d\ne 0$, a basis for each Yetter-Drinfeld submodule $M_{\beta}$ (if $d=0$, then we choose a basis as for $\gamma_{\theta}$):
\begin{align*}
M_{1^22^234}&: \big\{[x_{12},x_{\un{1}234}]_c\big\}, &
M_{1^22^23^24}&: \big\{ [x_{123},x_{\un{1}234}]_c \big\}, 
\\ 
M_{1^22^33^24}&: \big\{ [x_{12^23},x_{1234}]_c , [x_{\un{12^23}},x_{1234}]_c \big\}, &
M_{1^22^43^24}&: \big\{ [x_{12^23},x_{12^234}]_c \big\}, 
\\ 
M_{1^22^43^34}&: \big\{ [x_{12^23},x_{12^23^24}]_c \big\}, &
M_{1234}&: \big\{ x_{1234} , x_{\un{1}234}  \big\}, 
\\ 
M_{12^234}&: \big\{ [x_{1234},x_2]_c , [x_{\un{1}234},x_2]_c \big\}, &
M_{1^22^43^34^2}&: \big\{ [x_{12^234},x_{12^23^24}]_c \big\}, 
\\ 
M_{12^23^24}&: \big\{ [x_{1234},x_{23}]_c , [x_{\un{1}234},x_{23}]_c  \big\}, &
M_{2^234}&: \big\{ x_{\un{2}234} \big\}, 
\\ 
M_{12^33^24}&: \big\{ [x_{12^23^24},x_2]_c , [x_{\un{12^23^24}},x_2]_c \big\}, &
M_{2^23^24}&: \big\{ [x_{\un{2}3} , x_{234}]_c \big\}, 
\\ 
M_{234}&: \big\{ x_{234} , x_{\un{2}34} \big\}, &
M_{34}&: \big\{ x_{34} \big\}.
\end{align*}
where $x_{\beta}$ is the first vector fixed for $M_{\beta}$, while for $\dim M_{\beta}=2$, we denote by $x_{\un{\beta}}$ the second vectors in the order fixed above.
\end{description}

\begin{remark}\label{rem:degree-PBW-generators}
Let $\beta=1^{a_1}2^{a_2}\cdots \theta^{a_{\theta}}\in\Delta_+^M$, $g_{\beta}:=g_1^{a_1}g_2^{a_2}\cdots g_\theta^{a_{\theta}}\in G$. 
\begin{enumerate}[leftmargin=*]
\item If $\dim M_{\beta}=1$, then $x_{\beta}$ has $G$-degree $g_{\beta}$.
\item If $\dim M_{\beta}=2$, then $x_{\beta}$ has $G$-degree $g_{\beta}$ and $x_{\un{\beta}}$ has $G$-degree $g_{\beta}\kappa$.
\end{enumerate}
\end{remark}

We omit the details of the proof that either $\{x_{\beta}\}$ or $\{x_{\beta},x_{\un{\beta}}\}$ is a basis of $M_{\beta}$. The first step is to check that $M_{\beta}$ is spanned by $\{x_{\beta}\}$, respectively $\{x_{\beta},x_{\un{\beta}}\}$, using the defining relations of $\htoba(M)$; this can be done recursively on the convex (total) order. If $\dim M_{\beta}=2$ we see that $x_{\beta},x_{\un{\beta}}$ have different $G$-degree by Remark \ref{rem:degree-PBW-generators}, so they are linearly independent.

\begin{prop}\label{prop:Z-polynomial-generators}
The subalgebra $\Zc(M)=\htoba(M)^{\co \pi}$ of coinvariants under the canonical projection $\pi\colon\htoba(M)\twoheadrightarrow\toba(M)$ is a Hopf subalgebra of $\htoba(M)$. It is a skew-polynomial algebra in variables 
\begin{align}\label{eq:xbeta2}
&x_{\beta}^2, && \beta \in \Delta_+^{M}-\{\alpha_i\}, 
\\\label{eq:xbeta-xbetaprime}
&[x_{\beta},x_{\un{\beta}}]_c, && \beta \in \Delta_+^{M}-\{\alpha_i\} \text{ such that }\dim M_{\beta}=2.
\end{align}
\end{prop}
\pf
We proceed in several steps. In what follows, $\beta\in\Delta_+$ is not simple.

\begin{stepv}
For each $\beta\in\Delta_+$ such that $\dim M_{\beta}=1$, we have $x_{\beta}^2 \in M_{\beta}^2 \cap \Zc(M)$.

For each $\beta\in\Delta_+$ such that $\dim M_{\beta}=2$, we have $x_{\beta}^2$, $x_{\un{\beta}}^2$, $[x_{\beta},x_{\un{\beta}}]_c\in \Zc(M)$.
\end{stepv}

The subalgebra spanned by $M_{\beta}$ is isomorphic to the Nichols algebra $\toba(M_{\beta})$ by \cite{HS-right-coideal}. 
Assume first that $\dim M_{\beta}=1$. The braiding of $M_{\beta}$ satisfies that $c(x_{\beta}\ot x_{\beta})=-x_{\beta}\ot x_{\beta}$, so $x_{\beta}^2=0$ in $\toba(M)$. Thus $x_{\beta}^2\in \ker \pi \cap M_{\beta}^2$, and applying $F_\sigma$ we get
\begin{align*}
F_{\sigma}(x_{\beta}^2)\in F_{\sigma} \left(\ker \pi \cap M_{\beta}^2\right) =  \ker \pi \cap F_{\sigma} \left(M_{\beta}\right)^2 = \Zc(W).
\end{align*}
Hence $x_{\beta}^2\in \Zc(M)$, since $F_{\sigma}$ leaves the coalgebra structure unchanged.

Now, if $\dim M_{\beta}=2$, then the braiding of $M_{\beta}$ satisfies
\begin{align*}
c(x_{\beta}\ot x_{\beta})&=-x_{\beta}\ot x_{\beta}, &
c(x_{\un{\beta}}\ot x_{\un{\beta}})&=-x_{\un{\beta}}\ot x_{\un{\beta}}, &
c^2(x_{\beta}\ot x_{\un{\beta}})&=x_{\beta}\ot x_{\un{\beta}}.
\end{align*}
By a similar argument, $x_{\beta}^2$, $x_{\un{\beta}}^2$, $[x_{\beta},x_{\un{\beta}}]_c\in \Zc(M)$. \qed

\begin{stepv}
For each $\beta\in\Delta_+$ such that $\dim M_{\beta}=1$, $\{\xt_{\beta}\}$ is a basis of $F_{\sigma}(M_{\beta})$.

For each $\beta\in\Delta_+$ such that $\dim M_{\beta}=2$, $\{\xt_{\beta}, \xt_{\un{\beta}}\}$ is a basis of $F_{\sigma}(M_{\beta})$.
\end{stepv}

The statement certainly holds for simple roots, so we fix a non-simple root $\beta$. 
For types $\alpha_{\theta}$, $\delta_{\theta}$, $\epsilon_{\theta}$, we always have $\dim M_{\beta}=2$ and $M_{\beta}=[M_{\beta_1},M_{\beta_2}]_c$ for some $\beta_1, \beta_2\in\Delta_+$. Notice that $(\ad_c x_i)x_{\un{j}}=0$ for all $i,j\in\I$ since $i$, $\un{j}$ belong to different connected components of the Dynkin diagram of type $X_{\theta}\times X_{\theta}$. Hence $[\xt_{\gamma}, \xt_{\un{\delta}}]_c$ for all $\gamma,\delta\in \Delta_+$, and arguing recursively,
\begin{align*}
F_{\sigma} & (M_{\beta}) 
=F_{\sigma}\left([M_{\beta_1},M_{\beta_2}]_c\right)
=\big[ F_{\sigma}(M_{\beta_1}),F_{\sigma}(M_{\beta_2}) \big]_c
\\
=&\Big[ \Bbbk\xt_{\beta_1}+\Bbbk\xt_{\un{\beta_1}},\Bbbk\xt_{\beta_2}+\Bbbk\xt_{\un{\beta_2}} \Big]_c
= \Bbbk\Big[\xt_{\beta_1},\xt_{\beta_2}\Big]_c+\Bbbk\Big[\xt_{\un{\beta_1}},\xt_{\un{\beta_2}} \Big]_c
=\Bbbk\xt_{\beta}+\Bbbk\xt_{\un{\beta}}.
\end{align*}
For multiply-laced types, the proof follows similarly, case-by-case. \qed

\begin{stepv}
There exist $\N_0^{\theta}$-homogeneous elements
\begin{itemize}
\item $y_{\beta}\in M_{\beta}^2$ of $G$-degree $g_{\beta}^2$ when $\dim M_{\beta}=1$,
\item $y_{\beta}$, $y_{\un{\beta}}\in M_{\beta}^2$ of $G$-degree $g_{\beta}^2$, respectively $g_{\beta}^2\kappa$, when $\dim M_{\beta}=2$,
\end{itemize}
 that $q$-commute with every $G$-homogeneous element of $\htoba$. Moreover, $\Zc(M)$ is a skew-polynomial algebra in these variables.
\end{stepv}

Assume that $\dim M_{\beta}=2$. Note that, for all $i\in\I$, the elements $\xt_i^2$, $\xt_{\un{i}}^2$ are linear combinations of elements of $G$-degree $g_i^2$ and $g_i^2\kappa$ with 
non-trivial components on each degree. Hence, the elements $\xt_{\beta}^2$, $\xt_{\un{\beta}}^2$ are written as linear combinations of elements of $G$-degree $g_{\beta}^2$ and $g_{\beta}^2\kappa$ with non-trivial components on each degree, since they are obtained applying Lusztig's isomorphisms to appropriate $\xt_i^2$, $\xt_{\un{i}}^2$. Thus there exist $\yt_{\beta}$ and $\yt_{\un{\beta}}$ of $G$-degrees $g_{\beta}^2$ and $g_{\beta}^2\kappa$ respectively, that span the same as $\xt_{\beta}^2$ and $\xt_{\un{\beta}}^2$.
When $\dim M_{\beta}=1$ we may choose $\yt_{\beta}=\xt_{\beta}^2$. 

As $\Zc(W)$ is a skew-polynomial algebra in variables $x_{\beta}^2$, $x_{\un{\beta}}^2$ and each element $x_{\beta}^2$, $x_{\un{\beta}}^2$ is skew-central, the same holds with respect to $\yt_{\beta}$, $\yt_{\un{\beta}}$.

Let $\beta\in\Delta_+$ be such that $\dim M_{\beta}=2$. We set
\begin{align*}
y_{\beta} &:=F_{\sigma}^{-1}(\yt_{\beta})\in M_{\beta}^2, & 
y_{\un{\beta}}&:=F_{\sigma}^{-1}(\yt_{\un{\beta}})\in M_{\beta}^2.
\end{align*}
Then $y_{\beta}$, $y_{\un{\beta}}\in \Zc(M)$ since $F_{\sigma}$ preserves the coalgebra structure.
Note that 
\begin{align*}
F_{\sigma} (x_iy_{\beta})&=\sigma(g_i,g_{\beta}^2)F_{\sigma}(x_i)\yt_{\beta}, & 
F_{\sigma} (y_{\beta}x_i)&=\sigma(g_{\beta}^2,g_i)\yt_{\beta}F_{\sigma}(x_i),
\end{align*}
and these two elements differ up to a non-zero scalar for all $i\in\I$, thus $y_{\beta}$ is skew-central. The same happens for $y_{\un{\beta}}$, and for $y_{\beta}$ when $\dim M_{\beta}=1$.
In particular, the image under $F_{\sigma}$ of a multiplication of various $y_{\beta}$'s, $y_{\un{\beta}}$'s is the multiplication to the corresponding $\yt_{\beta}$'s, $\yt_{\un{\beta}}$'s up to a non-zero scalar, and the Step follows. \qed

\begin{stepv}
For each $\beta\in\Delta_+$ such that $\dim M_{\beta}=1$, we have $\Bbbk y_{\beta}=\Bbbk  x_{\beta}^2$.

For each $\beta\in\Delta_+$ such that $\dim M_{\beta}=2$, we have
$\Bbbk y_{\beta}= \Bbbk x_{\beta}^2$, $\Bbbk y_{\un{\beta}}=\Bbbk[ x_{\beta},  x_{\un{\beta}}]_c$.
\end{stepv}

Assume first that $\dim M_{\beta}=1$. Then $\dim M_{\beta}^2=1$ in $\htoba(M)$, and the claim follows since both $y_{\beta}$ and $x_{\beta}^2$ are generators of $M_{\beta}^2$.

Now assume that $\dim M_{\beta}=2$. In this case $\dim M_{\beta}^2=3$ in $\htoba(M)$ since $\dim F_{\sigma}(M_{\beta}^2)=3$ in $\htoba(W)$. On the other hand, $\dim M_{\beta}^2=1$ in $\toba(V)$: it is generated by $x_{\un{\beta}}x_{\beta}$ since $x_{\beta}^2=x_{\un{\beta}}^2=[ x_{\beta},  x_{\un{\beta}}]_c=0$. Notice that
$\dim F_{\sigma}(M_{\beta})^2 \cap \Zc(W)=2$, and $F_{\sigma}(M_{\beta})^2 \cap \Zc(W)$ contains elements with non-trivial components in degrees $g_{\beta}^2$ and $g_{\beta}^2\kappa$. Hence $\dim M_{\beta}^2 \cap \Zc(M)=2$, with one-dimensional homogeneous components of degrees $g_{\beta}^2$ and $g_{\beta}^2\kappa$.
Thus $\Bbbk y_{\beta}= \Bbbk x_{\beta}^2= \Bbbk x_{\un{\beta}}^2$ and $\Bbbk y_{\un{\beta}}=\Bbbk[ x_{\beta},  x_{\un{\beta}}]_c$, as claimed. \qed

\smallbreak
Hence Step 3 shows that $\Zc(M)$ is a skew-polynomial algebra, and Step 4 assures that we can choose generating variables as stated. 
\epf

\subsection{A presentation of the Nichols algebra}\label{subsec:presentation-general}
Here we put together the results obtained in \S \ref{subsec:distinguished}, \ref{subsec:coinvariants} to get a presentation, a PBW basis and the Hilbert series for the Nichols algebra.

\begin{theorem}\label{thm:Nichols-presentation-PBW-basis}
\begin{enumerate}[leftmargin=*,label=\rm{(\roman*)}]
\item\label{item:Nichols-PBW-basis} A set of PBW generators for $\toba(M)$ is given by 
\begin{align}\label{eq:PBW-generators}
& x_{\beta}, \, \beta\in\Delta_+; & &x_{\un{\beta}}\text{ when }\dim M_{\beta}=2.
\end{align}
The height of $x_{\beta}$, $x_{\un{\beta}}$ is $2$ for all $\beta\in \Delta_+$.

\item\label{item:Nichols-presentation} The Nichols algebra $\toba(M)$ is presented by generators $x_i$, $i\in \I$, $x_{\un{j}}$, $j\in\I_{\ell}$, and relations
\eqref{eq:dist-pre-NicholsMi-relations}, \eqref{eq:dist-pre-NicholsMi-dim1}, \eqref{eq:dist-pre-Nichols-aij=-1-relations},
\eqref{eq:dist-pre-Nichols-aij=0-relations}, \eqref{eq:dist-pre-Nichols-aij=0-relations-bis}, \eqref{eq:dist-pre-Nichols-rels-gamma},
\eqref{eq:dist-pre-Nichols-structureMjijk}, \eqref{eq:xbeta2} and \eqref{eq:xbeta-xbetaprime}.
\end{enumerate}
\end{theorem}

\pf
\ref{item:Nichols-PBW-basis} 
By \cite[Theorem 2.6]{HS-rank2-1}, the multiplication map
\begin{align*}
& \bigotimes_{\beta\in\Delta_+^M} \toba(M_{\beta}) \to \toba(M)
\end{align*}
is an isomorphism of $\Z^\theta$-graded objects in $\ydG$. If $\dim M_{\beta}=1$, then $M_{\beta}$ has braiding $-\id$ and $1,x_{\beta}$ is a basis of $\toba(M_{\beta})$. If $\dim M_{\beta}=2$, then $M_{\beta}$ has braiding as in \eqref{eq:ADE-self-braiding}: that is, $\toba(M_{\beta})$ is a quantum plane with basis $1,x_{\beta},x_{\un{\beta}},x_{\beta}x_{\un{\beta}}$, and the claim follows.

\ref{item:Nichols-presentation}
By Proposition \ref{prop:pre-Nichols}, relations \eqref{eq:dist-pre-NicholsMi-relations}, \eqref{eq:dist-pre-NicholsMi-dim1}, \eqref{eq:dist-pre-Nichols-aij=-1-relations}, \eqref{eq:dist-pre-Nichols-aij=0-relations}, \eqref{eq:dist-pre-Nichols-aij=0-relations-bis}, \eqref{eq:dist-pre-Nichols-rels-gamma} and
\eqref{eq:dist-pre-Nichols-structureMjijk} hold in $\toba(M)$. Also, \eqref{eq:xbeta2} and \eqref{eq:xbeta-xbetaprime} hold in $\toba(M)$ since the subalgebra generated by $M_{\beta}$ is isomorphic to $\toba(M_{\beta})$ as an algebra. Therefore, if $\toba$ denotes the quotient of $T(M)$ by all these relations, then there exists a canonical projection $\toba\twoheadrightarrow \toba(M)$ of graded Hopf algebras. Moreover,
\begin{align*}
\toba = \htoba(V)/ \langle \Zc(V)^+ \rangle,
\end{align*}
so $\mH_{\htoba(V)}=\mH_{\Zc(V)}\mH_{\toba}$ by \cite[Lemma 2.4]{ACS}. By Propositions \ref{prop:Z-polynomial-generators} and \ref{prop:pre-Nichols},
\begin{align*}
\mH_{\toba} = \left(\prod_{\beta \in \Delta_+^{V}:\dim M_{\beta}=2 } (1+t^{\beta})^2\right) 
\left(\prod_{\beta \in \Delta_+^{V}:\dim M_{\beta}=1 } (1+t^{\beta})\right)=\mH_{\toba(V)},
\end{align*}
and we deduce that $\toba=\toba(M)$.
\epf

\begin{remark}
The order of the elements in the PBW basis in Theorem \ref{thm:Nichols-presentation-PBW-basis} \ref{item:Nichols-PBW-basis} is given by the expression of the element $w_0$ of maximal length fixed below. For example, for type $\phi_{4}$ we have the following PBW basis:
\begin{align}\label{eq:phi4-PBW-basis}
\begin{aligned}
& x_{4}^{a_{1}} x_{34}^{a_{2}} x_{\un{234}}^{a_{3}} x_{234}^{a_{4}} x_{2^23^24}^{a_{5}} x_{\un{23}}^{a_{6}} 
x_{23}^{a_{7}} x_{2^234}^{a_{8}} x_{2^23}^{a_{9}} x_{\un{2}}^{a_{10}} x_{2}^{a_{11}} x_{\un{12^33^24}}^{a_{12}} 
x_{12^33^24}^{a_{13}} x_{\un{12^23^24}}^{a_{14}} \\
& x_{12^23^24}^{a_{15}} x_{1^22^43^34^2}^{a_{16}} x_{\un{12^234}}^{a_{17}} x_{12^234}^{a_{18}} 
x_{\un{1234}}^{a_{19}} x_{1234}^{a_{20}} x_{34}^{a_{21}} x_{1^22^43^34}^{a_{22}} x_{1^22^43^24}^{a_{23}} x_{\un{1^22^33^24}}^{a_{24}} 
\\
& x_{1^22^33^24}^{a_{25}} x_{\un{12^23}}^{a_{26}} x_{12^23}^{a_{27}} x_{1^22^23^24}^{a_{28}}  x_{\un{123}}^{a_{29}} x_{123}^{a_{30}} 
x_{1^22^234}^{a_{31}} x_{1^22^23}^{a_{32}} x_{\un{12}}^{a_{33}} x_{12}^{a_{34}} x_{\un{1}}^{a_{35}} x_{1}^{a_{36}} , \, a_i\in\I_{0,1}.
\end{aligned}
\end{align}

\end{remark}

\subsection{Rigidity of Nichols algebras}
We briefly discuss rigidity for Nichols algebras of types $\alpha_{\theta}$, $\gamma_{\theta}$, $\delta_{\theta}$, $\epsilon_{\theta}$ and $\phi_4$, inspired by \cite{AKM} where rigidity for finite-dimensional Nichols algebras over abelian groups is studied. This will come in handy in the next section, where  we study the liftings of these Nichols algebras.

\medbreak

Let $\R_M\subset T(M)$ be the set of defining relations as in Theorem \ref{thm:Nichols-presentation-PBW-basis}.
We start by describing the Yetter-Drinfeld structure of $\Bbbk\R_M$. Recall the index $\ell$ introduced in \eqref{eq:ell}.

\begin{remark}\label{rem:YD-structure-relations}
\begin{enumerate}[leftmargin=*]
\item The $G$-degree of $x_{i\un{i}}$ is $g_i^2\kappa$. Using Example \ref{ex:simple-YD-dim2}, it is easy to see that
\begin{align}\label{eq:character-xibari}
&\un{\chi}_{i}:G\to \Bbbk^{\times}, &
&\un{\chi}_{i}(h):=\begin{cases} 
\chi_i(hg_j^{-1}hg_j)  & h\in G^{g_i},
\\
\chi_i(\kappa)\chi_i(h^2) & h\notin G^{g_i},
\end{cases}
\end{align}
is a character, and the $G$-action on $x_{i\un{i}}$ is given by $\un{\chi}_{i}$.

\item The $G$-degree of $x_i^2$ and $x_{\un{i}}^2$ is $g_i^2$. The $G$-action when $i\le \ell$ is
\begin{align*}
g\cdot x_i^2 &=\begin{cases}
\chi_i^2 (g)x_i^2, & g\in G^{g_i}, \\
\chi_i^2 (g_j^{-1}g)x_{\un{i}}^2, & g\notin G^{g_i}; 
\end{cases}
&
g\cdot x_{\un{i}}^2 &=\begin{cases}
\chi_i^2 (g)x_{\un{i}}^2, & g\in G^{g_i}, \\
\chi_i^2 (gg_j)x_i^2, & g\notin G^{g_i}. 
\end{cases}
\end{align*}
If $i>\ell$, then the action on $x_i^2$ is given by $\chi_i^2$.

\item Let $i<j\le \ell$ be such that $a_{ij}=-1$, and set
\begin{align*}
r_1 &:=x_{i\un{j}}+\chi_j(g_i^2) x_{\un{i}j}, & r_2&:=x_{\un{i}\un{j}}+\chi_i(\kappa) x_{ij}.
\end{align*}
The $G$-degrees of $r_1$ and $r_2$ are $g_ig_j\kappa$ and $g_ig_j$. 
The $G$-action is given by

\begin{center}
\begin{tabular}{|c|c|c|}
\hline
$h\in$ & $h\cdot r_1$ &  $h\cdot r_2$
\\ 
\hline
$G^{g_i}\cap G^{g_j}$ & $\chi_i(h)\chi_j(h) r_1$ & $\chi_i(h)\chi_j(h) r_2$
\\
\hline
$G^{g_i}-G^{g_j}$ & $\chi_i(\kappa h)\chi_j(g_ih) r_2$ & $\chi_i(\kappa h)\chi_j(g_i^{-1}h) r_2$
\\
\hline
$G^{g_j}-G^{g_i}$ & $\chi_i(g_j^{-1}h)\chi_j(\kappa h) r_2$ & $\chi_i(hg_j)\chi_j(\kappa h) r_1$
\\
\hline
$G-(G^{g_i}\cup G^{g_j})$ & $\chi_i(hg_j)\chi_j(g_ih) r_1$ & $\chi_i(hg_j^{-1})\chi_j(g_i^{-1}h) r_2$
\\
\hline
\end{tabular}
\end{center}

\item Let $i<j\le \ell$ be such that $a_{ij}=0$. Both $x_{ij}$ and $x_{\un{i}\un{j}}$ have $G$-degree $g_ig_j$, while the $G$-degree of $x_{i\un{j}}$ and $x_{\un{i}j}$ is $g_ig_j\kappa$. For the action, as $G^{g_i}$ and $G^{g_j}$ are both subgroups of index 2, there are two possibilities. 
If $G^{g_i} \neq G^{g_j}$\footnote{This occurs, for example, in type $\alpha_\theta$ whenever $j-i\ge3$. We may choose $g_a=g_{j-1}\ne g_b=g_{i+1}$.}, then we choose $g_a\in G^{g_i}-G^{g_j}$ and $g_b\in G^{g_j}-G^{g_i}$. By Example \ref{ex:simple-YD-dim2}, the $G$-action is given by:
\smallbreak
\hspace{-35pt}
\begin{adjustbox}{max width=\textwidth}
\begin{tabular}{|c|c|c|c|c|}
\hline
$h\in$ & $h\cdot x_{ij}$ &  $h \cdot x_{\un{i}\un{j}}$ & $h\cdot x_{i\un{j}}$ &  $h\cdot x_{\un{i}j}$
\\ 
\hline
$G^{g_i}\cap G^{g_j}$ & $\chi_i(h)\chi_j(h) x_{ij}$ & $\chi_i(g_b^{-1}hg_b)\chi_j(g_a^{-1}hg_a) x_{\un{i}\un{j}}$ 
& $\chi_i(h)\chi_j(g_a^{-1}hg_a)x_{i\un{j}}$ & $\chi_i(g_b^{-1}hg_b)\chi_j(h)x_{\un{i}j}$
\\
\hline
$G^{g_i}-G^{g_j}$ & $\chi_i(h)\chi_j(g_a^{-1}h) x_{i\un{j}}$ & $\chi_i(g_b^{-1}hg_b)\chi_j(hg_a) x_{\un{i}j}$ 
& $\chi_i(h)\chi_j(hg_a)x_{ij}$ & $\chi_i(g_b^{-1}hg_b)\chi_j(g_a^{-1}h)x_{\un{i}\un{j}}$
\\
\hline
$G^{g_j}-G^{g_i}$ & $\chi_i(g_b^{-1}h)\chi_j(h) x_{\un{i}j}$ & $\chi_i(hg_b)\chi_j(g_a^{-1}hg_a) x_{i\un{j}}$ 
& $\chi_i(g_b^{-1}h)\chi_j(g_a^{-1}hg_a)x_{\un{i}\un{j}}$ & $\chi_i(hg_b)\chi_j(h)x_{ij}$
\\
\hline
$G-(G^{g_i}\cup G^{g_j})$ & $\chi_i(g_b^{-1}h)\chi_j(g_a^{-1}h) x_{\un{i}\un{j}}$ & $\chi_i(hg_b)\chi_j(hg_a) x_{ij}$ 
& $\chi_i(g_b^{-1}h)\chi_j(hg_a)x_{\un{i}j}$ &  $\chi_i(hg_b)\chi_j(g_a^{-1}h)x_{i\un{j}}$
\\
\hline
\end{tabular}
\end{adjustbox}
\smallbreak
Hence $x_{ij}$, $x_{\un{i}\un{j}}$, $x_{i\un{j}}$, $x_{\un{i}j}$ span a 4-dimensional irreducible Yetter-Drinfeld submodule.

In the case $G^{g_i}=G^{g_j}$ we choose $g_a=g_b \notin G^{g_i}$, and the action is described by the second and last rows of the previous table.
There are two Yetter-Drinfeld submodules, spanned by $\{x_{ij}, x_{\un{i}\un{j}}\}$, and by $\{x_{i\un{j}},x_{\un{i}j}\}$.

\item Let $i\le \ell<j$ be such that $a_{ij}=0$. The $G$-degree of $x_{ij}$ is $g_ig_j$, while the $G$-degree of $x_{\un{i}j}$ is $g_ig_j\kappa$. 
For the action, notice that $g_j\in Z(G)$: if we pick $\mathbf{g}\notin G^{g_i}$, then
\begin{center}
\begin{tabular}{|c|c|c|c|c|}
\hline
$h\in$ & $h\cdot x_{ij}$ &  $h\cdot x_{\un{i}j}$
\\ 
\hline
$G^{g_i}$ & $\chi_i(h)\chi_j(h) x_{ij}$ & $\chi_i(\mathbf{g}^{-1}h\mathbf{g})\chi_j(h)x_{\un{i}j}$
\\
\hline
$G-G^{g_i}$ & $\chi_i(\mathbf{g}^{-1}h)\chi_j(h) x_{\un{i}j}$ &  $\chi_i(h\mathbf{g})\chi_j(h)x_{ij}$
\\
\hline
\end{tabular}
\end{center}

\item For $i,j$ such that $a_{ij}=-2$, the $G$-degree of $(\ad_c x_j) x_{ij\un{i}}$ is $g_i^2g_j^2\kappa$, where $G$ acts via $\un{\chi}_1\chi_j^2$.

\item Let $i<j<k$ be such that $a_{ji}=a_{jk}=-1$. Consider
\begin{align*}
\mathbf{r}_1&:=(\ad_c x_j) x_{ijk}, & \mathbf{r}_2&:=(\ad_c x_j) x_{\un{i}jk},
\end{align*} 
which have $G$-degrees $g_ig_j^2g_k\kappa$ and $g_ig_j^2g_k$ respectively. The $G$-action is given by
\begin{center}
\begin{tabular}{|c|c|c|}
\hline
$h\in$ & $\un{\chi}_j(h)^{-1}h\cdot \mathbf{r}_1$ &  $\un{\chi}_j(h)^{-1}h\cdot \mathbf{r}_2$
\\ 
\hline
$G^{g_i}\cap G^{g_k}$ & $\chi_i(h)\chi_k(h) \mathbf{r}_1$ & $\chi_i(g_j^{-1}hg_j)\chi_k(h)\mathbf{r}_2$
\\
\hline
$G^{g_i}-G^{g_k}$ & $\chi_i(h)\chi_k(g_j^{-1}h)\mathbf{r}_2$ & $\chi_i(g_j^{-1}hg_j)\chi_k(g_j^{-1}h)\mathbf{r}_1$
\\
\hline
$G^{g_k}-G^{g_i}$ & $\chi_i(g_j^{-1}h)\chi_k(h) \mathbf{r}_2$ & $\chi_i(hg_j)\chi_k(h) \mathbf{r}_1$
\\
\hline
$G-(G^{g_i}\cup G^{g_k})$ & $\chi_i(g_j^{-1}h)\chi_k(g_j^{-1}h) \mathbf{r}_1$ &  $\chi_i(hg_j)\chi_k(g_j^{-1}h)\mathbf{r}_2$
\\
\hline
\end{tabular}
\end{center}

\item Let $\beta \in \Delta_+^{M}-\{\alpha_i\}$. The $G$-degrees of $x_{\beta}^2$ and $[x_{\beta},x_{\un{\beta}}]_c$ are $g_{\beta}^2$ and $g_{\beta}^2\kappa$, respectively, where $g_{\beta}$ is as in Remark \ref{rem:degree-PBW-generators}. We define accordingly 
$\un{\chi}_{\beta}:=\un{\chi}_1^{a_1} \cdots \un{\chi}_{\theta}^{a_\theta}$, and $G$ acts on $x_{\beta}^2$ and $[x_{\beta},x_{\un{\beta}}]_c$ via $\un{\chi}_{\beta}$. 
\end{enumerate}
\end{remark}

\begin{remark}\label{rem:centralizer-gigj-aij=-1}
Let $i,j\in\I$ such that $G^{g_i}\ne G^{g_j}$ and $g_i,g_j\notin Z(G)$. Then
$$ G^{g_ig_j}=(G^{g_i}\cap G^{g_j}) \cup \big(G-(G^{g_i}\cup G^{g_j})\big) $$ 
and the following rule defines a character
\begin{align}\label{eq:centralizer-g1g2-character}
\chi_{ij}: & G^{g_ig_j}\to \Bbbk^{\times}, &
\chi_{ij}(h) &:= \begin{cases}
\chi_i(h)\chi_j(h), &  h\in G^{g_i}\cap G^{g_j},\\
\chi_i(hg_j)\chi_j(hg_i), &  h\notin G^{g_i}\cup G^{g_j}.
\end{cases} 
\end{align}
\end{remark}

\begin{theorem}\label{thm:trivial-hom-rels-V}
Let $M\in\ydG$ of type $\alpha_{\theta}$, $\gamma_{\theta}$, $\delta_{\theta}$, $\epsilon_{\theta}$ or $\phi_4$. Then 
$$\Hom^{\Bbbk G}_{\Bbbk G}(\Bbbk\R_M,M)=0  .$$
\end{theorem}
\pf
Let $\mathtt{r}\in\R_M$ be $G$-homogeneous of degree $\mathtt{g}\in G$. By direct computation, $\mathtt{g}\cdot\mathtt{r}=\mathtt{r}$. On the other hand, $M$ has a basis
$\{x_i|i\in\I\}\cup \{x_{\un{j}}|j\in\Jb\}$, where $\Jb=\{i\in\I|\dim M_i=2\}$. Here $x_i$ has degree $g_i$, while $x_{\un{j}}$ has degree $g_j\kappa$, and
\begin{align*}
g_i\cdot x_i&=-x_i, \quad i\in\I; & g_j\kappa\cdot x_{\un{j}}&=-x_{\un{j}},\quad j\in\Jb.
\end{align*}
Hence the claim follows. 
\epf

Recall that a graded braided bialgebra is \emph{rigid} if it has no non-trivial graded deformations. 
See \cite[\S 2]{AKM} and the references therein for details. Next we address rigidity for $\toba(M)$.

\begin{theorem}\label{thm:rigidity}
Let $M\in\ydG$ be of type either $\alpha_{\theta}$, $\gamma_{\theta}$, $\delta_{\theta}$, $\epsilon_{\theta}$ or $\phi_4$.
Then $\toba(M)$ is rigid.
\end{theorem}
\pf
The category $\ydG$ is semisimple and $\Hom^{\Bbbk G}_{\Bbbk G}(\Bbbk\R_M,M)=0$ by Theorem \ref{thm:trivial-hom-rels-V}.
Hence \cite[Theorem 5.3]{AKM} applies for $\toba(M)$.
\epf

\begin{remark}
The previous notion of rigidity is related to another one introduced in \cite{Meir} coming from the action of an appropriate algebraic group
on the Nichols algebra (viewed as a braided Hopf algebra). In fact, the notion of rigidity in loc. cit. is equivalent to generation in degree
one, which holds by Theorem \ref{thm:gen-degree-one}. This gives a different proof of Theorem \ref{thm:rigidity}, independent of Theorem \ref{thm:trivial-hom-rels-V}. Anyway, we need Theorem \ref{thm:trivial-hom-rels-V} to compute liftings.
\end{remark}

\section{Liftings of Nichols algebras}\label{sec:liftings}

We describe all \emph{liftings} for Nichols algebras of types $\alpha_{\theta}$, $\gamma_{\theta}$, $\delta_{\theta}$, $\epsilon_{\theta}$ and $\phi_4$. Even when the braided vector space is of diagonal type (that is, when $\kappa$ acts trivially) we cannot invoke \cite{AnG} since the Yetter-Drinfeld realizations considered here are not \emph{principal}. 
Nevertheless we will perform an adaptation of the strategy developed in \cite{AAGMV,AnG}. 

We study the lowest rank type $\alpha_2$ first, with a double purpose. On the one hand, we will only show all the details in this case, with explicit formulas for the defining relations. 
On the other hand, it will be the starting point to prove the general case, in which we will conclude that all liftings are cocycle deformations of the associated graded Hopf algebras.

\smallbreak

Recall that a lifting of $M$ over $G$ is a finite-dimensional Hopf algebra $H$ with coradical $\Bbbk G$ and infinitesimal braiding $M$. Hence $\gr H \simeq \toba(M) \# \Bbbk G$ by Theorem \ref{thm:gen-degree-one}.

The family of liftings of $M$ over $G$ will be indexed by a set $\cR_{M} \subseteq \Bbbk^K$ of \emph{deformation parameters}, where $K$ is the number of suitable chosen Yetter-Drinfeld submodules of the subspace spanned by a minimal set $\cG$ of generators for the ideal defining $\toba(M)$.

For each $\bsl\in\cR_{ M }$ and $i\in \I_K$, we define $\bsl^{(i)},  \bsl^{(-i)}\in\Bbbk^K$ by
\begin{align}\label{eq:lambda-i}
(\bsl^{(i)})_j&:=\begin{cases}
0 & j\ne i, \\ \lambda_i & j=i;
\end{cases}
&
(\bsl^{(-i)})_j&:=\begin{cases}
\lambda_i & j\ne i, \\ 0 & j=i;
\end{cases}
& j&\in\I_K.
\end{align}
The aforementioned strategy starts by choosing a \emph{good} stratification $\cG= \cG_0 \sqcup \cG_1 \sqcup \cdots  \sqcup \cG_l$, meaning that 
the vector space spanned by $\cG_k$ is a Yetter-Drinfeld submodule of $\toba(M)$ and the elements of $\cG_k$ are primitive in the braided Hopf algebra $\toba_k := T (M) / \langle \cup_{j=0}^{k-1} \cG_j \rangle$, $k  \in \I_{l+1}$, with one possible exception: we do not requite primitiveness for the last step.

\subsection{Liftings of type \texorpdfstring{$\alpha_2$}{}}\label{subsec:liftings-alpha2}
Let $M\in\ydG$ of type $\alpha_2$. Let $\cR_{M}$ be the set of tuples
$ \bsl=(\mu_1,\mu_2,\lambda_{1},\lambda_{2},\lambda_{12},\mu_{12},\mu_{12}')\in\Bbbk^7$
that satisfy the constraints 
\begin{align}\label{eq:lambda-conditions-alpha2}
\begin{aligned}
\mu_{i} &= 0 \text { if either } \chi_i^2 \neq \varepsilon \text { or } g_i^2 = 1, \, i\in\I_2,
\\
\lambda_{i} &= 0 \text { if either } \un{\chi}_i \neq \varepsilon \text { or } g_i^2 = \kappa, \, i\in\I_2,
\\
\lambda_{12} &= 0 \text { if } \chi_{12} \neq \varepsilon,
\\
\mu_{12} &= 0 \text { if either } \un{\chi}_1\un{\chi}_2 \neq \varepsilon \text{ or } (g_1g_2)^2=1.
\\
\mu_{12}' &= 0 \text { if either } \un{\chi}_1\un{\chi}_2 \neq \varepsilon \text{ or } (g_1g_2)^2=\kappa.
\end{aligned}
\end{align}
The definition of $\un{\chi}_i$ and $\chi_{ij}$ was given in \eqref{eq:character-xibari} and \eqref{eq:centralizer-g1g2-character}. This subsection is devoted to prove the following:
\begin{theorem}\label{thm:liftings-alpha2} 
Let $M\in\ydG$ be of type $\alpha_2$. 
For each $\bsl \in \cR_{M}$, let $\cL(\bsl)$ be the quotient of $T(M)\# \Bbbk G$ by the following set of relations:
\begin{align*}
&z_i^2 - \mu_i ( 1 - g_i^2), \qquad \qquad \qquad  z_{i\un{i}}-\lambda_i(1-g_i^2\kappa),
\\
& z_{\un{1}\un{2}}+\chi_1(\kappa)z_{12}- \lambda_{12} (1 - g_1g_2 ),
\\
&z_{12}^2+\lambda_1\mu_2(1-g_1^2\kappa)g_2^2-\chi_1(\kappa)\mu_1\lambda_2(1-g_1^2)g_2^2\kappa -\mu_{12}(1-g_1^2g_2^2\kappa), 
\\
&[z_{12},z_{\un{1}2}]_c+2(1+\chi_1\chi_2(\kappa))\mu_1\mu_2(1-g_1^2)g_2^2-\lambda_1\lambda_2(\kappa-g_1^2)g_2^2-\mu_{12}'(1-g_1^2g_2^2),
\end{align*}
where we changed the labels $(x_i, x_{\un{i}})_{i \in I_2}$ of the generators of $T(M)$ to $(z_i, z_{\un{i}})_{i \in I_2}$. Then:
\begin{enumerate}[leftmargin=*,label=\rm{(\alph*)}]
\item\label{item:liftings-alpha2-biGalois} $\cL(\bsl)\simeq L(\cA(\bsl),\toba(M)\#\Bbbk G)$.
\item\label{item:liftings-alpha2-lifting} $\cL(\bsl)$ is a lifting of $M$ over $\Bbbk G$.
\item\label{item:liftings-alpha2-cocycle} $\cL(\bsl)$ is a cocycle deformation of $\toba(M)\#\Bbbk G$.
\end{enumerate}

Conversely, if $L$ is lifting of $M$ over $\Bbbk G$, then there exist $\bsl \in \cR_{M}$ 
such that $L \simeq \cL(\bsl)$.
\end{theorem}

Fix a Yetter-Drinfeld module $M$ over $G$ of type $\alpha_2$.
As
$x_i^2$, $x_{\un{i}}^2$, $x_{i\un{i}}$, $x_{1\un{2}}+\chi_2(g_1^2) x_{\un{1}2}$ and $x_{\un{1}\un{2}}+\chi_1(\kappa) x_{12}$ are primitive in $T(M)$ and
\begin{align*}
\Delta & (x_{12}^2) = x_{12}^2 \ot 1+ x_{1\un{1}} g_2^2 \ot x_2^2 
- \chi_{1}(\kappa) x_1^2 g_2^2 \kappa \ot x_{2\un{2}}+ g_1^2g_2^2\kappa \ot x_{12}^2,
\\ 
\Delta & ([x_{12}, x_{\un{1}2} ]_c) = [x_{12}, x_{\un{1}2} ]_c \ot 1 - \chi_{1}^{-1}(g_2^2) x_{1\un{1}} g_2^2\kappa \ot x_{2\un{2}} +\chi_{1}(g_2^2)x_1^2 g_2^2 \ot x_2^2
\\ & +\chi_{1}(g_2^2)x_{\un{1}}^2 g_2^2 \ot x_{\un{2}}^2 + \chi_1\chi_2(\kappa) x_{\un{1}}^2 g_2^2 \ot x_2^2
+x_1^2 g_2^2 \ot x_{\un{2}}^2 + g_1^2g_2^2 \ot [x_{12}, x_{\un{1}2} ]_c,
\end{align*}
we may choose the following stratification:
\begin{align*}
\cG_0=& \{x_i^2, \, x_{\un{i}}^2, \, x_{i\un{i}} \}  ,& \cG_1&= \{ x_{1\un{2}}+\chi_2(g_1^2) x_{\un{1}2}, \,  x_{\un{1}\un{2}}+\chi_1(\kappa) x_{12} \}, &
\cG_2 =& \{ x_{12}^2 , \, [x_{12}, x_{\un{1}2} ]_c\} .
\end{align*}
The Yetter-Drinfeld structure for each stratum is given in Remark \ref{rem:YD-structure-relations}.

Let $\mH_k:=\toba_k\# \Bbbk G$. Next we introduce a family of cleft objects of $\mH_k$ parametrized by the set $\cR_{M}$. 
Given $\bsl\in \cR_{M}$, define $\cE_0(\bsl) =\toba_0=T(M)$, $\cE_1(\bsl) =\toba_1$, but we change the labels of the generators to $(y_i, y_{\un{i}})_{i\in\I_2}$ in order to differentiate with generators $(x_i,x_{\un{i}})_{i\in\I_2}$ of the pre-Nichols algebras $\toba_k$. Let
\begin{align*}
\cE_1(\bsl) & := \cE_0(\bsl) / \left\langle y_i^2 - \mu_i, y_{\un{i}}^2 - \mu_i, y_{i\un{i}} - \lambda_i \colon i\in\I_2\right\rangle,
\\
\cE_2(\bsl) & := \cE_1(\bsl) / \left\langle y_{1\un{2}}+\chi_2(g_1^2) y_{\un{1}2}- \lambda_{12}, \, y_{\un{1}\un{2}}+\chi_1(\kappa) y_{12}- \lambda_{12} \right\rangle, 
\\
\cE_3(\bsl) &:= \cE_2(\bsl)/ \left\langle y_{12}^2-\mu_{12}, [y_{12},y_{\un{1}2}]_c-\mu_{12}' \right\rangle.
\end{align*} 

Each $\cE_{i}(\bsl)$ is a $\Bbbk G$-module algebra since the ideal is stable under the $G$-action by \eqref{eq:lambda-conditions-alpha2}. Thus we may introduce $\cA_i(\bsl):=\cE_i(\bsl) \# \Bbbk G$.

\begin{lemma}\label{lem:cleft-objects-alpha2}
Let $k\in\I_{3}$. Then $\cE_k(\bsl)\neq 0$ and each $\cA_k(\bsl)$ is an $\mH_k$-cleft object. There exists an $\mH_k$-colinear section $\gamma_k \colon \mH_k \to \cA_k$ that restricts to an algebra map $(\gamma_k) _{|\Bbbk G} \in \Alg (\Bbbk G, \cA_k)$.
\end{lemma}
\pf
Fix $\bsl\in\cR_M$; to simplify the notation we suppress $\bsl$ and put $\cE_k=\cE_k(\bsl)$, $\cA_k=\cA_k(\bsl)$. We prove the claim recursively on $k$. 

\smallbreak

For $k=1$, we notice that $\cE_1 \neq 0$ (and a fortiori $\cA_1\neq 0$) by \cite[Lemma 5.16]{AAGMV}.
Notice that $g_j(y_i^2-\mu_i)g_j^{-1}=y_{\un{i}}^2-\lambda_i$ if $i\ne j$, so in $\cA_1$ we have
\begin{align*}
\langle y_i^2-\mu_i, y_{\un{i}}^2-\mu_i : i=1,2 \rangle
&=\langle y_i^2-\mu_i : i=1,2 \rangle.
\end{align*}
We may refine the stratification and proceed in four steps, quotient out  first by $x_1^2$, then by $x_2^2$, now by $x_{1\un{1}}$ and finally by $x_{2\un{2}}$. At each step we consider the subalgebra $Y'$ generated by the relation $r$ in the corresponding pre-Nichols algebra, note that $Y'$ is isomorphic to a polynomial ring in one variable since $r \in \mP(T(M))_{g}-0$, and for this $g$ we have $g\cdot r=r$. Consider $Y=\Ss(Y')$. 
As $Y$ is a polynomial algebra generated by $rg^{-1}$, there exists an algebra map $\phi:Y\to\cA$ such that
$\phi(rg^{-1})=rg^{-1}-\lambda g^{-1}$, $\lambda\in\Bbbk$, which is $\mH$-colinear. Applying repeatedly \cite[Theorem 8]{Gu} as in \cite[Proposition 5.19]{AAGMV},
$\cA_1$ is a $\mH_1$-cleft object and the existence of the desired section $\gamma_1$ follows by \cite[Proposition 5.8]{AAGMV}.

\smallbreak

For $k=2$, it is enough to show that $\cA_2\ne 0$. Indeed, in that case \cite[Theorem 8]{Gu} assures that $\cA_2$ is 
an $\mH_2$-cleft object. Now \cite[Proposition 5.8]{AAGMV} provides a section $\gamma_2$ such that $(\gamma_2) _{|\Bbbk G} \in \Alg (\Bbbk G, \cA_2)$.
As in \cite[Lemma 3.4]{AnG}, non-vanishing of $\cA_2$ would follow from 
$$ \cE_2(\bsl^{(5)}) =\toba_1 / \left\langle y_{1\un{2}}+\chi_2(g_1^2) y_{\un{1}2}- \lambda_{12}, \, y_{\un{1}\un{2}}+\chi_1(\kappa) y_{12}- \lambda_{12} \right\rangle \ne 0. $$ 
Indeed, if $\varpi_1\colon \cA_1(\bsl^{(5)})=\mH_1\twoheadrightarrow \cE_2(\bsl^{(5)})\# \Bbbk G$ is the canonical projection, then
the composition of the algebra map $(\rho_1)_{|\cE_1}\colon\cE_1\to \cA_1\ot\mH_1$ with $\id\ot\varpi_1$ factors through $\cE_2=\cE_2(\bsl)$.

To check that $\cE_2(\bsl^{(5)})\neq0$ we use that $(\mH_1)_{\sigma}$ is the bosonization of a pre-Nichols algebra of diagonal type by $G$, and that the 
$(1,g_1g_2\kappa)$- and $(1,g_1g_2)$-primitive elements $y_{1\un{2}}+\chi_2(g_1^2) y_{\un{1}2}$ and $y_{\un{1}\un{2}}+\chi_1(\kappa) y_{12}$ span the same subspace as $\xt_{1\un{2}}$ and $\xt_{\un{1}2}$, see the proof of Proposition \ref{prop:pre-Nichols}. The quotient 
$(\mH_1)_{\sigma}/\langle \xt_{1\un{2}}-\lambda_{12}, \xt_{\un{1}2}-\lambda_{12} \rangle$ is not zero by \cite{AAG}, and
\begin{align*}
(\mH_1)_{\sigma}/\langle \xt_{1\un{2}}-\lambda_{12}, \xt_{\un{1}2}-\lambda_{12} \rangle \simeq F_{\sigma}\left(\cA_2(\bsl^{(5)})\right),
\end{align*}
which implies $\cA_2(\bsl^{(5)})\ne 0$.
Notice that 
$$ g_1 \left(y_{1\un{2}}+\chi_2(g_1^2) y_{\un{1}2}- \lambda_{12}\right)g_1^{-1}= -\chi_1(\kappa)\chi_2(g_1^2) \left(y_{\un{1}\un{2}}+\chi_1(\kappa) y_{12}- \lambda_{12}\right), $$
so in $\cA_2$,
$$ \langle y_{1\un{2}}+\chi_2(g_1^2) y_{\un{1}2}- \lambda_{12}, \, y_{\un{1}\un{2}}+\chi_1(\kappa) y_{12}- \lambda_{12} \rangle
=\langle y_{\un{1}\un{2}}+\chi_1(\kappa) y_{12}- \lambda_{12} \rangle. $$

\smallbreak

Finally, as $\mH_2=\htoba(M)\# \Bbbk G$, $\mH_3=\toba(M)\# \Bbbk G$, we have $\mH_2^{\co \pi_2}=\Zc(M)$, a skew-polynomial algebra in variables 
$x_{12}^2$, $[y_{12},y_{\un{1}2}]_c$ by Proposition \ref{prop:Z-polynomial-generators}. Hence \cite[Theorem 4]{Gu} applies and $\cA_3$ is $\mH_3$-cleft. The claim about $\gamma_3$ follows from \cite[Proposition 5.8]{AAGMV}.
\epf

\pf[Proof of Theorem \ref{thm:liftings-alpha2}]
Now follows by the same procedure as in \cite[Theorem 5.6]{AnS}, using Theorem \ref{thm:trivial-hom-rels-V}. Indeed, if we define $\cL_0(\bsl)=\mH_0$, 
\begin{align*}
\cL_1(\bsl) & = \cL_0(\bsl) / \left\langle z_i^2 - \mu_i ( 1 - g_i^2), z_{i\un{i}}-\lambda_i(1-g_i^2\kappa) \right\rangle,
\\
\cL_2(\bsl) & = \cL_1(\bsl) / \left\langle  z_{\un{1}\un{2}}+\chi_1(\kappa)z_{12}- \lambda_{12} (1 - g_1g_2 ) \right\rangle,
\end{align*}
and $\cL_3(\bsl)=\cL(\bsl)$, we can prove recursively that $\cL_i(\bsl)\simeq L(\cA_i(\bsl),\mH_i)$.
\epf

\subsection{The general case}\label{subsec:liftings-general}
Let $M \in\ydG$ of type $\alpha_{\theta}$, $\gamma_{\theta}$, $\delta_{\theta}$, $\epsilon_{\theta}$ or $\phi_4$.
The definition of $\un{\chi}_i$ and $\chi_{ij}$ was given in \eqref{eq:character-xibari} and \eqref{eq:centralizer-g1g2-character}. 
The set $\cR_{M}$ of deformation parameters contain tuples $\bsl$ satisfying the following constraints:
\begin{align}\label{eq:lambda-conditions-general}
\begin{aligned}
\mu_{i} &= 0 \text { if either } \chi_i^2 \neq \varepsilon \text { or } g_i^2 = 1, \, i\in\I_{\theta};
\\
\lambda_{i} &= 0 \text { if either } \un{\chi}_i \neq \varepsilon \text { or } g_i^2 = \kappa, \, i\in\I_{\theta};
\\
\lambda_{ij} &= 0 \text { if } i<j\le \ell, \, a_{ij}=-1, \, \chi_{ij} \neq \varepsilon;
\\
\lambda_{ij} &= 0 \text { if } i<j\le \ell, \, a_{ij}=0, \, \chi_{ij} \neq \varepsilon;
\\
\lambda_{ij}' &= 0 \text { if } i<j\le \ell, \, a_{ij}=0, \, \chi_{ij} \neq \varepsilon;
\\
\lambda_{ij} &= \lambda_{ij}' \text { if } i<j\le \ell, \, a_{ij}=0, \, G^{g_i}\ne G^{g_j};
\\
\lambda_{ij} &= 0 \text { if } i\le \ell<j, \, a_{ij}=0, \, \chi_{ij} \neq \varepsilon;
\\
\lambda_{ijk} &= 0 \text { if } i<j<k, \, a_{ji}=a_{jk}=-1, \, \un{\chi}_{j}\chi_{ik} \neq \varepsilon;
\\
\lambda_{ji} &= 0 \text { if } i<j, \, a_{ij}=-2, \, \un{\chi}_i\chi_{j}^2 \neq \varepsilon;
\\
\mu_{\beta} &= 0 \text { if either } \un{\chi}_{\beta} \neq \varepsilon \text{ or } g_{\beta}^2=1 \, (\beta \in \Delta_+^{V}-\{\alpha_i\});
\\
\mu_{\beta}' &= 0 \text { if either } \un{\chi}_{\beta} \neq \varepsilon \text{ or } g_{\beta}^2=\kappa \, (\beta \in \Delta_+^{V}-\{\alpha_i\}, \, \dim M_{\beta}=2).
\end{aligned}
\end{align}

In this subsection we prove our last main result:
\begin{theorem}\label{thm:liftings-general} 
Let $M \in \ydG$ of type $\alpha_{\theta}$, $\gamma_{\theta}$, $\delta_{\theta}$, 
$\epsilon_{\theta}$ or $\phi_4$.
For each $\bsl \in \cR_{V}$, see \eqref{eq:lambda-conditions-general}, let $\cL(\bsl)$ be the quotient of $T(M)\#\Bbbk G$, where we change the labels of the generators to $(z_i)_{i\in\I_{\theta+\ell}}$ by the following relations
\begin{align*}
& z_i^2 - \mu_i ( 1 - g_i^2), \qquad i\in\I_{\theta}; 
\\ 
& z_{i\un{i}}-\lambda_i(1-g_i^2\kappa), \qquad i\in\I_{\ell};
\\
& z_{\un{i}\un{j}}+\chi_i(\kappa)z_{ij}- \lambda_{ij} (1 - g_ig_j ), \qquad i<j\le \ell,  \, a_{ij}=-1;
\\
& z_{ij}-\lambda_{ij}(1-g_ig_j), z_{i\un{j}}-\lambda_{ij}'(1-g_ig_j\kappa), \qquad i<j\le \ell, \, a_{ij}=0; 
\\
& z_{ij}-\lambda_{ij}(1-g_ig_j), \qquad i\le \ell<j, \, a_{ij}=0;
\\
& (\ad_c z_j) z_{ij\un{i}}-\left(2\right)_{\chi_j(g_i^2\kappa)} \mu_j\lambda_{i} ( 1 - g_j^2) g_i^2+ \mu_{i} \lambda_j(1-g_j^2\kappa) g_i^2
-\lambda_{ji}(1-g_i^2g_j^2), \qquad a_{ij}=-2;
\\
& (\ad_c z_j) z_{ijk}-\left(2\right)_{\chi_j(g_ig_k\kappa)} \mu_j\lambda_{ik}' ( 1 - g_j^2) g_ig_k\kappa + \lambda_{ik} \lambda_j(1-g_j^2\kappa) g_ig_k
-\lambda_{ijk}(1-g_ig_j^2g_k\kappa), 
\\ 
& \qquad \qquad i<j<k, \, a_{ji}=a_{jk}=-1;
\\
& z_{\beta}^2-\mathbf{z}_{\beta}-\mu_{\beta}(1-g_{\beta}^2), \qquad \beta \in \Delta_+^{M}-\{\alpha_i\};
\\
& \left[z_{\beta},z_{\un{\beta}}\right]_c-\mathbf{z}_{\beta}-\mu_{\beta}'(1-g_{\beta}^2\kappa), \qquad \beta \in \Delta_+^{M}-\{\alpha_i\}, \dim M_{\beta}=2,
\end{align*}
where $\mathbf{z}_{\beta}$, $\mathbf{z}_{\beta}'\in T(M)\#\Bbbk G$ are defined recursively on $\beta \in \Delta_+^{M}-\{\alpha_i\}$ such that $z_{\beta}^2-\mathbf{z}_{\beta}$ is $(g_{\beta}^2,1)$-primitive  and $\left[z_{\beta},z_{\un{\beta}}\right]_c-\mathbf{z}'_{\beta}$ is $(g_{\beta}^2\kappa,1)$-primitive in the quotient of $T(M)\#\Bbbk G$ by the previous relations. Then:

\begin{enumerate}[leftmargin=*,label=\rm{(\alph*)}]
\item\label{item:liftings-general-deformed-PRV} $\cL(\bsl)\simeq L(\cA(\bsl),\toba(M)\#\Bbbk G)$.
\item\label{item:liftings-general-lifting} $\cL(\bsl)$ is a lifting of $M$ over $\Bbbk G$.
\item\label{item:liftings-general-cocycle} $\cL(\bsl)$ is a cocycle deformation of $\toba(M)\#\Bbbk G$.
\end{enumerate}

Conversely, if $L$ is lifting of $M$ over $\Bbbk G$, then there exist $\bsl \in \cR_{M}$ 
such that $L \simeq \cL_5(\bsl)$.
\end{theorem}

Let $M \in\ydG$ of type $\alpha_{\theta}$, $\gamma_{\theta}$, $\delta_{\theta}$, $\epsilon_{\theta}$ or $\phi_4$.
We choose first a stratification $\cG= \sqcup_{i=0}^4 \cG_i$ on the set of defining relations found in Theorem \ref{thm:Nichols-presentation-PBW-basis}
\begin{align*}
\cG_0=& \{x_i^2, \, x_{\un{i}}^2, \, x_{i\un{i}}|i\in\I_{\ell} \}\cup \{x_i^2 | i>\ell\};
\\
\cG_1=& \{x_{i\un{j}}+\chi_{j}(g_i^2) x_{\un{i}j}, \,  x_{\un{i}\un{j}}+\chi_i(\kappa) x_{ij}| i<j\le \ell,  \, a_{ij}=-1 \};
\\
\cG_2=& \{x_{ij}, \, x_{i\un{j}}, \, x_{\un{i}j}, \, x_{\un{i}\un{j}} | i<j\le \ell, \, a_{ij}=0 \} \cup \{x_{ij}, \, x_{\un{i}j} | i\le \ell<j, \, a_{ij}=0\};
\\
\cG_3 =& \{ r_{ijk}:=(\ad_c x_j) x_{ijk}, \, \un{r}_{ijk}:=(\ad_c x_j) x_{\un{i}jk} | i<j<k, \, a_{ji}=a_{jk}=-1 \}
\\
& \qquad \cup \{(\ad_c x_j) x_{ij\un{i}} | a_{ij}=-2\};
\\
\cG_4 =& \{x_{\beta}^2 | \beta \in \Delta_+^{V}-\{\alpha_i\}\} \cup  
\{[x_{\beta},x_{\un{\beta}}]_c | \beta \in \Delta_+^{V}-\{\alpha_i\}, \, \dim M_{\beta}=2\}.
\end{align*}
This stratification is good since
\begin{align*}
\Delta(r_{ijk}) =&r_{ijk}\ot 1 -\chi_j(g_ig_k\kappa) x_j^2 g_ig_k\kappa \ot x_{i\un{k}} -\chi_j(g_i^2g_k^2) x_j^2 g_ig_k\kappa \ot x_{\un{i}k}
\\ & + x_jx_{\un{j}} g_ig_k \ot x_{ik} + \chi_j(g_ig_k\kappa) x_{\un{j}}x_j g_ig_k \ot x_{\un{i}\un{k}}+g_ig_j^2g_k \kappa \ot r_{ijk},
\\
\Delta(\un{r}_{ijk}) =&\un{r}_{ijk}\ot 1 -\chi_j(g_ig_k\kappa) x_j^2 g_ig_k \ot x_{\un{i}\un{k}} -\chi_j(g_i^2g_k^2) x_j^2 g_ig_k \ot x_{ik}
\\ & + x_jx_{\un{j}} g_ig_k\kappa \ot x_{\un{i}k} + \chi_j(g_ig_k\kappa) x_{\un{j}}x_j g_ig_k\kappa \ot x_{i\un{k}}+g_ig_j^2g_k \ot \un{r}_{ijk}.
\end{align*}
Set $\mH_k:=\toba_k\# \Bbbk G$. The Yetter-Drinfeld structure of each stratum is given in Remark \ref{rem:YD-structure-relations}.

Let $\bsl\in \cR_{M}$. Define $\cE_0(\bsl) =\toba_0=T(M)$, but we change the labels of the generators to $(y_i, y_{\un{i}})_{i\in\I}$ to differentiate from the generators $(x_i,x_{\un{i}})_{i\in\I}$ of the pre-Nichols algebras $\toba_k$. Let
{\small
\begin{align*}
\cE_1(\bsl) & := \cE_0(\bsl) / \left\langle y_i^2 - \mu_i, y_{\un{i}}^2 - \mu_i, y_{i\un{i}} - \lambda_i \colon i\in\I_{\theta}\right\rangle,
\\
\cE_2(\bsl) & := \cE_1(\bsl) / \left\langle y_{i\un{j}}+\chi_j(g_i^2) y_{\un{i}j}- \lambda_{ij}, \, y_{\un{i}\un{j}}+\chi_i(\kappa) y_{ij}- \lambda_{ij} | i<j\le \ell,  \, a_{ij}=-1 \right\rangle, 
\\
\cE_3(\bsl) & := \cE_2(\bsl) / \left\langle 
\begin{array}{c}
y_{ij}-\lambda_{ij}, \, y_{i\un{j}}-\lambda_{ij}', \, y_{\un{i}j}-\lambda_{ij}', \, y_{\un{i}\un{j}}-\lambda_{ij}, \, i<j\le \ell, \, a_{ij}=0; \\
y_{ij}-\lambda_{ij}, \, y_{\un{i}j}-\lambda_{ij}, \, i\le \ell<j, \, a_{ij}=0
\end{array} \right\rangle, 
\\
\cE_4(\bsl) & := \cE_3(\bsl) / \left\langle 
\begin{array}{c}
(\ad_c y_j) y_{ijk}-\lambda_{ijk}, \, (\ad_c y_j) y_{\un{i}jk}-\lambda_{ijk}, \, a_{ji}=a_{jk}=-1,
\\
(\ad_c y_j) y_{ij\un{i}}-\lambda_{ji}, \, a_{ij}=-2
\end{array}
\right\rangle, 
\\
\cE_5(\bsl) & := \cE_4(\bsl) / \left\langle 
\begin{array}{c}
y_{\beta}^2-\mu_{\beta}, \, \beta \in \Delta_+^{V}-\{\alpha_i\},
\\
\left[y_{\beta},y_{\un{\beta}}\right]_c -\mu_{\beta}', \, \beta \in \Delta_+^{V}-\{\alpha_i\}, \, \dim M_{\beta}=2
\end{array}
\right\rangle.
\end{align*} 
}
Each $\cE_{i}(\bsl)$ is a $\Bbbk G$-module algebra since each defining ideal above is stable under the $G$-action by \eqref{eq:lambda-conditions-general}. Thus we may introduce $\cA_i(\bsl):=\cE_i(\bsl) \# \Bbbk G$.

\begin{lemma}\label{lem:cleft-objects-general}
Let $k\in\I_{5}$. Then $\cE_k(\bsl)\neq 0$ and each $\cA_k(\bsl)$ is an $\mH_k$-cleft object. There exists an $\mH_k$-colinear section $\gamma_k \colon \mH_k \to \cA_k$ that restricts to an algebra map $(\gamma_k) _{|\Bbbk G} \in \Alg (\Bbbk G, \cA_k)$.
\end{lemma}
\pf
Fix $\bsl\in\cR_M$; again, we simplify the notation and write $\cE_k=\cE_k(\bsl)$, $\cA_k=\cA_k(\bsl)$.
The proof is analogous to that of Lemma \ref{lem:cleft-objects-alpha2}, recursively on $k$. 

\smallbreak

When $k<5$, the key step is to prove that $\cE_k\ne 0$, which implies that $\cA_k\ne 0$: 
if so, then \cite[Theorem 8]{Gu} applies again to conclude that $\cA_k$ is 
$\mH_k$-cleft, hence there exists a section $\gamma_k$ as in the statement by \cite[Proposition 5.8]{AAGMV}.

To show that $\cE_k\ne 0$, it is enough to verify non vanishing when we \emph{deform just one submodule of relations}; 
that is, to consider the case $\bsl=\bsl^{(i)}$ for each $i$ and then proceed as in \cite[Lemma 3.4]{AG}. Indeed, if $\varpi_k:\cA_k(\bsl^{(i)})=\mH_{k-1}\twoheadrightarrow \cE_k(\bsl^{(I)})\# \Bbbk G$ is the canonical projection,
then the composition of the algebra map $(\rho_{k-1})_{|\cE_{k-1}}:\cE_{k-1}\to \cA_{k-1}\ot\mH_{k-1}$ with $(\id\ot\varpi_k)$ factors through $\cE_k=\cE_k(\bsl)$.

To verify that $\cE_k(\bsl^{(i)})\neq 0$ when the submodule of relations to be deformed is neither $\{(\ad_c x_j) x_{ijk}, \, (\ad_c x_j) x_{\un{i}jk}\}$, where $a_{ji}=a_{jk}=-1$, nor $\{(\ad_c x_j) x_{ij\un{i}}\}$ with $a_{ij}=-2$ we may use Lemma \ref{lem:cleft-objects-alpha2}.
For these two exceptions, we adapt the argument given in Lemma  \ref{lem:cleft-objects-alpha2} for relations $y_{1\un{2}}+\chi_2(g_1^2) y_{\un{1}2}$ and $y_{\un{1}\un{2}}+\chi_1(\kappa) y_{12}$, then use the cocycle $\sigma$ to reduce to deformations of Nichols algebras of diagonal type, so the result follows by \cite[Proposition 4.2]{AAG}.

\smallbreak

Finally, for $k=5$ we have that $\mH_4^{\co \pi_4}=\Zc(M)$ is a skew-polynomial algebra in variables 
$x_{\beta}^2$, $[x_{\beta},x_{\un{\beta}}]_c$ by Proposition \ref{prop:Z-polynomial-generators}, thus \cite[Theorem 4]{Gu} assures that
$\cA_5$ is $\mH_5$-cleft. The section $\gamma_5$ can be chosen so that $(\gamma_k) _{|\Bbbk G} \in \Alg (\Bbbk G, \cA_k)$ by \cite[Proposition 5.8]{AAGMV}.
\epf

We are ready to prove the main theorem of this section.

\pf[Proof of Theorem \ref{thm:liftings-general}]
We proceed as in \cite[Theorem 5.6]{AnS}, using correspondingly Theorem \ref{thm:trivial-hom-rels-V} and Lemma \ref{lem:cleft-objects-general}. Indeed, starting with $\cL_0(\bsl)=\mH_0$, we define succesive quotients $\cL_i(\bsl)$, $i\in \I_5$, where $\cL_4$ is the quotient by all the relations of $\cL(\bsl)$ except the last two sets (parametrized by $\beta\in\Delta_+$). Each $\cL_k(\bsl)$, $k<5$, is a Hopf algebra since one if obtained from the previous one by recursively quotient by skew-primitive elements. Working as in \cite[Theorem 1.6]{AAG}, $\cL_4(\bsl)\simeq L(\cA_4(\bsl),\mH_4)$, and there exist $\mathbf{z}_{\beta}\in \cL_4(\bsl)$ and $\mathbf{z}_{\beta}'\in \cL_4(\bsl)$ as stated below; moreover $\cL(\bsl)\simeq L(\cA_5(\bsl),\mH_5)$. The proof that these are all the liftings follows exactly as in \cite[Theorem 5.6]{AnS}\epf

\subsection{Foldings of liftings}

The folding construction in \cite[Part 1]{Len12}  was formulated in the following general setting: Let $H$ be a Hopf algebra and $H_\sigma,\;\sigma\in\hat{\Sigma}$ a group of biGalois objects with coherent choice of isomorphisms $\iota_{\sigma,\tau}:H_{\sigma\tau}\cong H_\sigma \Box H_\tau$.
By \cite[Theorem 1.6]{Len12}  the direct sum of algebras
$$\tilde{H}:=\bigoplus \sum_{\sigma\in\hat{\Sigma}} H_\sigma$$
can be endowed with the structure of a Hopf algebra with coproduct $\bigoplus_{\sigma,\tau}\iota_{\sigma,\tau}$. 

Conversely by \cite[Theorem 3.6]{Len12}, any Hopf algebra $\tilde{H}$ with $\Sigma$ a central subgroup is a 
folding of $H=\tilde{H}/\Sigma^+\tilde{H}$ by $\Sigma$. The biGalois objects are quotients of $\tilde{H}$ 
associated to a central character on $\Sigma$. The folding data in Section \ref{sec:folding} was formulated 
specifically for the situation $H=\toba(M)\#\Bbbk\Gamma$ and for biGalois objects arising from  $2$-cocycles 
$\sigma$ on the group $\Gamma$, trivially extended to $H$, and twisted Yetter-Drinfeld isomorphisms 
$\mathbf{u}:\toba(M)_\sigma\to \toba(M)$, extended by the identity on $G$ to $H$. In Theorem      
\ref{thm_NicholsFolding} we have stated the folding solely in terms of $\sigma,\mathbf{u}$, while in Theorem 
\ref{thm:folding} we have stated the folding with these specific choices of biGalois objects as above.

\smallbreak
We now discuss the following alternative systematic way to understand the liftings of folded Nichols algebras,  which we constructed in the previous section: Let $H'$ be a lifting of $H=\toba(M)\#\Bbbk\Gamma$ for a diagonal Nichols algebra, which are classified in \cite{AAG,AnG}. Let again $(H'_\sigma)_{\sigma\in \hat{\Sigma}}$ be a group of biGalois objects over the lifting, then we have a folding $\tilde{H'}$, whose graded algebra is the folding $\tilde{H}=\toba(\tilde{M})\#\Bbbk G$ of $H$. One source for such biGalois objects could be again folding data $(\sigma,\mathbf{u})$ where in addition $\mathbb{u}$ is compatible with the lifting $H'$, more precisely, leaves a lifting cocycle invariant. But there are also other possibilities, namely the $2$-cocycle $\sigma$ over $\Gamma$ could be nontrivially extended to $\toba(M)\#\Bbbk\Gamma$, which would cause a folding of $H$ that is a lifting of $\tilde{H}$ with values in the new center.

Conversely we obtain in this way all liftings of $\tilde{H}$ where $\Sigma$ is central, acting trivially on $\tilde{M}$. We have already shown for each Nichols algebra in Theorem \ref{thm:Doi-Twist}, that this trivial action can always be achieved by a Doi twist; however, it is not a-priori clear that these Doi twists carry over to the lifting. We will now use this tool to analyze the smallest example:

\begin{exa}[Case ${^2}A_2^2$]
We consider the Nichols algebra of type $\alpha_2$ defined in Section \ref{subsubsec:type-ADE} over a group $G$ generated by $g_1,g_2$ with $g_1g_2=\kappa g_2g_1$ and $\kappa$ central of order two, which is a central extension of the abelian group $\Gamma$ generated by $\bar{g}_1,\bar{g}_2$. We computed all its liftings in Section \ref{subsec:liftings-alpha2}. We can conveniently take the group action from Remark \ref{rem:YD-structure-relations} and the braiding from Section \ref{subsubsec:type-ADE} to replace $z_{ij}$ again by the braided commutator. For example, the first three relations depending on the parameters $\mu_i,\lambda_i,\lambda_{12}$ read
\begin{align*}
&z_i^2 = \mu_i ( 1 - g_i^2),
\qquad \qquad 
z_{i}z_{\un{i}}+\chi_i(\kappa)z_{\un{i}}z_i = \lambda_i(1-g_i^2\kappa),
\\
& z_{\un{1}}z_{\un{2}}-\chi_2(\kappa g_1^2)z_{2}z_{\un{1}}
+\chi_1(\kappa)z_1z_2-\chi_1(\kappa)z_{\un{2}}z_{1}
= \lambda_{12} (1 - g_1g_2 ).
\end{align*}
Note that acting with a group element on a relation may produce more relations, as we explained in the proof of Lemma \ref{lem:cleft-objects-alpha2} for the cleft objects.
For example, acting with $g_j$ on the first relation produces the relation $z_{\un{i}}^2 = \mu_i ( 1 - g_i^2)$.

The associated Nichols algebra and its liftings are foldings iff $\kappa$ is a central element in the Hopf algebra, i.e. $\chi_1(\kappa)=\chi_2(\kappa)=1$ (which we saw it is true up to Doi twist). In this case we saw in Remark \ref{rem:diagonalizeBraiding} that the braiding diagonalizes in the basis
\begin{align*}
\xt_{i} &= z_i + q_{ij} z_{\un{i}}, &
\xt_{\un{i}} &= z_i - q_{ij} z_{\un{i}}.
\end{align*}
The elements $\xt_{i}$, $\xt_{\un{i}}$ are not $G$-homogeneous, but $\Gamma$-homogeneous with degrees $\bar{g}_i$, $i\in \I_2$. On the other hand they are $G$-eigenvectors with $g_i,g_j$ acting on $\xt_{i}$ with eigenvalues $-1,-q_{ji}$ and on $\xt_{\un{i}}$ with eigenvalues $-1,q_{ji}$. 
The diagonal braiding matrix is of type $A_2\times A_2$ 
$$\begin{pmatrix}
-1 & -q_{12} & -1 & q_{12} \\
-q_{21} & -1 & q_{21} & -1 \\
-1 & -q_{12} & -1 & q_{12} \\
-q_{21} & -1 & q_{21} & -1 \\
\end{pmatrix}$$
and a twisted symmetry switching the two copies. In the folding construction $z_i,z_{\un{i}}$ arise as eigenvalues of this symmetry. 

\smallbreak
We now rewrite the relations in this basis, starting with those involving just one orbit $\xt_i,\xt_{\un{i}}$, which is a diagonal Nichols algebra of type $A_1\times A_1$:
\begin{align*}
&\tfrac{1}{4}(\xt_{i}+\xt_{\un{i}})^2= \mu_i(1 - g_i^2), \qquad 
\tfrac{q_{ij}^{-2}}{4}(\xt_{i}-\xt_{\un{i}})^2 = \mu_i ( 1 - g_i^2),
\\
& \tfrac{q_{ij}^{-1}}{4}((\xt_{i}+\xt_{\un{i}})(\xt_{i}-\xt_{\un{i}})+(\xt_{i}-\xt_{\un{i}})(\xt_{i}+\xt_{\un{i}})) = \lambda_i(1-g_i^2\kappa).
\end{align*}
These relations rewrite to 
\begin{align*}
&\xt_{i}^2= (1+q_{ij}^2)\mu_i(1 - g_i^2)+q_{ij}\lambda_i(1-g_i^2\kappa), &
&\xt_{i}\xt_{\un{i}}+\xt_{\un{i}}\xt_{i}=2(1-q_{ij}^2)\mu_i(1 - g_i^2)=0,
\\
&\xt_{\un{i}}^2= (1+q_{ij}^2)\mu_i(1 - g_i^2)-q_{ij}\lambda_i(1-g_i^2\kappa), &
&\tfrac{1}{2}(\xt_{i}^2-\xt_{\un{i}}^2) = q_{ij}\lambda_i(1-g_i^2\kappa),
\end{align*}
where we have to take into account that Section \ref{subsec:liftings-alpha2} states that $\mu_i\neq 0$, respectively $\lambda_i\neq 0$, only if $q_{ji}^2=1$, so the anti-commutator vanishes.
 
This is consistent with the possible liftings of diagonal $A_1\times A_1$: 
\begin{itemize}[leftmargin=*]
\item The anti-commutator relation admits a non-trivial lifting if $\chi_{i}\chi_{\un{i}}=\epsilon$, but in our case $\chi_{i}(g_j)\chi_{\un{i}}(g_j)=-1$.
\item The truncation relations admit non-trivial liftings if $\chi_i^2=\epsilon$, which is the case if and only if $q_{ji}^2=1$. If the respective lifting parameters are equal, then this lifting datum is compatible with a folding using the group $2$-cocycle. This produces the symmetric lifting depending on $\mu_i$. 
\item On the other hand the antisymmetric lifting depending on $\lambda_i$ requires a lifting cocycle that is non-trivially extended to the Nichols algebra. The corresponding non-trivial biGalois object is determined by plugging the non-trivial central character $\kappa\mapsto -1$.
\end{itemize}

We now turn to the relation involving $\lambda_{12}$:
\begin{align*}
\lambda_{12} (1 - g_1g_2 )&=\tfrac{q_{12}^{-1}}{2}(\xt_{1}-\xt_{\un{1}})
\tfrac{q_{21}^{-1}}{2}(\xt_{2}-\xt_{\un{2}})
-q_{12}^2
\tfrac{1}{2}(\xt_{2}+\xt_{\un{2}})
\tfrac{q_{12}^{-1}}{2}(\xt_{1}-\xt_{\un{1}})\\
&+\tfrac{1}{2}(\xt_{1}+\xt_{\un{1}})
\tfrac{1}{2}(\xt_{2}+\xt_{\un{2}})
-\tfrac{q_{21}^{-1}}{2}(\xt_{2}-\xt_{\un{2}})
\tfrac{1}{2}(\xt_{1}+\xt_{\un{1}})
\\
&=\tfrac{1}{2}(\xt_{1}\xt_{\un{2}}-q_{12}\xt_{\un{2}}\xt_1)+\tfrac{1}{2}(\xt_{\un{1}}\xt_{{2}}+q_{12}\xt_2\xt_{\un{1}})
\lambda_{12} (1 - g_1g_2\kappa)
\\
&=-\tfrac{q_{12}}{2}(\xt_{1}\xt_{\un{2}}-q_{12}\xt_{\un{2}}\xt_1)
+\tfrac{q_{12}}{2}(\xt_{\un{1}}\xt_{{2}}+q_{12}\xt_2\xt_{\un{1}}).
\end{align*}
Section \ref{subsec:liftings-alpha2} with $\chi_{12}$ in \eqref{eq:centralizer-g1g2-character} applied to $g_i^2$ and $g_1g_2$ states that $\lambda_{12}\neq 0$ only if $1=q_{ii}^2q_{ij}^2$  and $1=q_{11}q_{21}^2\cdot q_{22}q_{12}^2$, which is again equivalent to $q_{ij}^2=q_{ji}^2=1$. Possibly reversing $1,2$ we may assume we are in the case $q_{12}=1,q_{21}=-1$, then adding and subtracting the previous relations returns:
\begin{align*}
\xt_{\un{1}}\xt_{{2}}+\xt_2\xt_{\un{1}}&=2\lambda_{12}\left(1-g_1g_2\tfrac{\kappa+1}{2}\right), &
\xt_{{1}}\xt_{\un{2}}-\xt_{\un{2}}\xt_{{1}}&=2\lambda_{12}\left(1-g_1g_2\tfrac{\kappa-1}{2}\right).
\end{align*}
On the other hand, the diagonal Nichols algebra $A_2\times A_2$ has such liftings of
\begin{itemize}
    \item $\xt_{\un{1}}\xt_{{2}}-q_{ji}\xt_2\xt_{\un{1}}$ if $\chi_{\un{1}}\chi_{{2}}=\epsilon$, which is the case for $q_{12}=1,q_{21}=-1$.
    \item $\xt_{{1}}\xt_{\un{2}}+q_{ij}\xt_{\un{2}}\xt_{{1}}$ if  $\chi_{1}\chi_{\un{2}}=\epsilon$, which is the case for $q_{12}=-1,q_{21}=1$.
\end{itemize}
Altogether, there is no $\mathbf{u}$-symmetric lifting of this type, and the solution we find starts with a lifting $H'$ for one of these relation, again visible at the central character $\kappa\mapsto 1$, and the other of these relations appears in the biGalois object which is nontrivially extended from the group $2$-cocycle. 

We refrain from discussing the last two relations in a similar manner.
\end{exa}

\section{Future directions}\label{sec_outlook}
We conclude by some outlook questions that naturally arise from our work.

\begin{question}
Is there a modified folding construction that produces the remaining Nichols algebras in Heckenberger--Vendramin classification?
\end{question}

\begin{question}\label{q_decomposable}
Several folded Nichols algebras in \cite{Len14} do not appear in \cite{HV-rank>2} because their support is too small. More precisely, these are the cases ${^2}D_n$ and ${^3}D_4$ and $^{2}A_1^2$ familiar from Lie theory, as well as unfamiliar cases ${^2}A_2$ at a third root of unity and several cases involving other diagonal Nichols algebras.
We expect that our methods can be applied in these cases.
\end{question}

\begin{question}
Which modular tensor categories can be constructed from the new pointed Hopf algebras described here?

From the categorical perspective, there is a rather unique $\Sigma$-graded extensions of tensor categories with a $\Sigma$-crossed braiding \cite{ENOM}. Since the operations of $\Sigma$-graded extension and taking Hopf algebra representations commute, this extension could be computed by taking a $\Sigma$-symmmetric Nichols algebra over an abelian group, which can then be folded to Nichols algebra over the known $\Sigma$-extension of the abelian group. To get a braiding, this would require a non-trivial associator (an effect familiar for quantum group of even order root of unity), and for $\Sigma=\mathbb{Z}_2$ conjecturally involve a Tambara-Yamigami category. 
\end{question}

\begin{question}
The Logarithmic Kazhdan--Lusztig Correspondence, see e.g. \cite{FGST, Len17} conjectures the existence of a vertex algebra, realized as subalgebra of a free field algebra, whose tensor category of representations is equivalent to representations of a small quantum group. The folding construction and the previous problem suggests an extension of this construction, where the free field algebra is replaced by an orbifold model.     
\end{question}

\appendix

\section{Proof of Theorem \ref{thm:Doi-Twist}}\label{sec:appendix}
Here we complete the proof of Theorem \ref{thm:Doi-Twist}, which states that Nichols algebras of types $\alpha_{\theta}$, $\gamma_{\theta}$, $\delta_{\theta}$, $\epsilon_{\theta}$ and $\phi_4$  become of diagonal type when an appropriate twist is performed.
The cases that remain unsolved are $\alpha_2$, $\alpha_3$, $\gamma_3$, $\gamma_4$, $\delta_4$ and $\phi_4$, which will be dealt with in Proposition \ref{prop:DoiTwist}. The other cases are diagonal by Lemma \ref{lem:trivial-action-kappa-big-rank}.
\footnote{The corresponding folded Nichols algebras are those with symplectic root system \cite[\S 4]{Len14} of dimension $\theta$ and radical dimension $r=\theta-2$.}

\subsection{Group Cohomology Tools}
We start by collecting some useful group extensions and group cohomology statements for later use, following \cite[Chapter 7]{Len12}. Recall the definition of $\Gamma_{u,v,\kappa}$ in \eqref{eq:defn-gamma-uvkappa}.

\newcommand{\formalw}{\mathtt{t}}
\newcommand{\formalf}{\sqrt[r]{f}}

\begin{definition}
Let $\Gamma$ be an abelian group, $\formalw$ a generator of $\Z$. 
For each $w\in \Gamma$, $r\in\N$ we consider the abelian group 
$$ \Gamma(\sqrt[r]{w}):=\Gamma\times \Z /\langle (w,\formalw^{-r}) \rangle.$$
We shall identify $g\in\Gamma$, $\formalw^k$, $k\in\N$, with their images $\overline{(g,e)}$, $\overline{(e,\formalw^k)}$ in $\Gamma(\sqrt[r]{w})$. The defining relation becomes $w=\formalw^r$. We think of $\Gamma(\sqrt[r]{w})$ as the set $\{g\formalw^k | g\in\Gamma, 0\le k<r\}$ with product
\begin{align*}
g\formalw^j \cdot h\formalw^k =\begin{cases}
gh \formalw^{j+k}, & \text{ if }j+k<r, \\ ghw \formalw^{j+k-r}, & \text{ if }j+k\ge r, 
\end{cases}
&& g,h\in\Gamma, 0\le j,k<r.
\end{align*}
\end{definition}

\begin{obs}\label{rem:extensions-trivial-facts}
Fix $\Gamma$ an abelian group, $u,v,w,\kappa\in\Gamma$, where $\kappa^2=e$, $r\in\mathbb{N}$. 
\begin{enumerate}[leftmargin=*,label=\rm{(\roman*)}]
\item The inclusion $\Gamma\hookrightarrow \Gamma(\sqrt[r]{w})$ extends to an injective map
$$ \Gamma_{u,v,\kappa} \hookrightarrow \Gamma(\sqrt[r]{w})_{u,v,\kappa}.$$
\item Fix also $z\in\Gamma$, $s\in\N$. There is a canonical isomorphism
$\Gamma(\sqrt[r]{w})(\sqrt[s]{z}) \cong \Gamma(\sqrt[s]{z})(\sqrt[r]{w}),$ 
which in turn induces an isomorphism
$$\Gamma(\sqrt[r]{w})(\sqrt[s]{z})_{u,v,\kappa} \cong \Gamma(\sqrt[s]{z})(\sqrt[r]{w})_{u,v,\kappa}.$$
\end{enumerate}
\end{obs}

\begin{prop}\label{prop_cohomology_extension}
Let $\Gamma$ be an abelian group, $u,v,w,\kappa\in\Gamma$, where $\kappa^2=e$, $r\in\mathbb{N}$. 
Given $\sigma\in H^2(\Gamma_{u,v,\kappa},\Bbbk^\times)$, consider $f:=f_{\sigma}\colon \Gamma\to\Bbbk^{\times}$ given by
\begin{align*}
f(g) &=\sigma(g,w)/\sigma(w,g), & &g\in\Gamma.
\end{align*}
Then $\sigma$ lifts to a $2$-cocycle $\widetilde{\sigma}\in H^2 \big( \Gamma(\sqrt[r]{w})_{u,v,\kappa},\Bbbk^\times \big)$ if and only if there exists $\ft\in\widehat{\Gamma_{u,v,\kappa}}$ such that $\ft^r=f$, $\ft(w)=1$. In this case, 
\begin{align}\label{eq:cohom-extension-f}
\widetilde{\sigma}(g,\formalw)\widetilde{\sigma}^{-1}(\formalw,g)&=\ft(g), & \text{for all } &g\in\Gamma.
\end{align}
If $w=1$, then any $\sigma\in H^2(\Gamma_{u,v,\kappa},\Bbbk^\times)$ lifts to $\widetilde{\sigma}\in H^2 \big( \Gamma(\sqrt[r]{w})_{u,v,\kappa},\Bbbk^\times \big)$.
\end{prop}
\pf
Set $\underline{\Gamma}:=\Gamma_{u,v,\kappa}$, $G=\Gamma(\sqrt[r]{w})_{u,v,\kappa}$.
We can write $G$ as the following central extension:
$$ 1\to \langle w\formalw^{-r}\rangle \to \underline{\Gamma}\times \langle \formalw\rangle \to  G \to 1.$$
As $\Bbbk^\times$ is divisible, the map $H^1(\langle \formalw\rangle,\Bbbk^\times) \to H^1( \langle w \formalw^{-r}\rangle,\Bbbk^\times)$ is surjective and $H^2(\langle \formalw\rangle,\Bbbk^\times)= H^2(\langle w\formalw^{-r}\rangle,\Bbbk^\times)=1$. Hence the exact sequence in \cite[\S 1]{IM64} associated to the central extension below is
\begin{align*}
1 &\to H^1(G,\Bbbk^\times)\to H^1( \underline{\Gamma}\times \langle \formalw\rangle,\Bbbk^\times) \to H^1( \langle w\formalw^{-r}\rangle,\Bbbk^\times) \to\\
&\to H^2(G,\Bbbk^\times)\to H^2(\underline{\Gamma}\times \langle \formalw \rangle,\Bbbk^\times) \to
\mathrm{Pair}(\underline{\Gamma}\times \langle \formalw\rangle,\langle w\formalw^{-r}\rangle).
\end{align*}
Using the results above, the K\"unneth formula for the cohomologies of the direct product and decomposing the pairings we get 
\begin{align*}
1 & \longrightarrow  H^2(G,\Bbbk^\times)
\longrightarrow H^2(\underline{\Gamma},\Bbbk^\times)\times \mathrm{Pair}(\underline{\Gamma},\langle \formalw\rangle) \overset{\Phi}{\longrightarrow} 
\mathrm{Pair}(\underline{\Gamma},\langle w\formalw^{-r}\rangle) \times \mathrm{Pair}(\langle \formalw\rangle,\langle w\formalw^{-r}\rangle),
\end{align*}
where $\Phi$ is defined as follows:
\begin{itemize}[leftmargin=*]
\item for $\sigma\in H^2(\underline{\Gamma},\Bbbk^\times)$, we have $\Phi(\sigma)=(B_{\sigma},1)$, where $1\in \mathrm{Pair}(\langle \formalw\rangle,\langle w\formalw^{-r}\rangle)$ is the trivial pairing, and $B_{\sigma}\in\mathrm{Pair}(\underline{\Gamma},\langle w\formalw^{-r}\rangle)$ is given by 
\begin{align*}
B_{\sigma} (g,w\formalw^{-r}) &= \sigma(g,w)\sigma^{-1}(w,g), &
& g\in\underline{\Gamma}.
\end{align*}
\item for $P\in \mathrm{Pair}(\underline{\Gamma},\langle \formalw\rangle)$, we have $\Phi(P)=(F'_P, F''_P)$, where
\begin{align*}
F'(g,w\formalw^{-r}) &= P(g,\formalw)^{-r}, \, g\in\underline{\Gamma};
&
F''(\formalw,w\formalw^{-r}) &= P(w,\formalw). 
\end{align*}
\end{itemize}
Hence $(\sigma,P)\in\ker \Phi$ if and only if $\sigma(g,w)\sigma^{-1}(w,g)=P(g,\formalw)^{r}$ for all $g\in \underline{\Gamma}$ and $P(w,\formalw)=1$.
If $\sigma$ lifts to a $2$-cocycle $\widetilde{\sigma}\in H^2 \big( \Gamma(\sqrt[r]{w})_{u,v,\kappa},\Bbbk^\times \big)$, then set
$\ft$ as in \eqref{eq:cohom-extension-f}. Reciprocally, if there exists such $\ft$, we define $P\in \mathrm{Pair}(\underline{\Gamma},\langle \formalw\rangle)$, by
$P(g,\formalw)=\ft(g)$, and get $(\sigma,P)\in\ker\Phi$. By exactness of the sequence, $\ker\Phi$ is the image of the injective map 
$H^2(G,\Bbbk^\times)\longrightarrow H^2(\underline{\Gamma},\Bbbk^\times)\times \mathrm{Pair}(\underline{\Gamma},\langle \formalw\rangle)$; thus $\sigma$ lifts to a 2-cocycle $\widetilde{\sigma}\in H^2 \big( \Gamma(\sqrt[r]{w})_{u,v,\kappa},\Bbbk^\times \big)$. Moreover, $P$ describes the values on the additional generator. The last statement is clear.
\epf

\begin{coro}\label{cor_cohomology_doubleextension}
Let $\Gamma$ be an abelian group, $u,v,w,z,\kappa\in\Gamma$, where $\kappa^2=e$, $r,s\in\mathbb{N}$. 
Then $\sigma\in H^2(\Gamma_{u,v,\kappa},\Bbbk^\times)$ lifts to $\widetilde{\sigma}\in H^2 \big(\Gamma(\sqrt[r]{w})(\sqrt[s]{z})_{u,v,\kappa},\Bbbk^\times \big)$ if and only if there exists $\ft, \gt\in\widehat{\Gamma_{u,v,\kappa}}$ such that 
$\ft^r=f=\gt^s$, $\ft(w)=1=\gt(w)$.
\end{coro}
\pf
Use the isomorphism in Remark \ref{rem:extensions-trivial-facts} and Proposition \ref{prop_cohomology_extension}.
\epf

Next we assume that $\Gamma$ splits as $\Gamma=\Lambda\oplus\Omega$, where $u\in\Lambda$, $v,\kappa\in\Omega$. Recall the extension $\Gamma \hookrightarrow  \Gamma_{u,v,\kappa} \twoheadrightarrow \Z_2\times \Z_2=\langle x,y\rangle$ from \eqref{eq:defn-gamma-uvkappa}. Consider 
\begin{itemize}
\item $\underline{\Lambda}$ the subgroup of $\Gamma_{u,v,\kappa}$ generated by $\Lambda$ and $x$,
\item $\underline{\Omega}$ the subgroup of $\Gamma_{u,v,\kappa}$ generated by $\Omega$ and $y$.
\end{itemize}
Hence $\Gamma_{u,v,\kappa}\simeq \underline{\Lambda} \ltimes \underline{\Omega}$, where $x$ acts on $\underline{\Omega}$ by
\begin{align*}
x\cdot y & = \kappa y, & x \cdot h = h, & &h\in\Omega.
\end{align*}
Let $\widetilde{H}^2(\Gamma_{u,v,\kappa},\Bbbk^\times)$ denote the kernel of the restriction map 
$$ H^2(\Gamma_{u,v,\kappa},\Bbbk^\times)\to H^2(\underline{\Lambda},\Bbbk^\times) .$$
By \cite[Theorem 2 (I)]{Tah72}, we have that
\begin{align*}
H^2(\Gamma_{u,v,\kappa},\Bbbk^\times) \simeq H^2(\underline{\Lambda},\Bbbk^\times) \oplus \widetilde{H}^2(\Gamma_{u,v,\kappa},\Bbbk^\times).
\end{align*}
By \cite[Theorem 2 (II)]{Tah72}, there exists an exact sequence
\begin{align}\label{eq:exact-sequence-H2tilde}
0 & \longrightarrow H^1 \big(\underline{\Lambda}, \widehat{\underline{\Omega}} \big) \longrightarrow \widetilde{H}^2(\Gamma_{u,v,\kappa},\Bbbk^\times)
\overset{res}{\longrightarrow} H^2(\underline{\Omega},\Bbbk^\times)^{\underline{\Lambda}}.
\end{align}
The image of the first map is the subspace of $H^2(\Gamma_{u,v,\kappa},\Bbbk^\times)$ of $2$-cocycles which are cohomologically trivial on $\underline{\Lambda}$ and $\underline{\Omega}$. Next we will characterize $H^1 \big(\underline{\Lambda}, \widehat{\underline{\Omega}} \big)$ and describe the shape of $2$-cocycles coming from this group.

\begin{prop}\label{prop_cohomology_semidirect}
Let $\mathcal{T}$ denote the set of triples $(P,\chi,\psi)\in \mathrm{Pair}(\Lambda,\Omega)\times \widehat{\Lambda}\times\widehat{\Omega}$ such that
\begin{align*}
\psi(v)\psi(\kappa)&=\chi(u); 
& P(g,\kappa)&=1, & \chi(g)^2&=P(g,v), &\psi(h)^2&=P(u,h), &\text{for all }g&\in\Lambda,   h\in\Omega.
\end{align*}
\begin{enumerate}[leftmargin=*,label=\rm{(\alph*)}]
\item\label{item:cohomology_semidirect-a} The map $H^1 \big(\underline{\Lambda}, \widehat{\underline{\Omega}} \big)\to \mathcal{T}$, $\phi \mapsto \left( P_{\phi}, \chi_{\phi}, \phi(x)_{|\Omega} \right)$, where
\begin{align}\label{eq:cohomology_semidirect}
P_{\phi}(g,h) &=\phi(g)(h), & \chi_{\phi}(g) &=\phi(g)(y), & & g\in\Lambda,h\in\Omega,
\end{align}
is bijective.

\item\label{item:cohomology_semidirect-b} 
The image $\sigma\in\widetilde{H}^2(\Gamma_{u,v,\kappa},\Bbbk^\times) \subseteq H^2(\Gamma_{u,v,\kappa},\Bbbk^\times)$ of a triple $(P,\chi,\psi)\in \mathcal{T}$ (viewed as an element of $H^1 \big(\underline{\Lambda}, \widehat{\underline{\Omega}} \big)$) under the map $\partial$ in \eqref{eq:exact-sequence-H2tilde}  satisfies
\begin{align*}
\frac{\sigma(g,h)}{\sigma(h,g)}&=P(g,h), & \frac{\sigma(g,y)}{\sigma(y,g)}&=\chi(g), & \frac{\sigma(x,h)}{\sigma(h,x)}&=\psi(h),
\end{align*}
for all $g\in\Lambda$, $h\in\Omega$. In particular we have that 
\begin{align*}
\frac{\sigma(g,\kappa)}{\sigma(\kappa,g)}&=1 \text{ for all }g\in\Gamma, & \frac{\sigma(\kappa,y)}{\sigma(y,\kappa)}&=\chi(\kappa), & \frac{\sigma(x,\kappa)}{\sigma(\kappa,x)}&=\psi(\kappa).
\end{align*}
\end{enumerate}
\end{prop}
\pf
\ref{item:cohomology_semidirect-a}
As $\Lambda$ acts trivially on $\widehat{\underline{\Omega}}$, 
each crossed morphism $\phi\in H^1 \big(\underline{\Lambda}, \widehat{\underline{\Omega}} \big)$ restricts on $\Lambda$ to a homomorphism, and hence to a pairing $\Lambda \times \underline{\Omega} \to \Bbbk^\times$, which we think as a pair $(P, \chi)\in \mathrm{Pair}(\Lambda,\Omega)\times \widehat{\Omega}$ as in \eqref{eq:cohomology_semidirect} such that $\chi(g)^2=P(g,v)$ for all $g\in\Lambda$ (because $y^2=v$). 
We set $\chi_{\phi}:=\phi(x)_{|\Omega}:\Omega\to\Bbbk$. As $x$ acts trivially on $\Omega$, $\chi$ is a group homomorphism. 
Hence we have an injective map
$$ H^1 \big(\underline{\Lambda}, \widehat{\underline{\Omega}} \big)\to \{ (P,\chi,\psi)\in \mathrm{Pair}(\Lambda,\Omega)\times \widehat{\Lambda}\times\widehat{\Omega}\ | \ \chi(g)^2=P(g,v), \text{ for all } g\in\Lambda\}.$$
If $\xi:=\phi(x)(y)\in\Bbbk$, then $\xi^2=\phi(x)(v)=\chi(v)$ and 
\begin{align}\label{eq:cohomology_semidirect-inverse-map}
\phi(gx^i)(hy^j) & =P(g,h)\psi(g)^j\chi(h)^i\xi^{ij}, & & g\in\Lambda,h\in\Omega, i,j\in\{0,1\}.
\end{align}

Reciprocally, given a triple $(P,\chi,\psi)$ as above, set $\phi:\underline{\Lambda}\to \widehat{\underline{\Omega}}$ as in \ref{eq:cohomology_semidirect-inverse-map}. Then $\phi$ is a crossed homomorphism if and only if for all $g\in\Lambda$, $h\in\Omega$,
\begin{align*}
\phi(gx)(hy)&=\phi(xg)(hy)=\phi(x)(hy)(x \cdot \phi(g)(hy))=P(g,\kappa h)\psi(g)\chi(h)\xi,
\\
\phi(u)(h) &=\phi(x^2)(h)=\phi(x)(h)(x \cdot \phi(x)(h))=\chi(h)^2,
\\
\phi(u)(hy) &=\phi(x^2)(hy)=\phi(x)(hy)(x \cdot \phi(x)(hy))=\chi(h)^2\xi^2\chi(\kappa).
\end{align*}
This means $P(g,\kappa)=1$ for all $g\in\Lambda$, $\chi(h)^2=P(u,h)$ and $P(u,h)\psi(u)=P(u,h)\chi(v)\chi(\kappa)$ for all $h\in\Omega$.

\noindent \ref{item:cohomology_semidirect-b}
This follows by explicit computation of the coboundary map $\partial$, see for example the proof of \cite[Theorem 2 (II)]{Tah72}.
\epf

The next result will allow us to reduce the question about the existence of a 2-cocycle just for groups of order a power of $2$.

\begin{prop}\label{prop_2group}
Let $\Gamma=\Gamma_2\times \Gamma_{odd}$, where $|\Gamma_2|=2^n$ for some $n\in\N$, and $|\Gamma_{odd}|$ is odd. Let $u=u_2\ut$, $u=v_2\vt$, with $u_2,v_2\in\Gamma_{2}$, $\ut,\vt\in\Gamma_{odd}$. Then
\begin{align*}
\Gamma_{u,v,\kappa} &\simeq (\Gamma_2)_{u_2,v_2,\kappa}\times \Gamma_{odd},
&
H^2(\Gamma_{u,v,\kappa},\Bbbk^\times) & \simeq H^2((\Gamma_2)_{u_2,v_2},\Bbbk^\times)\times H^2(\Gamma_{odd},\Bbbk^\times).
\end{align*}
\end{prop}
\pf
Let $m,n\in\N_0$ be such that $|\ut|=2m+1$, $|\vt|=2n+1$. We write $\tilde{x}$, $\tilde{y}$ for the extra generators of $(\Gamma_2)_{u_2,v_2,\kappa}$ and keep $x,y$ for those in $\Gamma_{u,v,\kappa}$. Then
\begin{align*}
\Gamma_{u,v,\kappa} &\to (\Gamma_2)_{u_2,v_2,\kappa}\times \Gamma_{odd}, 
\\
g_2\gt x^iy^j & \mapsto g_2\ut^{i(m-1)}\vt^{j(n-1)}\tilde{x}^i \tilde{y}^j\gt, & & g_2\in\Gamma_2,\gt\in\Gamma_{odd}, i,j\in\{0,1\},
\end{align*}
is a group isomorphism. 
The isomorphism between cohomology groups follows by K\"unneth's formula.
\epf

\subsection{Non-abelian groups and 2-cocycles}

Next we discuss how to apply Propositions \ref{prop_cohomology_extension} and \ref{prop_cohomology_semidirect} to the main classes of groups $\Gamma_{u,v,\kappa}$ that will appear in the proof below in order to obtain the desired $2$-cocycles. 

\subsubsection{}\label{subsubsec:example1} If $\Gamma=\Z_2=\langle\kt\rangle$ then $\Gamma_{1,1,\kt}$, $\Gamma_{1,\kt,\kt}$, $\Gamma_{\kt,1,\kt}$ are isomorphic to the dihedral group of order $8$. 
\footnote{On the other hand, $\Gamma_{\kt,\kt,\kt}$ is the quaternion group, which has trivial cohomology.}

\begin{enumerate}[leftmargin=*,label=\rm{\alph*)}]
\item\label{item:example1-a} We apply Proposition \ref{prop_cohomology_semidirect} to $\Gamma_{1,\kt,\kt}$ with $\Lambda=1$, $\Omega=\Gamma$, $P=1$, $\chi=1$, $\psi(\kt)=-1$; then $H^2(\Gamma_{1,\kt,\kt},\Bbbk^\times)\simeq \Z_2=\langle\sigma\rangle$, where 
\begin{align*}
\sigma(x,\kt)\sigma^{-1}(\kt,x)&=-1,  & \sigma(y,\kt)\sigma^{-1}(\kt,y)&=1. 
\end{align*}

\smallskip

\item\label{item:example1-b} We apply Proposition \ref{prop_cohomology_extension} to $\Gamma_{1,\kt,\kt}$, where $f(x)=-1$, $f(y)=1$. If either $w=1$ or $2\nmid r$, then $\sigma$ can be lifted to $\Gamma(\sqrt[r]{w})$; but when $w=\kt$, $2\mid r$, the lift does not exist. As smallest example, $\Gamma(\sqrt[2]{\kt})$ is the almost extraspecial group $2^{3+1}$, which has cohomology $\Z_2^2$: all $2$-cocycles are lifts of the trivial $2$-cocycle on $\Gamma_{1,\kt,\kt}$ with $f=1$ and $\ft=\pm1$.

\smallskip

\item\label{item:example1-c} More generally, if there exists a surjective map $\pi:\Gamma\to \Z_2$ such that $\pi(u)=1$, $\pi(v)=\kt=\pi(\kappa)$, then the pullback of $\sigma$ from $(\Z_2)_{1,\kt,\kt}$ to $\Gamma_{u,v,\kappa}$ can be lifted to $\Gamma(\sqrt[r]{w})$ either when $\pi(w)=1$ or $2\nmid r$. For $\pi(w)=\kappa$, $2\mid r$, it can also be lifted if there is a character $\ft:\Gamma\to \Bbbk^\times$ such that 
\begin{align*}
\ft(u) &\in \G_r, & 2\nmid & r/\mathrm{ord}(\ft(u)), & 
\ft(v) &\in \G_{r/2}', & \ft(\kappa)&=1.
\end{align*}
\end{enumerate}

\subsubsection{}\label{subsubsec:example2} Let $\Gamma=\Z_2\times \Z_2=\langle \ut\rangle \times \langle \kt \rangle$. Reordering generators, the non-trivial possibilities for $\Gamma_{u,v,\kappa}$ are
$\Gamma_{\ut,1,\kt}$ and $\Gamma_{\ut,\kt,\kt}$, which are groups of order $16$ with Gap Id $3,4$ and Hall-Senior number $\#_{16}9,\#_{16}10$, see \cite[Chapter 7]{Len12}.

\begin{enumerate}[leftmargin=*,label=\rm{\alph*)}]
\item\label{item:example2-a} We apply Proposition \ref{prop_cohomology_semidirect} to $\Gamma_{\ut,1,\kt}$ with $\Lambda=\langle\ut\rangle$, $\Omega=\langle \kt\rangle$, $P\equiv 1$, $\psi(\ut)=\chi(\kt)=-1$: we get a $2$-cocycle $\sigma$ such that 
\begin{align*}
\frac{\sigma(x,\kt)}{\sigma(\kt,x)} &=-1, &
\frac{\sigma(x,\ut)}{\sigma(\ut,x)} &= 1, &
\frac{\sigma(y,\kt)}{\sigma(\kt,y)} &= 1, &
\frac{\sigma(y,\ut)}{\sigma(\ut,y)} &=-1;
\end{align*}
for the non-trivial choice $P(\ut,\kt)=-1$ there is no suitable $\psi$.

For the decomposition $\Lambda=1$, $\Omega=\Gamma$, we have the non-trivial choice 
$P=1$, $\psi=1$, $\chi(\ut)=\chi(\kt)=-1$. We obtain a $2$-cocycle $\sigma'$ such that 
\begin{align*}
\frac{\sigma'(x,\kt)}{\sigma'(\kt,x)} &=-1, &
\frac{\sigma'(x,\ut)}{\sigma'(\ut,x)} &=-1, &
\frac{\sigma'(y,\kt)}{\sigma'(\kt,y)} &= 1, &
\frac{\sigma'(y,\ut)}{\sigma'(\ut,y)} &= 1.
\end{align*}

Accordingly it is known that $H^2( \#_{16}9,\Bbbk^\times)=\Z_2\times\Z_2$.

\smallskip

\item\label{item:example2-b} We apply Proposition \ref{prop_cohomology_extension} to $\Gamma_{\ut,1,\kt}$ and the $2$-cocycle $\sigma$ constructed above, where $\wt=\kt^s\ut^t$, $f(x)=(-1)^s$, $f(y)=(-1)^t$, $s,t\in\{0,1\}$. 

\smallskip

\item\label{item:example2-c} Assume that there exists a surjective map $\pi:\Gamma\to \Z_2\times\Z_2$ such that 
$\pi(u)=\ut$, $\pi(v)=1$, $\pi(\kappa)=\kt$. The pullback of $\sigma$ from $(\Z_2\times\Z_2)_{\ut,1,\kt}$ to $\Gamma_{u,v,\kappa}$ can be lifted to $\Gamma(\sqrt[r]{w})$ if either $\pi(w)=1$ or $2\nmid r$. If $\pi(w)=\kappa$, $2\mid r$, then the $2$-cocycle can be lifted to $\Gamma_{u,v,\kappa}$ if there exists $\ft\in\widehat{\Gamma}$ such that
$\ft(\kappa)=1$ and 
\begin{align*}
\ft(u) \, (\text{resp. }\ft(v))=\xi & \in \begin{cases}
\G_r, \, 2 \nmid r/\mathrm{ord}(\xi), & \text{ if }s=1 \, (\text{resp. }t=1), \\
\G_{r/2}', & \text{ if }s=0 \,  (\text{resp. }t=0).
\end{cases}
\end{align*}
\end{enumerate}

\subsubsection{}\label{subsubsec:example3} Fix $t\in\N$. Set $k=4t$, $\Gamma=\Z_k\times \Z_k$ with generators $\ut$, $\vt$, and consider $\kt:=\vt^{2t}$. The group $\Gamma_{\ut,\vt,\kt}$ has order $4k^2$, center $\Gamma$, and is presented by generators $x,y$ and relations $x^{2k}=y^{2k}=1$, $[x,y]=y^k$, where $\ut=x^2$, $\vt=y^2$, $\kt=y^{k}$. 

Moreover $\langle y \rangle\simeq \Z_{2k}$ is a normal subgroup, $\langle x \rangle\simeq \Z_{2k}$ and $\Gamma_{\ut,\vt,\kt} \simeq \Z_{2k}\ltimes \Z_{2k}$ with action $x \cdot y=y^{k+1}$ (notice that $x^2\cdot y=y$).

\begin{enumerate}[leftmargin=*,label=\rm{\alph*)}]
\item\label{item:example3-a}  We apply Proposition \ref{prop_cohomology_semidirect} to $\Gamma_{\ut,\vt,\kt}$ with $\Lambda=\langle \ut \rangle$, $\Omega=\langle \vt \rangle$. Fix $\xi\in\G_k$. 
\begin{itemize}[leftmargin=*]
\item The pairings $P:\Lambda\times\Omega$ such that $P(\ut,\kt)=1$ are given by
\begin{align*}
P(\ut,\vt)&=\xi^{2i} & \text{ for some }&i\in\I_{2t}.
\end{align*}
\item $\chi\in\widehat{\Lambda}$ satisfies $\chi(\ut)^2=P(\ut,\vt)$ if and only if $\chi(\ut)=p_1\xi^i$ for some $p_1\in\{\pm1\}$. Analogously, $\psi\in\widehat{\Omega}$ is given by $\psi(\vt)=p_2\xi^i$, $p_2\in\{\pm1\}$.
\item As $\chi(\kt)=(-1)^i$, the condition $\psi(\vt)\psi(\kt)=\chi(\ut)$ always holds when $i$ is even, and for $i$ odd there are two choices since we need $p_1p_2=-1$. 
\end{itemize}
Altogether, when $i$ is even we obtain $4t$ different 2-cocycles $\sigma_{ip_1p_2}$ with $\chi(\kt)=1$, 
and for $i$ odd we have $2t$ different 2-cocycles $\sigma_{ip_1p_2}$ with $\chi(\kappa)=-1$. Set $\sigma:=\sigma_{1+-}$. 

\smallskip

\item\label{item:example3-b} We apply Proposition \ref{prop_cohomology_extension} to $\Gamma_{\ut,\vt,\kt}$, $\sigma$ as constructed above and $\wt=\ut^s\vt^t$: 
here, $f(x)=\psi(\vt^t)=(-\xi)^t$ and $f(y)=\chi(\ut^s)=\xi^s$.

\smallskip

\item\label{item:example3-c} More generally, fix a surjective map $\pi:\Gamma\to \Z_k\times \Z_k$, $u,v,\kappa\in\Gamma$ such that $\pi(u)=\ut$, $\pi(v)=\vt$, $\pi(\kappa)=\kt$,
$w=u^sv^t$, $s,t\in\I_k$, and $r\in\N$. The pullback of $\sigma$ from $(\Z_k\times \Z_k)_{\ut,\vt,\kt}$ to $\Gamma_{u,v,\kappa}$ can be lifted to $\Gamma(\sqrt[r]{w})$ if there exists $\ft \in\widehat{\Gamma}$ such that $\ft(x)^r=(-\xi)^t$, $\ft(y)^r=\xi^s$ and $\ft(\kappa)=1$.
\end{enumerate} 

\subsection{Proof of Proposition \ref{prop:DoiTwist}}\label{subsubssection:proof-lemma-Doi-twist}

We proceed case-by-case. We use the representations $M(g,\chi)$ coming from the Yetter-Drinfeld structure for each group in order to get a triple as in Proposition \ref{prop_cohomology_semidirect}, which in turn gives a $2$-cocycle, then use Proposition \ref{prop_cohomology_extension} when we need to extend the group accordingly. For each $M(g,\chi)$ we choose a basis given by centralizer coset representatives and give the corresponding matrices: When the matrix is a multiple of the identity we just write the corresponding scalar.

\subsubsection{Type $\alpha_2$}\label{subsubsec:proof-lemma-A2} Here $M=M(g_1,\chi_1)\oplus M(g_2,\chi_2)$, see \S \ref{subsubsec:alpha2-presentation}. 
The four possibilities for the parity vector $\mathtt{P}=(\chi_1(\kappa),\chi_2(\kappa))$ fall into two orbits under the Weyl groupoid action, namely $(1,1)$ and $\{(-1,1), (-1,-1),(-1,1)\}$. As the first one corresponds to trivial action of $\kappa$, we just need to study $(-1,1)$. Set $\paraA:=\chi_2(g_1^2)=-\chi_1(g_2^2)^{-1}$, $k=\ord \paraA$.
We compute $M(g_1,\chi_1)$, $M(g_2,\chi_2)$, $M(g_1g_2,\chi_{1}\chi_2)$. 

\begin{center}
\begin{tabular}{c|ccc|ccc}
& $\kappa$ & $g_1$ & $g_2$ 
& $g_1^2=u$ & $g_2^2=v$ &  \\
\hline
$M(g_1,\chi_1)$
& $-1$
& $\begin{pmatrix}  -1 & 0 \\ 0 & -1 \end{pmatrix}$
& $\begin{pmatrix}  0 & -\paraA^{-1} \\ 1 & 0 \end{pmatrix}$
& $1$ & $-\paraA^{-1}$ & 
\\
$M(g_2,\chi_2)$
& $1$
& $\begin{pmatrix}  0 & \paraA \\ 1 & 0 \end{pmatrix}$
& $\begin{pmatrix}  -1 & 0 \\ 0 & -1 \end{pmatrix}$
& $\paraA$ & $1$ & 
\\
\hline
$M(g_1g_2,\chi_{1}\chi_2)$
& $-1$
& $\begin{pmatrix}  0 & \paraA \\ 1 & 0 \end{pmatrix}$
& $\begin{pmatrix}  0 & 1 \\ -\paraA^{-1} & 0 \end{pmatrix}$
& $\paraA$ & $-\paraA^{-1}$ & 
\end{tabular}
\end{center}
From the actions above we can read off the group $G^{\min}$ explicitly:
\begin{align*}
G^{\min} &= \left\langle g_1,g_2,\kappa \;\mid\; 
[g_1,g_2]=\kappa,\;
[g_1,\kappa]=[g_2,\kappa]=(g_1^2\kappa)^k=(g_2^2)^k=\kappa^2=1,\right.
\\
&\left.\hspace{3cm}\text{if $2\mid k$ we add } (g_1^2\kappa)^{k/2}=\kappa \right\rangle.
\end{align*}

By Proposition \ref{prop_2group} we reduce to $2$-groups, so we have three cases:
\begin{itemize}[leftmargin=*]
\item $k=1$, that is $\paraA=1$. Then $x^2=1, y^2=\kappa$ and the group $G^{\min}$ is the dihedral group of order $8$, see \S \ref{subsubsec:example1} \ref{item:example1-a}.

\item $k=2$, that is $\paraA=-1$. Here $x^2=u$ has order $2$, $y^2=v=1$, $(xy)^2=\kappa u$. Then $G^{\min}$ is the group $\#_{16}9$ and such 2-cocycle exists, see \S \ref{subsubsec:example2} \ref{item:example2-a}.
\item $k=2^{n}$, $n\geq 2$. In this case $x^{2^n}=y^{2^n}=1$ and $y^{2^{n-1}}=\kappa$. Then 
$|G^{\min}|=2^{2+2n}$, the center is $\Z_{2^n}\times \Z_{2^n}$, and moreover $G^{\min} \simeq \Z_{2^{n+1}}\ltimes \Z_{2^{n+1}}$ as in \ref{subsubsec:example3} \ref{item:example3-a}, so there exists such 2-cocycle.
\end{itemize}


\subsubsection{Type $\alpha_3$}\label{subsubsec:proof-lemma-A3}
Here $M=M(g_1,\chi_1)\oplus M(g_2,\chi_2) \oplus M(g_3,\chi_3)$, see \S \ref{subsubsec:alpha3-presentation}. 
As $\chi_1(\kappa)=\chi_3(\kappa)$, there are four choices of $\mathtt{P}=(\chi_i(\kappa))_{i\in\I_3}$, which fall into two orbits under the Weyl groupoid action, namely 
\begin{align*}
& \{(1,1,1)\} && \text{and} &&\{(-1,1,-1),(-1,-1,-1),(1,-1,1)\}.
\end{align*}
Thus we just need to study $\mathtt{P}=(-1,1,-1)$. We compute the representations $M(g_i,\chi_i)$, which contains the previous case $A_2$. Set 
\begin{align*}
\paraA &:=\chi_2(g_1^2)=-\chi_1(g_2^2)^{-1}, & \paraB&:=\chi_3(g_1)=\chi_1(g_3)^{-1}, & \paraC&=\chi_2(g_1g_3^{-1}).
\end{align*}
Hence the action of $G$ on $M$ is given by:
\begin{center}
\begin{tabular}{c|cccc}
& $\kappa$ & $g_1$ & $g_2$ & $g_3$ 
\\
\hline
$M(g_1,\chi_1)$
& $-1$
& 
$\begin{pmatrix}  -1 & 0 \\ 0 & 1 \end{pmatrix}$
& 
$\begin{pmatrix}  0 & -\paraA^{-1} \\ 1 & 0 \end{pmatrix}$
& 
$\begin{pmatrix}  \paraB^{-1} & 0 \\ 0 & -\paraB^{-1} \end{pmatrix}$
\\
$M(g_2,\chi_2)$
& $1$
& 
$\begin{pmatrix}  0 & \paraA \\ 1 & 0 \end{pmatrix}$
& 
$\begin{pmatrix}  -1 & 0 \\ 0 & -1 \end{pmatrix}$
& 
$\begin{pmatrix}  0 & \paraC\paraA \\ \paraC & 0 \end{pmatrix}$ 
\\
$M(g_3,\chi_3)$
& $-1$
& 
$\begin{pmatrix}  \paraB & 0 \\ 0 & -\paraB \end{pmatrix}$
& 
$\begin{pmatrix}  0 & -(\paraC^2\paraA)^{-1} \\ 1 & 0 \end{pmatrix}$
& 
$\begin{pmatrix}  -1 & 0 \\ 0 & 1 \end{pmatrix}$
\end{tabular}
\end{center}

Using that $g_1g_2=\kappa g_2g_1$, $g_3g_2=\kappa g_2g_3$, $g_1g_3=g_3g_1$, $\kappa$ is central and $\kappa^2=1$, cf. \S \ref{subsubsec:alpha3-presentation}, the subgroup $\widetilde{\Gamma}$ generated by \footnote{We set the generators of $\widetilde{\Gamma}$ according with generators for the symplectic root system $n=3$, $r=1$ in \cite[Thm. 4.5]{Len14}.}
\begin{align*}
u &:= g_1^2, & v &:= g_2^2\kappa, & \formalw & := g_3g_1^{-1}\kappa, & &\kappa,
\end{align*}
is contained in $Z(G^{\min})$, and $G^{\min}/\widetilde{\Gamma}$ has four elements: $\widetilde{\Gamma}$, $g_1\widetilde{\Gamma}$, $g_2\widetilde{\Gamma}$, $g_1g_2\widetilde{\Gamma}$. From here we check that $G^{\min}\simeq \widetilde{\Gamma}_{u,v,\kappa}$, where $x\mapsto g_1$, $y\mapsto g_2$. 

Let $\Gamma$ be the subgroup generated by $u$, $v$, $\kappa$. Set $N_1=\operatorname{lcm}(\ord \paraA,\ord \paraB^2)$, $N_2=\operatorname{lcm}(\ord \paraA,\ord \paraC^2)$, $N_3=\operatorname{lcm}(\ord \paraB,\ord \paraC)$. 
As $\formalw^{N_3}=\id$, we can define
\begin{align*}
r & := \min \{s\in\I_{N_3} | \formalw^s \in \Gamma \}
\\
& = \min \{s\in\I_{N_3} | \exists m\in\I_{N_1},n\in\I_{N_2}: \paraC^{2n}=\paraB^{2m+2r}, \paraA^{n}\paraB^{r}=1, \paraA^{m}=\paraC^{r}  \}.
\end{align*}
We have that $\widetilde{\Gamma}\simeq \Gamma(\sqrt[r]{w})$ for $w=\formalw^r$.
The action of $u$, $v$, $\formalw$ on $M_1$, $M_2$, $M_3$ is given, respectively, by the following scalars: 
\begin{align*}
&(1,\paraA,\paraB^{-2}), &&(\paraA^{-1},1,\paraA^{-1}\paraC^{-2}), &&(\paraB^{-1},\paraC,\paraB^{-1}). 
\end{align*}
We also set $k_1=\ord \paraA$, $k_2=\ord \paraB$, $k_3=\ord \paraC$. Then
\begin{align*}
\Gamma & \simeq \big\langle \kappa,u,v \mid 
\kappa^2=1,u^{N_1}=1, v^{N_2}=1, \kappa=(u^i v)^{k_1/2}\text{ if } 2\mid k_1, \paraB^{2i}=\paraC^{-k_1}
\big \rangle.
\end{align*}
As $k_1|N_1,N_2$, there exists a surjective map $\Gamma\twoheadrightarrow \Gamma'$, where $\Gamma'$ is an abelian group as in \S \ref{subsubsec:proof-lemma-A2}, so there exists a $2$-cocycle $\sigma'$ for $\Gamma'$ such that $\tfrac{\sigma'(g_2,\kappa)}{\sigma'(\kappa,g_2)}=1$, $\tfrac{\sigma'(g_1,\kappa)}{\sigma'(\kappa,g_1)}=-1=\tfrac{\sigma'(g_3,\kappa)}{\sigma'(\kappa,g_3)}$. Let $\sigma$ be the pullback of $\sigma'$ on $\Gamma$: we look for a lift on $G^{\min} \simeq \Gamma(\sqrt[r]{w})_{u,v,\kappa}$, so we look for a character $\ft$ as in Proposition \ref{prop_cohomology_extension}.

By Proposition \ref{prop_2group} it is enough to solve the case in which the three $k_i$ are powers of $2$. We split in three cases as in \S \ref{subsubsec:proof-lemma-A2}:
\begin{itemize}[leftmargin=*]
\item $k=1$, that is $\paraA=1$. 
Either $w=1$ or else $k_2>k_3$, $w=\kappa$ and $r=k_2/2>1$. In the second case we construct $\ft\in\widehat{\Gamma}$ in $M(g_1,\chi_1)^*\otimes M(g_3,\chi_1)$: that is, $\ft(u)=\paraB^2 \in\G_r'$, $\ft(v)=\paraC^{-2}\in\G_{r/2}'$, $\ft(\kappa)=1$. Then we apply \S \ref{subsubsec:example1} \ref{item:example1-c}.

\item $k=2$, that is $\paraA=-1$. If $k_2>k_3$, then $r=k_2/2$, $w=\kappa$; if $k_2<k_3$, then $r=k_3/2$, $w=u$; otherwise $k_2=k_3$, $r=k_2/2$ and $w=u\kappa$. In any case we construct $\ft\in\widehat{\Gamma}$ in $M(g_1,\chi_1)^*\otimes M(g_3,\chi_1)$ as in \S \ref{subsubsec:example2} \ref{item:example2-c}, and there exists such a lift.

\item $k=2^{n}$, $n\geq 2$. Here $r=\max(k_2/k_1,\; k_3/k_1,\;1)$, and we construct $\ft\in\widehat{\Gamma}$, again in $M(g_1,\chi_1)^*\otimes M(g_3,\chi_1)$:
\begin{align*}
\ft(\kappa)&=1, & \ft(u)&=\paraB^2, & \ft(v)&=\paraC^2.
\end{align*}
This $\ft$ satisfies the conditions in \S \ref{subsubsec:example3} \ref{item:example3-c}, and there exists such a lift.
\end{itemize}

\subsubsection{Type $\delta_4$}\label{subsubsec:proof-lemma-D4} 
As $\chi_1(\kappa)=\chi_3(\kappa)=\chi_4(\kappa)$, the choices of $\mathtt{P}=(\chi_i(\kappa))_{i\in\I_4}$ fall into two orbits under the Weyl groupoid action, namely 
\begin{align*}
& \{(1,1,1,1)\} && \text{and} &&\{(-1,1,-1,-1),(-1,-1,-1,-1),(1,-1,1,1)\},
\end{align*}
where the second entry denotes the center node in the Dynkin diagram. Now we study $\mathtt{P}=(-1,1,-1,-1)$. We fix the central elements 
\footnote{According with generators for the symplectic root system $n=4$, $r=2$ in \cite{Len14}.} $z=g_3g_1^{-1}$, $z'=g_4g_1^{-1}$.
Then this case can be achieved by combining the previous result of extending $\alpha_2$ to $\alpha_3$ by $z$ and by $z'$, see Corollary \ref{cor_cohomology_doubleextension}. 

\subsubsection{Type $\gamma_3$}\label{subsubsec:proof-lemma-C3} 
Here $M=M(g_1,\chi_1)\oplus M(g_2,\chi_2)\oplus M(g_3,\chi_3)$, with $g_3\in Z(G)$, $\chi_3\in\widehat{G}$. As $\kappa=[g_1,g_2]$, we have that $\chi_3(\kappa)=1$. 
The possible $\mathtt{P}=(\chi_i(\kappa))_{i\in\I_3}$ fall into two orbits under the Weyl groupoid action, namely 
$\{(1,1,1)\}$ and $\{(-1,1,1),(-1,-1,1),(1,-1,1)\}$. Fix $\mathtt{P}=(-1,1,1)$ and set
\begin{align*}
\paraA &:=\chi_2(g_1^2)=-\chi_1(g_2^{-2}), \, \paraB :=\chi_3(g_1)=\chi_1(g_3^{-1}), \, \paraC:=\chi_3(g_2)=-\chi_2(g_3^{-1}). 
\end{align*}
Now we compute the representations $M(g_i,\chi_i)$: 
\begin{center}
\begin{tabular}{c|cccc}
& $\kappa$ & $g_1$ & $g_2$ & $g_3$ 
\\
\hline
$M(g_1,\chi_1)$
& $-1$
& 
$\begin{pmatrix}  -1 & 0 \\ 0 & 1 \end{pmatrix}$
& 
$\begin{pmatrix}  0 & -\paraA^{-1}\\ 1 & 0 \end{pmatrix}$
& $\paraB^{-1}$
\\
$M(g_2,\chi_2)$
& $+1$
& 
$\begin{pmatrix}  0 & \paraA \\ 1 & 0 \end{pmatrix}$
& 
$\begin{pmatrix}  -1 & 0 \\ 0 & -1 \end{pmatrix}$
& $-\paraC^{-1}$
\\
$M(g_3,\chi_3)$
& $+1$
& $\paraB$
& $\paraC$ 
& $-1$ 
\end{tabular}
\end{center} 
If $u=g_1^2$, $v=g_2^2$, $\formalw=g_3$, then $G^{\min}\simeq \Gamma(\sqrt[r]{w})_{u,v,\kappa}$ for $\Gamma=\langle u,v,\kappa\rangle$ and appropriate $w\in\Gamma$, $r\in\N$. Moreover $G^{\min}$ is isomorphic to the one in \S \ref{subsubsec:proof-lemma-A3}, hence there exists a $2$-cocycle as in Proposition \ref{prop:DoiTwist}. 

\subsubsection{Type $\gamma_4$}\label{subsubsec:proof-lemma-C4} 
Here we have central elements $z=g_3g_1^{-1}$, $z'=g_4$, and this case is solved by combining the previous result of extending $A_2$ to $A_3$ by $z$ and by $z'$, see Corollary \ref{cor_cohomology_doubleextension}. 

\subsubsection{Type $\phi_4$}\label{subsubsec:proof-lemma-F4} 
Here $M=M(g_1,\chi_1)\oplus M(g_2,\chi_2)\oplus M(g_3,\chi_3)\oplus M(g_4,\chi_4)$, with $g_3,g_4\in Z(G)$, $\chi_3,\chi_4\in\widehat{G}$, and $\chi_3(\kappa)=\chi_4(\kappa)=1$. 
The possible $\mathtt{P}=(\chi_i(\kappa))_{i\in\I_4}$ fall into two orbits under the Weyl groupoid action: 
$\{(1,1,1,1)\}$ and $\{(-1,1,1,1),(-1,-1,1,1),(1,-1,1,1)\}$. Fix $\mathtt{P}=(-1,1,1,1)$ and set
{\small
\begin{align*}
\paraA &:=\chi_2(g_1^2)=-\chi_1(g_2^{-2}), & \paraB &:=\chi_3(g_1)=\chi_1(g_3^{-1}), & \paraC&:=\chi_3(g_2)=-\chi_2(g_3^{-1}),
\\
\paraB'&:=\chi_4(g_1)=\chi_1(g_4^{-1}), & \paraC' &:=\chi_4(g_2)=\chi_2(g_4^{-1}), & \paraD&:=\chi_4(g_3)=-\chi_3(g_2^{-1}).
\end{align*}
}
Now we compute the representations $M_i:=M(g_i,\chi_i)$: 
\begin{center}
\begin{tabular}{c|ccccc}
& $\kappa$ & $g_1$ & $g_2$ & $g_3$ & $g_4$ 
\\
\hline
$M(g_1,\chi_1)$
& $-1$
& 
$\begin{pmatrix}  -1 & 0 \\ 0 & -1 \end{pmatrix}$
& 
$\begin{pmatrix}  0 & -\paraA^{-1}\\ 1 & 0 \end{pmatrix}$
& $\paraB^{-1}$
& $\paraB'^{-1}$
\\
$M(g_2,\chi_2)$
& $1$
& 
$\begin{pmatrix}  0 & \paraA \\ 1 & 0 \end{pmatrix}$
& 
$\begin{pmatrix}  -1 & 0 \\ 0 & -1 \end{pmatrix}$
& $-\paraC^{-1}$
& $\paraC'^{-1}$
\\
$M(g_3,\chi_3)$
& $1$
& $\paraB$
& $\paraC$ 
& $-1$
& $-\paraD^{-1}$
\\
$M(g_4,\chi_4)$
& $1$
& $\paraB'$
& $\paraC'$ 
& $\paraD$
& $-1$
\end{tabular}
\end{center}

Set $u=g_1^2$, $v=g_2^2$. We will construct a $2$-cocycle $\sigma$ on $G^{\min}$ by using appropiate $2$-cocycles from the previous cases:
\begin{itemize}[leftmargin=*]
\item Set $G_{12}=\langle g_1,g_2 \rangle$, $\Gamma=\langle \kappa,u,v\rangle$; let $\widetilde{G}_{12}, \Gamma_{12}\subset \End(M_1\oplus M_2)$ 
be the subgroups obtained by restriction. Then $\Gamma$, $\Gamma_{12}$ are central subgroups, $G_{12}=\Gamma_{u,v,\kappa}$, $\widetilde{G}_{12}=(\Gamma_{12})_{u,v,\kappa}$, with canonical projections $G_{12}\twoheadrightarrow \widetilde{G}_{12}$, $\Gamma\twoheadrightarrow \Gamma_{12}$, and $\widetilde{G}_{12}$, $\Gamma_{12}$ are as in \S \ref{subsubsec:proof-lemma-A2}. Hence there exists a $2$-cocycle as we need: the pullback $\sigma_{12}$ on $G_{12}$ satisfies
$\tfrac{\sigma_{12}(g_1,\kappa)}{\sigma_{12}(\kappa,g_1)}=-1$, $\tfrac{\sigma_{12}(g_2,\kappa)}{\sigma_{12}(\kappa,g_2)}=1$.

\item For $j=3,4$ we set $G_{12j}=\langle g_1,g_2,g_j \rangle$, $\Gamma_{12j}=\langle \kappa,u,v,g_j\rangle$. Then $\Gamma_{12j}$ is a central subgroup of the form
$\Gamma_{12j}\simeq \Gamma(\sqrt[r_j]{w_j})$ for appropiate $r_j\in\N$, $w_j\in\Gamma$: the proof for $j=3$ is the same as in \S \ref{subsubsec:proof-lemma-C3} since $M_1\oplus M_2\oplus M_3$ is of type $C_3$, and for $j=4$ we have the same structure (luckily). Using the same argument as in \S \ref{subsubsec:proof-lemma-C3} we check the existence of $2$-cocycle $\sigma_{12j}$ on $G_{12j}$ such that $\tfrac{\sigma_{12j}(g_1,\kappa)}{\sigma_{12j}(\kappa,g_1)}=-1$, $\tfrac{\sigma_{12j}(g_2,\kappa)}{\sigma_{12j}(\kappa,g_2)}=\tfrac{\sigma_{12j}(g_j,\kappa)}{\sigma_{12j}(\kappa,g_j)}=1$.

\item Finally, $G^{\min}\simeq \Gamma(\sqrt[r_3]{w_3})(\sqrt[r_4]{w_4})_{u,v,\kappa} \simeq \Gamma(\sqrt[r_4]{w_4})\Gamma(\sqrt[r_3]{w_3})_{u,v,\kappa}$. The existence of a $2$-cocycle $\sigma$ on $G^{\min}$ such that
\begin{align*}
\tfrac{\sigma(g_1,\kappa)}{\sigma(\kappa,g_1)}&=-1, & \tfrac{\sigma(g_2,\kappa)}{\sigma(\kappa,g_2)}= \tfrac{\sigma(g_3,\kappa)}{\sigma(\kappa,g_3)}=\tfrac{\sigma(g_4,\kappa)}{\sigma(\kappa,g_4)}=1
\end{align*}
follows from the $2$-cocycles $\sigma_{12j}$ on $G_{12j}$, $j=3,4$ and Corollary \ref{cor_cohomology_doubleextension}.
\end{itemize}

This concludes the proof of Proposition \ref{prop:DoiTwist}.
\qed

\subsection*{Acknowledgements} The authors thank Istvan Heckenberger for many useful discussions, suggestions and hospitality. In particular, for the discussion of Lemma \ref{lem:trivial-action-kappa-big-rank}.


\begin{thebibliography}{AAAA}

\bibitem[A]{A-leyva} N. Andruskiewitsch. \emph{An Introduction to Nichols Algebras}. In Quantization, Geometry and Noncommutative Structures in Mathematics and Physics. 
A. Cardona, P. Morales, H. Ocampo, S. Paycha, A. Reyes, eds., pp. 135--195, Springer (2017).

\bibitem[AA]{AA-diag-survey} N. Andruskiewitsch, I. Angiono. 
\emph{On finite-dimensional Nichols algebras of diagonal type.} 
Bull. Math. Sci. \textbf{7} (2017), 353--573. 


\bibitem[A+]{AAGMV} N. Andruskiewitsch, I. Angiono, A. Garc\'ia Iglesias, A. Masuoka,
C. Vay. \emph{Lifting via cocycle deformation}. J. Pure Appl. Alg. {\bf 218}  (2014), 684--703.

\bibitem[AAG]{AAG} N. Andruskiewitsch, I. Angiono, A. Garc\'ia Iglesias.
\emph{Liftings of Nichols algebras of diagonal type I. Cartan type A}. Int. Math. Res. Not. IMRN {\bf 2017} (9) (2017), 2793--2884.

\bibitem[ACG]{ACG} N. Andruskiewitsch, G. Carnovale, G. Garc\'ia.
\emph{Finite-dimensional pointed Hopf algebras over finite simple groups of Lie type V. Mixed classes in Chevalley and Steinberg groups}.  
Manuscripta Math., to appear.

\bibitem[AG]{AG} N. Andruskiewitsch, M. Gra\~na,
\emph{From racks to pointed Hopf algebras}. Adv. Math. \textbf{178} (2003), 177--243.

\bibitem[AS1]{AS-adv} N. Andruskiewitsch, H.-J. Schneider. 
\emph{Finite quantum groups and Cartan matrices}, Adv. Math. \textbf{154} (2000), 1--45.

\bibitem[AS2]{AS-cambr} \bysame. 
\emph{Pointed Hopf algebras}. In Recent developments in Hopf algebras Theory, 
MSRI Publ. \textbf{43} (2002), 1--68, Cambridge Univ. Pr.

\bibitem[AS3]{AS4} \bysame.
\emph{On the classification of finite-dimensional pointed Hopf algebras}, 
Ann. of Math. \textbf{171} (2010), 375--417.

\bibitem[An1]{A-standard} I. Angiono, \emph{On Nichols algebras with standard braiding},
Algebra Number Theory \textbf{1}  (2009), 35-106.

\bibitem[An2]{A-presentation} \bysame. \emph{On Nichols algebras of diagonal type}. 
J. Reine Angew. Math.  \textbf{683} (2013),  189--251.

\bibitem[An3]{A-jems} \bysame. \emph{A presentation by generators and relations of Nichols
algebras of diagonal type and convex orders on root systems}. 
J. Eur. Math. Soc.  \textbf{17} (2015), 2643--2671.

\bibitem[An4]{A-pre-Nichols}  \bysame. \emph{Distinguished pre-Nichols algebras}. Transf. Groups \textbf{21} (2016), 1--33.

\bibitem[ACS]{ACS} I. Angiono, E. Campagnolo and G. Sanmarco.
\emph{Finite GK-dimensional pre-Nichols algebras of super and standard type}, \texttt{arXiv:2009.04863}.

\bibitem[AnG]{AnG} I. Angiono and A. Garc\'\i a Iglesias.
\emph{Liftings of Nichols algebras of diagonal type II. All liftings are cocycle deformations}. 
Selecta Math. \textbf{25} (2019), Paper No. 5, 95 pp.

\bibitem[AnGP]{AGP} I. Angiono, C.  Galindo, M. Pereira. \emph{De-equivariantization of Hopf algebras}. Algebr. Represent. Theory \textbf{17} (2014), 161--180. 

\bibitem[AnKM]{AKM} I. Angiono, M. Kochetov, M. Mastnak.
\emph{On rigidity of Nichols algebras}.
J. Pure Appl. Algebra \textbf{219} (2015), 5539--5559. 

\bibitem[AnS]{AnS} I. Angiono, G. Sanmarco.
\emph{Pointed Hopf algebras over non-abelian groups with decomposable braidings. I.} 
J. Algebra \textbf{549} (2020), 78-111.

\bibitem[AJS]{AJS} H.H. Andersen, J.C. Jantzen, W. Soergel, \emph{Representations of quantum groups at a pth root of unity and of semisimple groups in characteristic
p: independence of p}. Astérisque \textbf{220} (1994), 321.

\bibitem[CH]{CH} M. Cuntz, I. Heckenberger, \emph{Finite Weyl groupoids of rank three}. Trans. Amer. Math. Soc. \textbf{364} (2012), 1369--1393. 

\bibitem[DoT]{DT} Y. Doi, M. Takeuchi, \emph{Multiplication alteration by two-cocycles - the quantum version},
Comm. Algebra {\bf 22} (1994), 5715--5732.

\bibitem[Dr]{D-qg} V. Drinfeld, \emph{Quantum groups}, Proc. Int. Congr. Math., Berkeley 1986, Vol. \textbf{1} (1987), 798--820.

\bibitem[EGNO]{EGNO} P. Etingof, S. Gelaki, D. Nikshych, V. Ostrik, 
\emph{Tensor categories}. Mathematical Surveys and Monographs \textbf{205},
AMS, Providence, RI, 2015, xvi+343pp.

\bibitem[ENOM]{ENOM} P. Etingof, D. Nikshych, V. Ostrik, E. Meir, \emph{Fusion categories and homotopy theory}. Quantum Topol. \textbf{1} (2010), 209--273.

\bibitem[FBZ]{FBZ} E. Frenkel, D. Ben-Zvi, 
\emph{Vertex Algebras and Algebraic Curves}, 
Mathematical Surveys and Monographs 88, American Mathematical Society (2004).

\bibitem[FGST]{FGST} B. Feigin, A. Gainutdinov, A. Semikhatov, I. Tipunin,
\emph{Kazhdan-Lusztig correspondence for the representation category of the triplet W-algebra in logarithmic CFT}, Theor. Math. Phys. 148 (2006).


\bibitem[FK]{FK} S. Fomin, K. N. Kirillov, \emph{Quadratic algebras, Dunkl elements, and Schubert calculus}.
Progr. Math. {\bf 172} (1999), 146--182.

\bibitem[Fu]{Fu} J. Fuchs, \emph{Affine Lie algebras and quantum groups}, Cambridge Monographs Math. Phys. (1994).

\bibitem[Gr]{Gr00} M. Gra\~{n}a, \emph{On Nichols algebras of low dimension}, Contemp. Math. \textbf{267} (2000).

\bibitem[Gu]{Gu} R. G\"unther, \emph{Crossed products for pointed Hopf algebras}.
Comm.  Algebra, \textbf{27} (1999), 4389--4410.


\bibitem[H1]{H-inv} I. Heckenberger. 
\emph{The Weyl groupoid of a Nichols algebra of diagonal type}. 
Invent. Math. \textbf{164} (2006), 175--188 .

\bibitem[H2]{H-classif RS} \bysame 
\emph{Classification of arithmetic root systems}. 
Adv. Math. \textbf{220} (2009), 59--124 .

\bibitem[HLV]{HLV} I. Heckenberger, A. Lochmann, L. Vendramin, \emph{Braided racks, Hurwitz actions and Nicholsalgebras with many cubic relations}, Transform. Groups \textbf{17} (2012), 157--194.

\bibitem[HS1]{HS-rank2-1} I. Heckenberger, H.-J. Schneider.
\emph{Nichols algebras over groups with finite root system of rank two I}. 
J. Algebra \textbf{324} (2010), 3090--3114.

\bibitem[HS2]{HS-right-coideal} I. Heckenberger, H.-J. Schneider.
\emph{Right coideal subalgebras of Nichols algebras and the Duflo order on the Weyl groupoid}.
Isr. J. Math. \textbf{197} (2013), 139--187.

\bibitem[HS3]{HS-book} \bysame \emph{Hopf algebras and root systems}. Mathematical Surveys and Monographs \textbf{247}. 
Providence, RI: American Mathematical Society (AMS) (ISBN 978-1-4704-5232-2/hbk; 978-1-4704-5680-1/ebook). xix, 582 p. (2020). 

\bibitem[HV]{HV-rank>2} I. Heckenberger, L. Vendramin.
\emph{A classification of Nichols algebras of semi-simple Yetter-Drinfeld modules over non-abelian groups}. 
J. Eur. Math. Soc. \textbf{19} (2017), 299--356. 

\bibitem[HeY]{HY} I. Heckenberger, H. Yamane, \emph{A generalization of Coxeter groups, root systems, and Matsumoto's theorem}, Math. Z.
\textbf{259} (2008), 255--276.


\bibitem[IM]{IM64} N.  Iwahori, H.  Matsumoto,  \emph{Several  remarks  on  projective  representations  of finite  groups}. J. Faculty of Science Univ. Tokyo  \textbf{10} (1964), 129--146.

\bibitem[J]{Ji}  M. Jimbo, \emph{A $q$-difference analogue of $U(g)$ and the Yang Baxter equation}, Lett. Math. Phys. \textbf{10} (1985), 63--69.

\bibitem[KL]{KL} T. Kerler, V.V. Lyubashenko, \emph{Non-Semisimple Topological Quantum Field Theories for 3-Manifolds with Corners}, Springer Lect. Notes Math., vol. 1765, Springer-Verlag, New York (2001).

\bibitem[L1]{Len12} S. Lentner, \emph{Orbifoldizing Hopf and Nichols Algebras}. Dissertation, LMU München (2012), \url{https://edoc.ub.uni-muenchen.de/15363/}. 

\bibitem[L2]{Len14} S. Lentner, \emph{New large-rank Nichols algebras over non-abelian groups with commutator subgroup $\Z_2$}, J. Algebra \textbf{419} (2014), 1--33.

\bibitem[L3]{Len17} S.D. Lentner, \emph{Quantum groups and Nichols algebras acting on conformal field theories}, Adv. Math., to appear.

\bibitem[Lu1]{L - mod rep} G. Lusztig. \emph{Modular representations and quantum groups}. 
Contemp. Math. \textbf{82} 59--77, Amer. Math. Soc., Providence, RI, 1989.


\bibitem[Lu2]{Lu} G. Lusztig, \emph{Introduction to quantum groups}. Birkh\"auser (1993).

\bibitem[Ma]{Masuoka} A. Masuoka, \emph{Abelian and non-abelian second cohomologies of quantized
enveloping algebras}, J.  Algebra \textbf{320} (2008), 1--47.

\bibitem[Me]{Meir} E. Meir, \emph{Geometric perspective on Nichols algebras}. J. Algebra \textbf{601} (2022), 390--422.

\bibitem[MiS]{MS} A. Milinski, H.-J. Schneider, \emph{Pointed Indecomposable Hopf Algebras over Coxeter Groups}.
Contemp. Math. {\bf267} (2000), 215--236.

\bibitem[Mm]{Mombelli} M. Mombelli,
\emph{The Brauer-Picard group of the representation category of finite supergroup algebras}.
Rev. Uni\'on Mat. Argent. \textbf{55} (2014), 83--117. 

\bibitem[Mo]{Mo-libro} S. Montgomery.
\emph{Hopf algebras and their actions on rings}, CMBS \textbf{82},  
Amer. Math. Soc. (1993).

\bibitem[T]{Tah72} K. Tahara, \emph{On the Second Cohomology Groups of Semidirect Products}, Math. Z. \textbf{129} (1972), 365--379.

\bibitem[Tur]{Tur} V, Turaev, \emph{Quantum invariants of knots and 3-manifolds}, de Gruyter Studies in Mathematics 18 Walter de Gruyter \& Co (1994).

\end{thebibliography}
\end{document}